\documentstyle{amsppt}
\magnification=\magstep1
\baselineskip=14pt
\parindent=0pt
\parskip=14pt

\vsize=7.7in
\voffset=-.4in
\hsize=5.5in


\overfullrule=0pt

\define\A{{\Bbb A}}
\define\C{{\Bbb C}}

\define\R{{\Bbb R}}
\define\Q{{\Bbb Q}}

\define\Z{{\Bbb Z}}


\redefine\H{\frak H}


\define\M{\Cal M}    

\define\p{\frak p}

\define\a{\alpha}
\redefine\b{\beta}

\define\e{\epsilon}
\redefine\l{\lambda}
\redefine\o{\omega}
\define\ph{\varphi}

\define\s{\sigma}
\redefine\P{\Phi}
\predefine\Sec{\S}
\redefine\L{\Lambda}
\redefine\S{{\Cal S}}


\define\nat{\natural}


\define\back{\backslash}

\define\lra{\longrightarrow}
\redefine\tt{\otimes}
\define\scr{\scriptstyle}
\define\liminv#1{\underset{\underset{#1}\to\leftarrow}\to\lim}
\define\limdir#1{\underset{\underset{#1}\to\rightarrow}\to\lim}

\define\isoarrow{\ {\overset{\sim}\to{\longrightarrow}}\ }
\define\lisoarrow{\ {\overset{\sim}\to{\longleftarrow}}\ }


\define\nass{\noalign{\smallskip}}


\define\CH{\widehat{CH}}

\define\degH{\widehat{\text{\rm deg}}}

\redefine\ord{\text{\rm ord}}
\define\Ei{\text{\rm Ei}}
\redefine\O{\Omega}
\predefine\oldvol{\vol}
\redefine\vol{\text{\rm vol}}
\define\pr{\text{\rm pr}}



\define\Spec{\text{\rm Spec}\,}

\define\Sym{\text{\rm Sym}}
\define\tr{\text{\rm tr}}

\font\cute=cmitt10 at 12pt
\font\smallcute=cmitt10 at 9pt
\define\kay{{\text{\cute k}}}
\define\smallkay{{\text{\smallcute k}}}

%

%

\define\OB{\Cal O_B}
\define\End{\text{\rm End}}
\define\diag{\text{\rm diag}}
\redefine\Im{\text{\rm Im}}
\redefine\Re{\text{\rm Re}}
\define\SL{\text{\rm SL}}
\define\sgn{\text{\rm sgn}}
\define\und#1{\underline{#1}}


\define\GSpin{\text{\rm GSpin}}

\define\CT#1{\operatornamewithlimits{CT}_{#1}}

\define\E{\Bbb E}
\define\F{\Bbb F}
\define\Spf{\text{\rm Spf}\,}

\define\Lie{\text{\rm Lie}}
\define\Hom{\text{\rm Hom}}

\define\hor{{\text{\rm horiz}}}

\define\hfal{h_{\text{\rm Fal}}}

\define\Pic{\text{\rm Pic}}
\define\Aut{\text{\rm Aut}}
\define\Pich{\widehat{\Pic}}
\define\degh{\widehat{\deg}\ }

\define\tent#1{ \vphantom{\vbox to #1pt{}} }   

\define\ran{\,\operatorname{\rangle}}
\define\lan{\operatorname{\langle}\,}

\define\undd#1{\und{\und{#1}}}

\define\cha{\text{\rm char}}



\define\BorevichShafarevic{\bf1}
\define\bostbourb{\bf2}
\define\bost{\bf3}
\define\bostumd{\bf4}
\define\bostgilletsoule{\bf5}
\define\boutotcarayol{\bf6}
\define\chaifaltings{\bf7}
\define\cohen{\bf8}
\define\colmez{\bf9}
\define\delignemumford{\bf10}
\define\DR{\bf11}
\define\drinfeld{\bf12}
\define\Eichler{\bf13}
\define\faltings{\bf14}
\define\funkethesis{\bf15}
\define\funkecompo{\bf16}
\define\gelbart{\bf17}
\define\gsihes{\bf18}
\define\grossqc{\bf19}
\define\Hain{\bf20}
\define\hirzebruchzagier{\bf21}
\define\kottwitz{\bf22}
\define\splitting{\bf23}
\define\duke{\bf24}
\define\annals{\bf25}
\define\kbourb{\bf26}
\define\Bints{\bf27}
\define\krcrelle{\bf28}
\define\krinvent{\bf29}
\define\krsiegel{\bf30}
\define\tiny{\bf31}   
\define\kryIII{\bf32}   
\define\ky{\bf33}  
\define\kuehn{\bf34}
\define\lebedev{\bf35}
\define\naktag{\bf36}
\define\rao{\bf37}
\define\RZ{\bf38}
\define\raynaud{\bf39}
\define\shimura{\bf40}
\define\vistoli{\bf41}
\define\waldspurger{\bf42}
\define\weilI{\bf43}
\define\yangden{\bf44}
\define\zagier{\bf45}

\centerline{\bf Derivatives of Eisenstein series and Faltings heights
}

\centerline{by}

\centerline{Stephen S. Kudla\footnote{Partially
supported by NSF grant DMS-9970506 and by a Max-Planck Research Prize
from the Max-Planck Society and Alexander von Humboldt Stiftung. }}

\centerline{Michael Rapoport}

\centerline{Tonghai Yang\footnote{Partially
supported by NSF grant DMS-0070476}}

In a series of papers, \cite{\annals}, \cite{\krsiegel}, \cite{\krcrelle},
\cite{\krinvent}, \cite{\tiny}, \cite{\kbourb}, we showed that certain quantities from the
arithmetic geometry of Shimura varieties associated to orthogonal groups
occur in the Fourier coefficients of the derivative of suitable
Siegel-Eisenstein series. It was essential in these examples that this
derivative was the second term in the Laurent expansion of a
Siegel-Eisenstein series at the center of symmetry, and that the first
term in this Laurent expansion vanished ({\it incoherent case}). In the
present paper we prove a relation between a generating function for the
heights of Heegner cycles on the arithmetic surface associated to a
Shimura curve and the second term in the Laurent expansion at $s={1\over
2}$ of an Eisenstein series of weight ${3\over 2}$ for $SL_2$. It is
remarkable that $s={1\over 2}$ is not the center of symmetry and that the
first term of the Laurent expansion is non-zero. In fact, this nonzero
value has a geometric interpretation in terms of the Shimura curve over
the field of complex numbers. Considering the fact that the Eisenstein
series is a rather familiar classical object, it is surprising that this
interpretation of its Laurent expansion at $s={1\over 2}$ has not been
noticed before. As we will argue below in this introduction, we believe
that our result is part of a general pattern involving the heights of
divisors on arithmetic models of Shimura varieties associated to
orthogonal groups.

We now describe our results in more detail.

Let $B$ be an indefinite division quaternion algebra over $\Q$ and let $O_B$ be
a maximal order
in $B$. Let $D(B)$ be the product of all primes $p$ at which
$B$ is division. Let $\M$ be the moduli space of abelian varieties of dimension
$2$
with a (special) action of $O_B$. Then $\M$ is an integral model of the
Shimura curve attached to $B$; it is proper of relative dimension $1$ over
$\Spec(\Z)$,
with semi-stable reduction at all primes and is smooth at all primes $p$ at
which
$B$ splits, i.e., for $p\nmid D(B)$. Ignoring, for the moment, the fact that
$\M$ is only a stack,
we may consider $\M$ as an arithmetic surface in the sense of Arakelov theory,
\cite{\faltings}, \cite{\bost}, $\dots$.

For each $m\in \Z$ and for $v\in \R^\times_+$, we define a class in the
arithmetic
Chow group
$$\hat\Cal Z(m,v) = \big(\,\Cal Z(m), \Xi(m,v)\,\big)\in \CH^1(\M).\tag0.1$$
Here, for $m>0$, $\Cal Z(m)$ is the divisor on $\M$
corresponding to those $O_B$--abelian varieties which admit a special
endomorphism
$x$ with $x^2=-m$.
These cycles can be viewed as the Shimura curve analogues of the
cycles on the modular curve defined by elliptic curves with CM by the order
$\Z[\sqrt{-m}]$. For $m<0$, $\Cal Z(m)=\emptyset$. For all $m\ne0$,
$\Xi(m,v)$ is the (non-standard) Green's function introduced in \cite{\annals}.
The class $\hat\Cal Z(0,v)$ will be defined presently.

The moduli stack $\M$ carries a universal abelian variety $\Cal A/\M$, and the
Hodge
bundle $\o$ on $\M$ is defined by
$$\o = \wedge^2(\,Lie(\Cal A))^*.\tag0.2$$
We equip $\o$ with the metric $||\ ||$ which, for $z\in \M(\C)$
is given by
$$||\a||^2_z =
e^{-2C}\cdot\frac{1}{4\pi^2}\,\bigg|\int_{A_z(\C)}\a\wedge\bar\a\,\bigg|,\tag0.3$$
where
$$C= \frac12\big(\,\log(4\pi) +\gamma\,\big),\tag0.4$$
for Euler's constant $\gamma$.  The reason for this normalization will be
explained below.
We thus obtain a class $\hat\o = (\o,||\ ||)\in \Pich(\M)$, and we set
$$\hat\Cal Z(0,v) = -\bigg(\, \hat\o+(0,\log(v))\,\bigg) \in
\CH^1(\M).\tag0.5$$

Using the Gillet--Soul\'e height pairing $\langle\ ,\ \rangle$
between $\CH^1(\M)$ and $\Pich(\M)$, \cite{\gsihes}, we form the height
generating series
$$\phi_{\text{\rm height}}(\tau) = \sum_m \lan\,\hat\Cal
Z(m,v),\hat\o\ran\,q^m,\tag0.6$$
where, for $\tau = u+iv$ in the upper half plane $\H$, we have set $q= e(\tau) = e^{2\pi
i\tau}$.
The quantities $\lan\,\hat\Cal Z(m,v),\hat\o\ran$ can be thought of as
arithmetic degrees
\cite{\bostgilletsoule}, \cite{\bostbourb}.
At the same time, we can define the more
elementary generating series
$$\phi_{\text{\rm degree}}(\tau) = \sum_m \deg(\hat\Cal Z(m,v))\,q^m\ \ ,\tag0.7$$
where
$\deg(\hat \Cal Z(m,v)) = \deg(\Cal Z(m))$
is simply the usual (geometric) degree of the $0$--cycle $\Cal Z(m)_\C$ on the
complex Shimura curve $\M_\C$.

To see that these generating series are the $q$--expansions of modular forms,
we consider
a family of Eisenstein series.
In 1975, Zagier \cite{\zagier}, \cite{\cohen}, introduced a (non-holomorphic)
Eisenstein series of weight $\frac32$, whose Fourier expansion is given by
$$\Cal F(\tau) = -\frac1{12} + \sum_{m>0} H(m)\,q^m + \sum_{n} \frac1{16\pi}\,
v^{-\frac12}
\int_1^\infty e^{-4\pi n^2 v r}\, r^{-\frac32}\,dr\, q^{-n^2},\tag0.8$$
where $H(m)$ is the number of classes of
positive definite integral binary quadratic forms of discriminant $-m$.
This series, which played a key role in the work of Hirzebruch and Zagier
\cite{\hirzebruchzagier} on
generating functions for intersection numbers of curves on Hilbert
modular surfaces, can be viewed as the value at $s=\frac12$ of an Eisenstein
series $\Cal F(\tau,s)$, defined for $s\in \C$ and satisfying a functional
equation $\Cal F(\tau,-s) = \Cal F(\tau,s)$.  In fact, there is a whole
family of such series, $\Cal E(\tau,s;D)$, where $D$ is a square free positive
integer, whose values at $s=\frac12$, for $D>1$, are given by\footnote{Our
normalization
of these series differs slightly at $2$ from that used by Zagier,
so that our $\Cal E(\tau,\frac12;1)$ is not quite Zagier's function,
cf (8.24) below.}
$$\Cal E(\tau,\frac12;D) = -(-1)^{\ord(D)}\frac1{12}\,\prod_{p\mid D} (p-1)
+ \sum_{m>0} 2\,\delta(d;D)\,H_0(m;D)\,q^m.\tag0.9$$
Here $H_0(m;D)$ is a variant of the class number $H(m)$, and $\delta(d;D)$
is either $0$ or a power of $2$, cf. (8.19) and (8.20) respectively, and
$\ord(D)$ is the number of prime factors of $D$.
In the case $D=D(B)>1$, a simple calculation of $\deg(\Cal Z(m))$ proves the
(known) relation
$$\phi_{\text{\rm degree}}(\tau) = \Cal E(\tau,\frac12;D(B)),\tag0.10$$
so that the value of $\Cal E(\tau,s;D(B))$ at $s=\frac12$ is the degree
generating
function.
The main result of this paper asserts that the second term in the Laurent
expansion of
the Eisenstein series $\Cal E(\tau,s;D(B))$ at the point $s=\frac12$ contains
 information about the
arithmetic surface $\M$:
\proclaim{Theorem A} For $D(B)>1$,
$$\phi_{\text{\rm height}}(\tau) = \Cal E'(\tau,\frac12;D(B)) + \bold c$$
for some constant $\bold c$.
\endproclaim

This identity is proved by a direct computation of the two sides. The resulting
formulas,
cf.~Theorem~8.8, are quite
complicated. For example,
for $m>0$ in the case in which $\deg(\Cal Z(m)_\C)\ne0$, the coefficient of
$q^m$
in the derivative of the Eisenstein series is given by
$$\align
&2\,\delta(d;D(B))\,H_0(m;D(B))\cdot \bigg[ \,
\frac12\, \log(d)
 + \frac{L'(1,\chi_d)}{L(1,\chi_d)} -\frac12\log(\pi) -\frac12\gamma \tag0.11\\
\nass
\nass
{}&\hbox to .7in{}
+\frac12 J(4\pi mv)
+\sum_{p \atop p\nmid D(B)} \bigg(\ \log|n|_p - \frac{b_p'(n,0;D)}{b_p(n,0;D)}\
\bigg)
+\sum_{p \atop p\mid D(B)} K_p\,\log(p)
\ \bigg].\\
\endalign
$$
Here we write the discriminant of the order $\Z[\sqrt{-m}]$ as $4m = n^2 d$ for
a fundamental
discriminant $-d$, and the other notation is explained in Theorem~8.8. Theorem
A asserts that this
expression coincides with the height pairing $\lan \hat\Cal Z(m,v),\hat\o\ran$!
A point over $\bar\Q$ of $\Cal Z(m)$ corresponds to an $O_B$--abelian
surface $A$ over $\bar\Q$, equipped with an action of $\Z[\sqrt{-m}]$ commuting
with that of
$O_B$. Such a surface is {\it isogenous} to a product $E_d\times E_d$, where
$E_d$
is an elliptic curve with complex multiplication by the maximal order
$O_\smallkay$
in the imaginary quadratic field $\Q(\sqrt{-d})$. With our normalization of the
metric on the Hodge bundle, the Faltings height of $E_d\times E_d$ is given by
$$\hfal^*(E_d\times E_d) = 2\,\hfal^*(E_d) =
\frac12\log(d)+\frac{L'(1,\chi_d)}{L(1,\chi_d)}
- \frac12\log(\pi) -\frac12\gamma,\tag0.12$$
so that the geometric meaning of the first terms of (0.11),
and of our title, emerges.
The change in the Faltings height due to the isogeny is accounted for by the
term involving the
sum over $p\nmid D(B)$, where the logarithmic derivatives occurring there are given explicitly
in
Lemma~8.10.  The sum over $p\mid D(B)$ has the following geometric meaning. The
arithmetic surface
$\M$ has bad reduction at such primes and the cycle
$\Cal Z(m)$, defined as a moduli space, can include components of the special
fiber $\M_p$, i.e.,
vertical components \cite{\krinvent}. Their contribution to the height pairing
coincides with
the term $K_p\cdot\log(p)$ in (0.11), where $K_p$ is given explicitly in
Theorem~8.8. Finally, there is an additional
`archimedean' term in the height pairing, which arises from the fact that the
Green's function
$\Xi(m,v)$ is not orthogonal to the Chern form $\mu=c_1(\hat\o)$ of the Hodge
bundle. This contribution
coincides with the term involving $\frac12J(4\pi m v)$.

The ambiguous constant $\bold c$ occurs in Theorem A because we do not know the
exact value of the quantity $\lan\hat\o,\hat\o\ran$. More precisely,
for $m=0$, the constant term of the derivative of the Eisenstein series at
$s=\frac12$ with
$D=D(B)>1$ is given by
$$\Cal E'_0(\tau,\frac12;D(B))
=\zeta_D(-1)\,\bigg[\frac12\,\log(v) - 2\frac{\zeta'(-1)}{\zeta(-1)} -1 + 2C+
\sum_{p\mid D}\frac{p\log(p)}{p-1}\,\bigg],\tag0.13$$
where $\zeta_D(s) = \zeta(s)\prod_{p\mid D}(1-p^{-s})$.
On the other hand, by (0.5), the constant term of the generating function for heights is
given by
$$\phi_{\text{\rm height},0}(\tau) = \lan \hat\Cal Z(0,v),\hat\o\ran
= -\lan \hat\o,\hat\o\ran -\frac12\,\log(v)\,\deg(\o).\tag0.14$$
Noting that $\deg(\o) = -\zeta_D(-1)$, we see that the constant terms would
coincide as well, i.e., the constant $\bold c$ in Theorem A would vanish, if
$$\lan \hat\o,\hat\o\ran \overset{??}\to{=}\ \zeta_D(-1)\,\bigg[
2\frac{\zeta'(-1)}{\zeta(-1)}
+1 - 2C- \sum_{p\mid D}\frac{p\log(p)}{p-1}\,\bigg].\tag0.15$$
If we write $\hat\o_o = (\o,||\ ||_{\text{\rm nat}})$ for the Hodge bundle with
the more standard choice of
metric, cf. (10.15) below and \cite{\bost}, and if we delete the `extra' factor
of $\frac12$ which occurs due to the
fact that $\M$ is a stack, cf. section 4, then (0.15) amounts to
$$\lan \hat\o_o,\hat\o_o\ran^{\text{\rm nat}} \overset{??}\to{=}\
4\,\zeta_D(-1)\,
\bigg[ \frac{\zeta'(-1)}{\zeta(-1)} +\frac12 - \frac12\sum_{p\mid
D}\frac{p\log(p)}{p-1}\,\bigg].\tag0.16$$
If we take formally $D=D(B)=1$, so that $\M$ would be the modular curve and
$\hat\o_o$ the bundle of modular forms
of weight $2$ with the Petersson metric, then, indeed, by the result of Bost
and K\"uhn, \cite{\bostumd},\cite{\kuehn},
$$\lan \hat\o_o,\hat\o_o\ran^{\text{\rm nat}} = 4\,\zeta(-1)\,
\bigg[ \frac{\zeta'(-1)}{\zeta(-1)} +\frac12\,\bigg].\tag0.17$$
These considerations are one of the motivations for our choice of metric on the
Hodge bundle and our definition of $\hat\Cal Z(0,v)$.

We expect that Theorem A will continue to hold when $D(B)=1$, i.e., in the case of the
modular curve $\M$, where
(0.17) will allow us to eliminate the constant $\bold c$. There are, however, extra
complications.
The first is that the metric on $\hat\o$ becomes singular at the cusp. This
difficulty was overcome
by Bost \cite{\bost}, \cite{\bostumd} and K\"uhn \cite{\kuehn} by extending the definition of $\Pich(\M) \simeq \CH^1(\M)$
to allow
more general Green's functions. On these more general Chow groups, the
geometric degree map should be defined as
$$\deg:\CH^1(\M) \lra \R, \qquad (Z,g_Z) \mapsto \int_{\M(\C)} \o_Z,\tag0.18$$
where $\o_Z$ is the (now not necessarily smooth) form occuring on the right
hand side
of the Green's equation
$$dd^c\,g_Z + \delta_Z = \o_Z.\tag0.19$$
Note that this definition agrees with the previous one in the case $D(B)>1$.
With our previous definition of $\hat\Cal Z(m,v)$ for $m\ne0$ and $m=-n^2$ and with a slight
modification when $m=0$ or $m=-n^2$,
the result of Funke \cite{\funkethesis}, \cite{\funkecompo} shows that, indeed,
$$\align
\phi_{\deg{}}(\tau) &= \Cal E(\tau,\frac12;1)\\
\nass
{}&=
 -\frac1{12} + \sum_{m>0} 2\,H_0(m;1)\,q^m+
\sum_{n\in\Z}\frac{1}{8\pi}\,v^{-\frac12}\,\int_1^\infty e^{-4\pi n^2
vr}\,r^{-\frac32}\,dr \cdot q^{-n^2}\ .
\tag0.20\endalign
$$
In fact, if we had defined the cycles $\Cal Z(m)$ by imposing the action of
an order of discriminant $m$ (rather than $4m$) on our $O_B$--abelian surface,
then the
degree generating function would coincide exactly with Zagier's function (0.8)!
We defer the calculations of the additional terms which occur in the derivative
$\Cal E'(\tau,\frac12;1)$ and in the height generating
function for $D(B)=1$ to a sequel to this paper \cite{\kryIII}.

Relations like the ones proved here between the first (resp.
second) term of
the Laurent expansion of an Eisenstein series and the generating function for
degrees (resp. heights) should hold in much greater generality. More precisely,
suppose that
$V$ is a rational vector space with nondegenerate inner product of signature
$(n,2)$.
Let $H=\GSpin(V)$ and let $D$ be the space of oriented negative $2$--planes in
$V(\R)$.
Then, for each compact open subgroup $K\subset G(\A_f)$, there is a
quasiprojective variety
$X_K$, defined over $\Q$, with
$$X_K(\C) \simeq H(\Q)\back\bigg(\, D\times H(\A_f)/K\,\bigg).\tag0.21$$
For each integer $m>0$, and each $K$--invariant `weight function' $\ph\in
S(V(\A_f))^K$
in the Schwartz space of $V(\A_f)$, there is a divisor $Z(m,\ph)_K$ on $X_K$,
rational over $\Q$,
\cite{\duke}.
The variety $X_K$ comes equipped with a metrized line bundle $\hat\Cal L=(\Cal
L, ||\ ||)$, and it is proved in
\cite{\Bints} that, with the exception of the cases $n=1$, $V$ isotropic and $n=2$, $V$ split,
 the degree generating function
$$\phi_{\deg}(\tau;\ph) :=\vol(X_K)\,\ph(0) + \sum_{m>0}
\deg(Z(m,\ph))\,q^m\tag0.22$$
coincides with the value $E(\tau,\frac{n}2;\ph)$ of an Eisenstein series
$E(\tau,s;\ph)$
of weight $\frac{n}2+1$ associated to $\ph$. Here $\vol(X_K)$ (resp.
$\deg(Z(m,\ph))$)
is the volume of $X_K(\C)$, (resp. $Z(m,\ph)_K$) with respect to $\O^n$ (resp.
$\O^{n-1}$), where
$\O$ is the negative of the first Chern form of $\hat\Cal L$. We believe that
there
should be an analogue of Theorem A in this situation. To obtain such a result,
one needs,
first of all, suitable extensions $\Cal Z(m;\ph)$ of the cycles $Z(m,\ph)$ to
suitable integral models
$\frak X_K$
of the $X_K$'s. Next, since the varieties $X_K$ are, in general, not
projective, one needs
nice compactifications $\bar \frak X_K$ and, more importantly, an extension of
the Gillet--Soul\'e theory,
general enough to allow the singularities of the metric on the extension
$\hat\o$ of $\hat\Cal L^\vee$ to
the compactification, etc. Assuming all of this, one would have cycles
$$\widehat{\Cal Z}(m,v;\ph) = (\Cal Z(m;\ph),\Xi(m,v))\in \CH^1(\bar\frak
X)_K\tag0.23$$
and a class $\hat\o\in \CH^1(\bar\frak X)_K$. The analogue of Theorem A
would identify the height generating series
$$\phi_{\text{\rm height}}(\tau;\ph) : = \sum_{m} \lan\widehat{\Cal
Z}(m,\ph;v), \hat\o^n\ran\, q^m\tag0.24$$
with the derivative $\Cal E'(\tau,\frac{n}2;\ph)$ at $s=\frac{n}2$ of a
normalized (and possibly slightly modified, cf.
section~6 below)  version $\Cal E(\tau,s;\ph)$ of the Eisenstein series
$E(\tau,s;\ph)$.
However, it seems a challenge to go beyond the case considered in the
present paper and to obtain such results for more general level structures
(even for $n=1$) and for higher values of $n$, e.g.\ for $n=2$
(Hilbert-Blumenthal surfaces) or $n=3$ (Siegel threefolds).
Nonetheless,  the results of \cite{\Bints} provide some additional evidence
in favor of this picture for general $n$.

With hindsight, it may be said that the results in \cite{\tiny} support our
picture in the case $n=0$. This is the one case which is common to the
general picture developed here and the general picture of our papers \cite{\annals}, \cite{\krsiegel}, \cite{\krcrelle},
\cite{\krinvent}, \cite{\tiny}, \cite{\kbourb}.
In this rather
degenerate case,
the variety $X_K$ is zero dimensional, so that the cycles $Z(m;\ph)$ are
actually empty,
and the degree generating function $\phi_{\deg}(\tau;\ph)$ is identically zero.
On the other hand, the associated Eisenstein series $\Cal E(\tau,s;\ph)$ of
weight $1$ is incoherent, in the sense of \cite{\annals}, \cite{\kbourb}, so that $\Cal E(\tau,0;\ph)=0$
as well. The
main result of \cite{\tiny}, Theorem~3, may be interpreted as the identity
$$\phi_{\roman{height}}(\tau, \varphi_0)={\Cal E}'(\tau,0;\varphi_0)\ \
.\leqno(0.26)$$
Here, as in the present paper, $\varphi=\varphi_0$ is the characteristic
function of a certain standard lattice and $K$ is the maximal open compact subgroup.
To make the transition
to the result in \cite{\tiny} one has to take into account the following two
remarks. First, the arithmetic degree $\widehat{\roman{deg}}(\cdot)$ on
$\widehat{CH}^1({\frak X}_K)$ used in \cite{\tiny} may be viewed as
$$\widehat{\roman{deg}}(Z)=\langle Z,\hat\omega^0\rangle\ \ .\leqno(0.27)$$
Second, in defining the degree generating function in this case, we set
$$\hat Z (0,v)=\hat\omega +(0,{\roman{log}}\ v)\in\widehat{CH}^1({\frak X}_K)\ \ ,\leqno(0.28)$$
where $\hat\omega$ is the Hodge bundle on ${\frak X}$ (the moduli stack of
elliptic curves with complex multiplication by $\Cal O_{q}$ for a prime
$q\equiv 3\mod(4)$), with metric {\it normalized as in (10.16)} of the present paper.
Indeed, this particular choice of normalization, i.e., the choice of the
constant $C$ in (0.3), was motivated by the requirement that no ambiguous constant
like $\bold c$ in Theorem A should arise in (0.26).
Specifically, the constant term
which occurs in \cite{\tiny} is given by
$$\align
a_0(\phi,v) &= - 2\, h(\kay)\,\big(\,\frac12\,\log(v) + 2 \,\hfal(E) -
\log(2\pi) -\frac12\log(\pi) +\frac12\gamma +2\log(2\pi)\,\big)\\
\nass
{}&= - 2\, h(\kay)\,\big(\,\frac12\,\log(v) + 2 \,\hfal(E)
+\frac12\log(\pi)+\frac12\gamma+\log(2)\,\big)\tag0.29\\
\nass
{}&= - 2\, h(\kay)\,\big(\,\frac12\,\log(v) + 2 \,\hfal^*(E)\,\big).
\endalign
$$
Here the quantity $\hfal(E)-\frac12\log(2\pi)$ is the Faltings height in
the
normalization of Colmez \cite{\colmez}, which was used in (0.16) of
\cite{\tiny}, cf. Proposition~10.10 below.
We found it particularly striking that the normalization of the metric
on the Hodge bundle which eliminates any garbage constant in the case $n=0$
also gives a precise match in the positive Fourier coefficients in our Shimura
curve case ($n=1$).
Of course, this is perhaps not so surprising, given the fact that the
cycles in the case of signature $(n,2)$ are themselves (weighted) combinations
of Shimura varieties of the same type for signature $(n-1,2)$. Thus a main term
in
the arithmetic degrees
which occur in the positive Fourier coefficients of the height generating
function for the
signature $(n,2)$ case is the `arithmetic volume' occuring in the constant term
for the
$(n-1,2)$ case. This `explains' the relation between the present paper and the
results of
\cite{\tiny}. It should not be difficult to verify that the positive Fourier
coefficients of the
derivatives of Eisenstein series of weight $\frac{n}2+1$ at $s=\frac{n}2$
are related to the constant terms of the derivatives of those of weight
$\frac{n-1}2+1$ at $s=\frac{n-1}2$
in a similar way.

In a similar vein, we remark that the height generating function which,
according to our picture above, is related to the derivative of Eisenstein
series  on $\SL_2 = \text{\rm Sp}_1$ of weight ${n\over 2}+1$ at $s={n\over 2}$, is connected with the
singular Fourier coefficients of the derivative of Eisenstein series of genus $2$, i.e., on $\text{\rm Sp}_2$, of
weight ${n\over 2}+1$ at $s={n-1\over 2}$. In fact, this is how
we arrived at the height generating function considered in this paper. We
hope to elaborate on this point in a future paper.

The results of this paper are an outgrowth of a project begun during the first author's
visit to the Mathematische Institut of the University of Cologne in the fall of 1999.
He would like to thank the Institut for providing a congenial and stimulating working environment.
The second author thanks the department of mathematics
of the University of Maryland for its (by now almost customary) hospitality
during his sabbatical in the spring of 2001.

\subheading{Contents}
\medskip
\centerline{Introduction\hfill}
\centerline{Part I. Arithmetic geometry\hfill}
\centerline{1. The moduli stack $\M$\hfill}
\centerline{2. Uniformization\hfill}
\centerline{3. The Hodge bundle\hfill}
\centerline{4. Arithmetic Chow groups\hfill}
\centerline{5. Special cycles and the generating function\hfill}
\centerline{Part II. Eisenstein series\hfill}
\centerline{6. Eisenstein series of weight 3/2\hfill}
\centerline{7. The main identity \hfill}
\centerline{8. Fourier expansions and derivatives\hfill}
\centerline{Part III. Computations: geometric\hfill}
\centerline{9. The geometry of $Z(m)$'s\hfill}
\centerline{10. Contributions of horizontal components\hfill}
\centerline{11. Contributions of vertical components\hfill}
\centerline{\qquad Appendix to section 11: The case $p=2$\hfill}
\centerline{12. Archimedean contributions\hfill}
\centerline{13. Remarks on the constant term\hfill}
\centerline{Part IV. Computations: analytic\hfill}
\centerline{14. Local Whittaker functions: the non--archimedean case\hfill}
\centerline{15. Local Whittaker functions: the archimedean case\hfill}
\centerline{16. The functional equation\hfill}
\centerline{References\hfill}

\comment
\subheading{Notation}

$\hat{\Z}$, $\A_f$,

Let $\psi = \otimes \psi_p$ be the `canonical' additive character
of $\Bbb Q \backslash \A$ given by
$$
\psi_p(x) = \cases
e(x)  &\text{ if $ p = \infty$,}\\
e(-\lambda(x))  &\text{ if $ p < \infty$.}
\endcases
$$
Here $e(x) = e^{ 2 \pi i x}$ and $\lambda: \Q_p\rightarrow \Q_p/\Z_p \simeq \Z[\frac1{p}]/\Z\hookrightarrow  \Q/\Z $.

$\gamma$ is Euler's constant.

\endcomment

\subheading{Part I. Arithmetic geometry}

\subheading{1. The moduli stack $\M$}

Let $B$ be an indefinite quaternion algebra over $\Q$. We fix a maximal order
$O_B$ in $B$, and we let $D(B)$ be the product of the primes $p$ at which $B_p$
is a division algebra. For the moment, we allow the case $B=M_2(\Q)$ and $O_B=M_2(\Z)$, where $D(B)=1$.

We denote by $\M$ the stack over $\Spec\Z$ representing the following
moduli problem. The moduli problem associates to a scheme $S$ the category $\M(S)$
whose objects are pairs $(A,\iota)$, where $A$ is an abelian scheme over $S$ and
$$\iota:O_B\lra \End_S(A)$$
is a homomorphism such that, for $a\in O_B$,
$$\det(\iota(a);\Lie(A)) = Nm^o(a).\tag1.1$$
Here $Nm^o$ is the reduced norm on $B$ and, as usual, \cite{\kottwitz}, \cite{\RZ},
the identity (1.1) is meant as an identity of polynomial functions on $S$.
All morphisms in this category are isomorphisms.

\proclaim{Proposition 1.1} $\M$ is an algebraic stack in the sense of
Deligne-Mumford. Furthermore, $\M$ is proper over $\Spec\Z$ if
$B$ is a division algebra. The restriction of $\M$ to $\Spec\Z[D(B)^{-1}]$
is smooth of relative dimension $1$. Finally, if $p\mid D(B)$,
then $\M\times_{\Spec\Z}\Spec\Z_p$ has semi-stable reduction. \qed
\endproclaim

\subheading{2. Uniformization}

Let $H= B^\times$, considered as an algebraic group over $\Q$. Let
$$D = \Hom_\R(\C,B_\R),\tag2.1$$
be the set of homomorphisms of $\R$--algebras, taking $1$ to $1$,
with the natural conjugation action of $H(\R)$. This action is transitive. We fix an isomorphism
$B_\R\simeq M_2(\R)$, so that
$H(\R)\simeq GL_2(\R)$,  and a compatible isomorphism
$D\simeq \C\setminus \R$, the union of the upper and lower half planes.
Also let $K= \hat{O}_B^\times \subset H(\A_f)$, be the compact open subgroup determined by
$O_B$, where $\hat{O}_B = O_B\tt_\Z\hat\Z$. Then we have, as usual,
an isomorphism of Deligne-Mumford stacks over $\C$,
$$\M\times_{\Spec\Z}\Spec\C = \left[\,H(\Q)\back  D\times H(\A_f)/K\,\right],\tag2.2$$
where the right is to be understood in the sense of stacks, \cite{\delignemumford}, p.99.
The stack on the right hand side may be written in a simpler way using the fact
that $H(\A_f) = H(\Q)K$. Let
$$\Gamma = H(\Q)\cap K = O_B^\times.\tag2.3$$
Then
$$\M\times_{\Spec\Z}\Spec\C = [\Gamma\back D].\tag2.4$$
Note that $\Gamma$ acts on $D$ through its image $\bar{\Gamma} = \Gamma/\{\pm1\}$
in $\text{\rm PGL}_2(\R)$, with finite stabilizer groups. Instead of
considering $[\Gamma\back D]$ as an algebraic stack, it is more traditional
to view this quotient as an orbifold \cite{\Hain}.  Intuitively speaking, this means
that the quotient of $D$ by the action of $\Gamma$ is not carried out, but rather,
all information obtained from the action of $\Gamma$ on $D$ is stored. As
a particular instance, consider the hyperbolic volume form $\mu$ on $D$, normalized as
$$\mu=\frac{1}{2\pi}\,y^{-2}\,dx\wedge dy,\tag2.5$$
in standard coordinates on $\C\setminus \R$. Since this volume form is $\Gamma$--invariant,
it induces a volume form on the orbifold $[\Gamma\back D] = \M(\C)$.
The volume of the orbifold $\M(\C)$ is given by
$$\vol(\M(\C)) =\int_{[\Gamma\back D]} \mu= \frac12 \int_{\Gamma\back D} \mu,\tag2.6$$
where the extra factor of $\frac12$ in the second expression is due to the fact that the stablizer
in $\Gamma$ of a generic point of $D$, has order $2$.
Explicitly, we have, \cite{\Eichler},
$$\vol(\M(\C)) = \frac{1}{12}\prod_{p|D(B)} (p-1) ={} - \zeta_D(-1).\tag2.7$$
where
$$\zeta_D(s) = \zeta(s)\,\prod_{p\mid D} (1-p^{-s}).$$

We now turn to $p$--adic uniformization, \cite{\drinfeld}, \cite{\boutotcarayol},
\cite{\krinvent}. We fix a prime $p\mid D(B)$.
Let $B'$ be the definite quaternion algebra over $\Q$ whose invariants agree with those
of $B$ at all primes $\ell\ne p,\ \infty$. Let $H' = B^{\prime\,\times}$
considered as an algebraic group over $\Q$. We fix identifications
$H'(\A_f^p)\simeq H(\A_f^p)$ and $H'(\Q_p) = GL_2(\Q_p)$.
Let $\hat\O^2$ be the Deligne-Drinfeld formal scheme relative to $GL_2(\Q_p)$.
Then
$$\M\times_{\Spec\Z}\Spec W(\bar\F_p) \simeq \left[\,H'(\Q)\back
 \left(\,\hat\O^2\times_{\Spf\,\Z_p}\Spf\, W(\bar\F_p)\tent{10}\,\right)\times\Z\times H'(\A_f^p)/K^{\prime\,p}\,\right].\tag2.8
$$
Here $K^{\prime\,p}$ corresponds to $(O_B\tt\hat\Z^p)^\times$ under the
identification of $H(\A_f^p)$ with $H'(\A_f^p)$,
and $g\in H'(\Q)$ acts on the $\Z$ factor by shifting by $\ord_p(\det(g))$. Again, this
formula can be simplified since $H'(\Q)K_f^p$
maps surjectively onto $\Z \times H'(\A_f^p)$. Let
$$H'(\Q)^1 = \{\ g\in H'(\Q)\mid \ord_p(\det(g)) = 0\ \}.\tag2.9$$
Put $\Gamma' = H'(\Q)^1\cap K_f^p$. Then
$$\M\times_{\Spec\Z}\Spec W(\bar\F_p) \simeq \left[\ \Gamma' \back \hat\O^2\times_{\Spec \Z_p} \Spec W(\bar\F_p)\ \tent{10}\right],\tag2.10$$
where, again the right hand side is considered as (the algebraization of) a formal
Deligne-Mumford stack. The group $\Gamma'$ acts through $\bar\Gamma' = \Gamma'/\{\pm1\}\subset \text{\rm PGL}_2(\Q_p)$
with finite stabilizer groups.

\subheading{3. The Hodge bundle}

We denote by $(\Cal A,\iota)$ the universal abelian scheme over $\M$
and by $\e:\M\rightarrow \Cal A$, its zero section. The
{\it Hodge line bundle} on $\M$ is the following line bundle:
$$\o = \e^*(\Omega^2_{\Cal A/\M}) = \wedge^2 \Lie(\Cal A/\M)^*.\tag3.1$$
For convenience, we will refer to $\o$ as the Hodge bundle.
\demo{Remark 3.1} Assume that $(B,O_B) = (M_2(\Q),M_2(\Z))$. In this case,
$\M$ may be identified with the moduli stack of elliptic curves, and the
universal object $\Cal A$ with $(\Cal E^2,\iota_0)$ where $\Cal E$ is the universal
elliptic curve and $\iota_0:M_2(\Z)\rightarrow \End(\Cal E^2) = M_2(\End(\Cal E))$
is the natural embedding. In this case, $\Lie(\Cal A/\M) = \Lie(\Cal E/\M)^{\oplus2}$
and hence
$$\o = \o_{\Cal E/\M}^{\tt2}.\tag3.2$$
Here $\o_{\Cal E/\M} = \Lie(\Cal E/\M)^*$.
Recall \cite{\DR}, VI.4.5, that $\o_{\Cal E/\M}^{\tt2}$ can be identified with the module of
relative differentials of $\M/\Spec\Z$. The following proposition
generalizes this fact.
\qed\enddemo

\proclaim{Proposition 3.2}
The Hodge bundle $\omega$ is isomorphic to the relative dualizing sheaf
$\omega_{{\Cal M}/\Z}$.
\endproclaim

\demo{Proof}
Since the fibers of ${\Cal M}$ over $\Spec\, \Z$ are Gorenstein,
$\omega_{{\Cal M}/\Z}$ is an invertible sheaf. Since ${\Cal M}$ is regular
of dimension 2, it suffices to show that the restrictions of $\omega$ to
the smooth locus ${\Cal M}^{\roman{smooth}}$ is isomorphic to the
restriction of $\omega_{{\Cal M}/\Z}$ to ${\Cal M}^{\roman{smooth}}$,
i.e., to the sheaf of relative differentials $\Omega^1_{{\Cal
M}^{\roman{smooth}}/\Z}$.
\par
By deformation theory we have a canonical identification
$$\Omega^1_{{\Cal M}^{\roman{smooth}}/{\Z}}=  {Hom}_{O_B}(\Lie\, {\Cal A},
(\Lie(\hat{\Cal A})^*)\ \ ,\tag3.3$$
where $({\Cal A},\iota)$ is, as before, the universal object over ${\Cal
M}$ and where $\hat{\Cal A}$ denotes the dual abelian variety.
\par
This formula shows that it suffices to check the claimed equality
after passing to the completion $\Z_p$ for each prime $p$. For any prime $p\nmid D(B)$ our
identification problem
reduces to the situation considered in Remark 3.1. Hence, all we need to do
is to extend to ${\Cal M}^{\roman{smooth}}$ the isomorphism between $\omega$ and
$\Omega^1_{{\Cal M}/\Spec}$
over ${\Cal M}[D(B)^{-1}]$. Fix a prime
number $p\mid D(B)$. Denote by $R$ the involution on $O_B$
$$\alpha\longmapsto \alpha^R=\delta\,\alpha^\iota\,\delta^{-1},\tag3.4$$
where $\delta\in O_B$ satisfies $\delta^2=D(B)$, and where
$\alpha\mapsto\alpha^\iota$ is the main involution on $B$. After choosing a
$p$--principal polarization on ${\Cal A}$ whose Rosati involution induces
the involution $R$ on
$O_B$, we may identify $(\Lie\, \hat{\Cal A})^*$ with $(\Lie\, {\Cal A})^*$, in
such a way that $\alpha\in O_B$ acts on $(\Lie\, {\Cal A})^*$ as
$(\alpha^R)^*$. We write
$$O_{B_p}=\Z_{p^2}[\Pi]/(\Pi^2=p,\ \Pi a=a^{\sigma}\Pi,\ \ \forall
a\in\Z_{p^2})\ \ .\tag3.5$$
We may assume that the restriction of $R$ to $\Z_{p^2}$ is trivial and that $\Pi^R=-\Pi$.
After extending scalars from $\Z$ to $\Z_{p^2}$, we have the eigenspace
decomposition of $\Lie\, {\Cal A}$ as
$$\Lie\, {\Cal A}={\Cal L}_0\oplus {\Cal L}_1\ \ ,\tag3.6$$
such that the action of $\Pi$ on $\Lie\, {\Cal A}$ is of degree 1 with
respect to this $\Z/2$-grading. The condition (1.1) ensures that ${\Cal
L}_0$ and ${\Cal L}_1$ are both line bundles. It now follows,
via (3.3) and (3.6), that the
local sections of $\Omega^1_{{\Cal M}/\Z}$ are given by local
homomorphisms $\varphi_i: {\Cal L}_i\to {\Cal L}_i^*$ $(i=0,1)$ forming a
commutative diagram
$$\matrix
{\Cal L}_0
&
\buildrel\varphi_0\over\lra
&
{\Cal L}_0^*
\\
\llap{$\scriptstyle\Pi_0$}\big\downarrow
&&
\big\downarrow\rlap{$\scriptstyle -\Pi_1^*$}
\\
{\Cal L}_1
&
\buildrel\varphi_1\over\lra
&
{\Cal L}_1^*
\\
\llap{$\scriptstyle \Pi_1$}\big\downarrow
&&
\big\downarrow\rlap{$\scriptstyle -\Pi_0^*$}
\\
{\Cal L}_0
&
\buildrel\varphi_0\over\lra
&
{\Cal L}_0^*
\endmatrix
\tag3.7
$$
The pair $(\varphi_0,\varphi_1)$ therefore defines an injective
homomorphism
$$\alpha: \Omega^1_{{\Cal M}/\Z}\lra {\Cal L}_0^{\otimes(-2)}\oplus {\Cal
L}_1^{\otimes(-2)}\  \subset\  (\Lie\, {\Cal A})^*\otimes (\Lie\, {\Cal
A}).\tag3.8$$
On the generic fiber, ${\Cal L}_0^{\otimes(-2)}$ and ${\Cal
L}_1^{\otimes(-2)}$ both coincide with $\omega$ and $\alpha$ induces an
isomorphism of $\Omega^1_{{\Cal M}/\Z}$ with the diagonal. On the smooth
locus either $\Pi_0$ or $\Pi_1$ is an isomorphism locally around any given
point. Assume for instance that $\Pi_0$ is an isomorphism. Then
$\varphi_0$ determines $\varphi_1$ by the commutativity of the upper
square in (3.7),
$$\varphi_1=(-\Pi_1^*)\circ \varphi_0\circ \Pi_0^{-1}\ \ .\tag3.9$$
But then also the lower square commutes. Since $\Pi_0$ is an isomorphism,
it suffices to check this after premultiplying $\varphi_1$ with $\Pi_0$.
But
$$(-\Pi_0^*)\circ \varphi_1\circ \Pi_0=(-\Pi_0^*)\circ
(-\Pi_1^*)\circ\varphi_0=\varphi_0\circ \Pi_1\circ \Pi_0\ \ .\tag3.10$$
It follows that on the open sublocus of ${\Cal M}^{\roman{smooth}}$ where
$\Pi_0$ is an isomorphism, the first projection applied to (3.8) induces
an isomorphism between $\Omega^1_{{\Cal M}/\Z}$ and ${\Cal
L}_0^{\otimes(-2)}$. On the other hand on this open sublocus, $\omega$ can
be identified with ${\Cal L}_0^{\otimes(-2)}$, which proves the claim.
\qed
\enddemo

By base change to $\C$, the Hodge bundle induces a line bundle $\o_\C$
on $\M_\C = [\Gamma\back D]$. In the orbifold picture, we may view $\o_\C$ as being
given by a descent datum with respect to the action of $\Gamma$ on the pullback
of $\o_\C$ to $D$. At a point $z$ of $\M_\C$, a section $\a$ of $\o_\C$ corresponds a holomorphic $2$ form
on $\Cal A_z$, and so there is a natural norm \cite{\bost} on $\o_\C$ given by:
$$||\a_z||_{\text{\rm nat}}^2 = \bigg|\left(\frac{i}{2\pi}\right)^2 \int_{\Cal A_z(\C)} \a\wedge\bar\a\ \bigg|.\tag3.11$$

Equivalently, $\o_\C$ is given by the
automorphy factor $(c z + d)^2$, i.e., by the action of
$\Gamma$ on $D\times \C$ defined by
$$\gamma=\pmatrix a&b\\c&d\endpmatrix : (z,\zeta) \mapsto (\gamma(z), (cz+d)^2 \zeta\,).\tag3.12$$
More precisely, for $z\in D \simeq \C\setminus \R$, we have an isomorphism
$$B_\R \simeq M_2(\R) \isoarrow \C^2, \qquad u \mapsto u\cdot \pmatrix z\\1\endpmatrix = \pmatrix w_1\\w_2\endpmatrix,\tag3.13$$
and the corresponding abelian variety $\Cal A_z$ has $\Cal A_z(\C) = \C^2/\Lambda_z$, where $\Lambda_z$
is the image of $O_B$ in $\C^2$. The pullback of $\o_\C$ to $D$ is trivialized via the
section $\a= dw_1\wedge dw_2$.
The Petersson norm,
$||\ ||_{\text{\rm Pet}}$ on the bundle of modular forms of weight $2$, is defined by the $\Gamma$-invariant
norm on the trivial line bundle $D\times \C$ given by
$$||(z,\zeta)||^2_{\text{\rm Pet}} = |\zeta|^2\,(4\pi\Im(z))^2.\tag3.14$$
If $f$ is such a modular form, then we identify $f$ with the section, \cite{\chaifaltings}, pp141--2,
$$\a(f) = f(z)\, (2\pi i\, dw_1\wedge 2 \pi i\, dw_2) = - 4\pi^2 f(z)\,\a.$$

\proclaim{Lemma 3.3} The two norms on $\o_\C$ are related by
$$||\ ||_{\text{\rm nat}}^2= D(B)^2 \,||\ ||^2_{\text{\rm Pet}} .$$
\endproclaim
\demo{Proof}
The pullback to $M_2(\R)$ of the form
$\a\wedge\bar\a= dw_1\wedge dw_2\wedge d\bar w_1\wedge d \bar w_2$ above is $4\,\Im(z)^2$ times the
standard volume form, and so, via (3.11),
$$||\a||_{\text{\rm nat}}^2 =\frac1{4\pi^2}\cdot 4 \,\Im(z)^2\,\vol(M_2(\R)/O_B).\tag3.15$$
Then
$$||\a(f)||^2_{\text{\rm nat}} = 16\pi^4\,|f(z)|^2\cdot \frac{D(B)^2}{\pi^2}\,Im(z)^2.\qed$$
\enddemo

\proclaim{Definition 3.4} The metrized Hodge bundle $\hat\o$ is $\o$ equipped with the
metric
$$||\ || = e^{-C}\,||\ ||_{\text{\rm nat}},$$
where
$$C = \frac12\big(\,\log(4\pi)+\gamma\,\big).$$
Here $\gamma$ is Euler's constant.
\endproclaim

The motivation for this normalization is explained in the introduction.

The Chern form $c_1(\hat\o)$ for this metric is then
$$c_1(\hat\o) = -dd^c\log||\a||^2 = \mu\tag3.16$$
with $\mu$ as in (2.5), and so
$$\deg(\hat\o) = \int_{[\Gamma\back D]} c_1(\hat\o) = \vol(\M(\C)).\tag3.17$$

\subheading{4. The arithmetic Picard group and the arithmetic Chow group}

From now on, we assume that $D(B)>1$, so that $B$ is a division algebra
and $\M$ is proper over $\Spec \Z$. If we had imposed a sufficient level
structure, then $\M$ would be an arithmetic surface over $\Spec \Z$,
\cite{\gsihes}, \cite{\bost}, etc. Then the Chow groups (tensored with
$\Q$) $CH^r(\M)$ and arithmetic Chow groups $\CH^r(\M)$ would be defined,
with $\CH^1(\M)\simeq \Pich(\M)$, the group of isomorphism classes of
metrized line bundles, and these would be equipped with a height pairing
$$\langle\ ,\ \rangle:\CH^1(\M)\times\CH^{1}(\M) \lra \CH^2(\M)\tag4.1$$
and the arithmetic degree $$\degh:\CH^2(\M)\lra \C.\tag4.2$$ In this
section, we explain how to carry over (parts of) this formalism to our
DM--stack $\M$.

We begin with $\Pich(\M)$. There are two ways to define the concept of a metrized line bundle on $\M$.
First, one can define such an object to be a rule which
associates, functorially to any $S$--valued point $S\rightarrow \M$ of
$\M$, a line bundle $\Cal L_S$ on $S$ equipped with a $C^\infty$--metric on the
line bundle $\Cal L_{S,\C}$ on $S\times_{\Spec \Z}\Spec \C$.
Second, one can define a metrized line bundle on $\M$ to be an invertible
sheaf on $\M$ together with a $\Gamma$--invariant metric on the pullback of
$\Cal L_\C$ to $D$ under the identification of $\M_\C$
with the orbifold $[\Gamma\back D]$. These definitions are equivalent.
As usual, we denote the set of isomorphism classes of metrized line bundles
on $\M$ by $\Pich(\M)$. This is an abelian group under the tensor product operation.

Let $\hat\Cal L = (\Cal L,||\cdot ||)$ be a metrized line bundle on $\M$.
Then $\hat\Cal L$ determines a {\it height} of a one dimensional irreducible reduced
proper DM--stack $\Cal Z$ mapping to $\M$.
To define it , we are guided by the heuristic principle that, in a numerical formula, a
geometric point $x$ of
a stack counts with fractional multiplicity $1/|\text{\rm Aut}(x)|$.
Let $\tilde\Cal Z$ be the normalization and
$\nu:\tilde\Cal Z\rightarrow \M$ be the natural morphism.

If $\Gamma(\tilde\Cal Z, {\Cal O}_{\tilde\Cal Z})=O_K$ for a number field $K$, then
$$h_{\hat\Cal L}(\Cal Z)=\widehat{\deg}\,\nu^*(\hat{\Cal L})\ \ \tag4.3$$
Here the right hand side is defined by setting, for a meromorphic section
$s$ of $\nu^*({\Cal L})$,
$$\widehat{\deg}\,\nu^*(\hat{\Cal L})=\sum_p\bigg( \sum_{x\in \tilde \Cal Z(\bar\F_p)}
{\ord_x(s)\over \vert\Aut(x)\vert}\,\bigg)\cdot\log p-{1\over 2}\int_{\tilde
\Cal Z(\C)}\log\Vert s\Vert^2.\tag4.4$$
Here the integral is defined as
$$\int_{\tilde\Cal  Z(\C)}\log\Vert s\Vert^2= \sum_{z\in \tilde\Cal Z(\C)} {1\over
\vert\Aut(z)\vert} \cdot\log \Vert s(z)\Vert^2.\tag4.5$$
Also $\ord_x(s)$ is defined by noting that the strict henselization $\tilde\Cal O_{\tilde\Cal Z,x}$ of the
local ring ${\Cal O}_{\tilde\Cal Z, x}$ is a discrete valuation ring. Let us
check that the expression (4.4) is independent of the choice of $s$. This
comes down to checking for a function $f\in \Q(\tilde\Cal Z)^\times =K^\times$
that
$$0=\sum_p\bigg( \sum_{x\in \tilde\Cal Z(\bar\Bbb F_p)} {\ord_x(f)\over
\vert\Aut(x)\vert}\,\bigg) \cdot\log p-{1\over 2}\sum_{\sigma: K\to\C} {1\over
\vert\Aut(\sigma)\vert} \cdot\log|\sigma(f)|^2\tag4.6$$
For $x\in \tilde\Cal  Z(\bar\F_p)$, let $\und{x}$ be the corresponding geometric point of
the coarse moduli scheme $Z=\Spec\, O_K$ of $\tilde\Cal Z$. Then
$$\tilde{\Cal O}_{Z, \und{x}}=(\tilde{\Cal O}_{\tilde\Cal Z,
x})^{\Aut(x)/\Aut(\bar\eta)}\tag4.7$$
where $\bar\eta$ is any generic geometric point of $\tilde\Cal Z$, and
$\tilde{\Cal O}_{\tilde\Cal Z, x}$ is a totally ramified extension of degree $\vert\Aut(x)\vert
/\vert\Aut(\bar\eta)\vert$ of $\tilde{\Cal O}_{Z,\und{x}}$. Inserting this
into (4.6), we obtain for the right hand side the expression
$${1\over\vert\Aut(\bar\eta)\vert}\bigg( \sum_p\sum_{x\in (\Spec
O_K)(\bar\F_p)}\ord_x(f)\cdot\log
p-\sum_{\sigma}\log\vert\sigma(f)\vert\bigg)\tag4.8$$
which is zero by the product formula for $f\in K^\times$.
\par
If $\Gamma(\tilde\Cal Z, {\Cal O}_{\tilde\Cal Z})=\F_q$, we put
$$
h_{\hat{\Cal L}}(\Cal Z)
=\deg\nu^*({\Cal L})\cdot\log q
=\bigg( \sum_{x\in\tilde \Cal Z(\bar\F_p)} {\ord_x(s)\over
\vert\Aut(x)\vert}\bigg)\cdot\log p, \tag4.9
$$
where $s$ is a meromorphic section of $\nu^*({\Cal L})$. Here
$\deg\nu^*({\Cal L})$ coincides with the definition given in \cite{\DR}, V.4.3.
\par
Next we need to define the (arithmetic) Chow group of ${\Cal M}$. By a
prime divisor on ${\Cal M}$ we mean a closed substack $\Cal Z$ of ${\Cal M}$
which is locally for the \'etale topology a Cartier divisor defined by an
irreducible equation. Let $Z^1({\Cal M})$ be the free abelian group
generated by the prime divisors on ${\Cal M}$. Any rational function
$f\in\Q({\Cal M})^\times$ (i.e.\ a morphism ${\Cal U}\to\A^1$ defined on a
non-empty open substack ${\Cal U}$ of ${\Cal M}$) defines a principal
divisor
$${\roman{div}}(f)= \sum_{\Cal Z}\ord_{\Cal Z}(f)\cdot \Cal Z\ \ ,\tag4.10$$
where the sum is over the prime divisors $\Cal Z$ of ${\Cal M}$, and where we
note that the strict henselization of the local ring at $\Cal Z$, $\tilde{\Cal
O}_{{\Cal M}, \Cal Z}$, is a discrete valuation ring. The factor group of
$Z^1({\Cal M})$ by the group of principal divisors is the Chow group
$CH^1({\Cal M})$, comp.\ \cite{\vistoli}.
\par
Let $\Cal Z\in Z^1({\Cal M})$. Then the divisor $\Cal Z_{\C}$ of ${\Cal
M}_{\C}=[\Gamma\back D]$ is of the form $\Cal Z_{\C}=[\Gamma\back D_{\Cal Z}]$
for a unique $\Gamma$-invariant divisor $D_{\Cal Z}$ of $D$. By a Green's
function for $\Cal Z$ we mean a real current $g$ of degree 0 on $D$ which is
$\Gamma$-invariant and such that
$$\omega = dd^cg+\delta_{D_{\Cal Z}}\tag4.11$$
is $C^\infty$. We denote by $\hat Z^1({\Cal M})$ the group of Arakelov
divisors, i.e., of pairs $(\Cal Z,g)$ consisting of a divisor $\Cal Z$ on ${\Cal M}$
and a Green's function for $\Cal Z$, with componentwise addition. If $f\in
\Q({\Cal M})^\times$, then $f\vert{\Cal M}_{\C}$ corresponds to a
$\Gamma$-invariant meromorphic function $\tilde f_{\C}$ on $D$ and we
define the associated principal Arakelov divisor
$$\widehat{\roman{div}}(f)= (\,{\roman{div}}(f), -\log \vert\tilde
f_{\C}\vert^2\,)\ \ .\tag4.12$$
The factor group of $\hat Z^1({\Cal M})$ by the group of principal
Arakelov divisors is the arithmetic Chow group $\widehat{CH}^1({\Cal M})$.
The groups $\widehat{CH}^1({\Cal M})$ and $\widehat{\roman{Pic}}({\Cal
M})$ are isomorphic. Under this isomorphism, an element $\hat{\Cal L}$
goes to the class of
$$\big(\,\sum_{\Cal Z}\ord_{\Cal Z}(s)\,{\Cal Z},\, -\log\Vert s\Vert^2\,\big)\ \ ,\tag4.13$$
where $s$ is a meromorphic section of ${\Cal L}$. Conversely, if $(\Cal Z,
g)\in Z^1({\Cal M})$, then its preimage under this isomorphism is
$$({\Cal O}(\Cal Z), \Vert\ \Vert)\ \ ,\tag4.14$$
where $-\log\Vert{\bold 1}\Vert^2=g$, with ${\bold 1}$ the canonical
$\Gamma$-invariant section of the pullback of $\hat{\Cal O}(\Cal Z)$ to $D$.
\par
We define a pairing
$$\langle\ ,\ \rangle: \hat Z^1({\Cal M})\times
\widehat{\roman{Pic}}({\Cal M})\longrightarrow \C\tag4.15$$
by formula (5.11) in Bost \cite{\bost},
$$\langle\, (\Cal Z,g),\hat{\Cal L}\,\rangle = h_{\hat{\Cal L}}(\Cal Z)+{1\over
2}\int_{[\Gamma\setminus D]} g\cdot c_1(\hat{\Cal L}).\tag4.16$$
Here $c_1(\hat{\Cal L})$ is the $\Gamma$-invariant form on $D$ defined by
the pullback to $D$ of $\hat{\Cal L}$ (analogously to $c_1(\hat\omega)$ in
section 3 above). The integral is defined as
$$\int_{[\Gamma\setminus D]}g\cdot c_1(\hat{\Cal L})= [\Gamma
:\Gamma']^{-1}\cdot \int_{\Gamma\setminus D} g\cdot c_1(\hat{\Cal L})
,\tag4.17$$
where $\Gamma'= {\roman{ker}}(\Gamma\to\Aut(D))$.
\par
It seems very likely that under the identification $\widehat{CH}^1({\Cal
M})\simeq \widehat{\roman{Pic}}({\Cal M})$, the pairing (4.15) descends to
a symmetric bilinear pairing
$$\lan\ ,\ \ran:\widehat{CH}^1({\Cal M})\times \widehat{CH}^1({\Cal M})\longrightarrow
\C\ \ ,\tag4.18$$
as is the case for arithmetic surfaces.
For ease of expression we will proceed as if this were the case,
although we have not checked it. Since all we will actually use is
the paring (4.15), this will cause no harm.

\subheading{5. Special cycles and the generating function}

In this section, we will define for each $m\in \Z$ and each $v\in \R^\times_+$
a class
$$\hat\Cal Z(m,v) = (\Cal Z(m), \Xi(m,v))\in \CH^1(\M).\tag5.1$$

We first assume that $m>0$. Then we consider the DM--stack $\Cal Z(m)$
classifying triples $(A,\iota,x)$ where $(A,\iota)$ is an object of $\M$
and where $x$ is a {\it special endomorphism},
\cite{\annals},\cite{\krinvent}, with $x^2 = -m$, i.e.,
$$x\in \End(A,\iota), \qquad \tr^o(x) = 0, \qquad x^2 = -m.\tag5.2$$
Then $\Cal Z(m)$ maps to $\M$ by a finite unramified morphism. Furthermore,
$\Cal Z(m)$ is purely one dimensional, except in the following cases,
\cite{\krinvent} and the Appendix to section 11,
$$ \exists p\mid D(B), \ p\ne 2,\qquad\text{such that $m\in \Z_p^{\times, 2}$.}\tag5.3$$
In the cases covered by (5.3), we set $\hat\Cal Z(m,v)=0$.
In all other cases, we define a Green's function for
the unramified morphism $\Cal Z(m)\rightarrow\M$, in the sense of section 4,
as follows (\cite{\annals}).  Let
$$V= \{x\in B\mid \tr^o(x)=0\ \}\tag5.4$$
with quadratic form $Q(x) = - x^2 = N^o(x)$ given by the restriction of the
reduced norm and with
associated inner product $(x,y) = \tr^o(xy^\iota)$. Note that the signature of $V(\R)$ is
$(1,2)$.
As in \cite{\annals}, we can identify $D$ with the space of oriented negative $2$-planes
in $V(\R)$. For $x\in V(\R)$ and $z\in D$, let $\pr_z(x)$ be the projection of $x$ to $z$
and let
$$R(x,z) = - (\pr_z(x),\pr_z(x))\ge0.\tag5.5$$
This quantity vanishes precisely when $\pr_z(x)=0$, i.e., when $z\in D_x$ where
$$D_x= \{\ z\in D\mid (x,z) = 0\ \}.\tag5.6$$
Let $L = V(\Q)\cap O_B$, and let
$$L(m) = \{\ x\in L\mid Q(x)=m\ \}.\tag5.7$$
Then, for $m\in \Z_{\ne 0}$, and $v\in\R^\times_+$, let
$$\Xi(m,v) = \sum_{x\in L(m)} \xi(v^{\frac12}x,z)\tag5.8$$
where
$$\xi(x,z) = -\Ei(-2\pi R(x,z))\tag5.9$$
for the exponential integral
$$-\Ei(-t) =\int_1^\infty e^{-tr}\,r^{-1}\,dr.\tag5.10$$
The properties of this function are described in \cite{\annals}, section 11.
For $m>0$, $\Xi(m,v)$ is a $\Gamma$--invariant Green's function for the divisor
$$D_{\Cal Z(m)} := \coprod_{x\in L(m)} D_x\tag5.11$$
in $D$.

When $m<0$,  $\Xi(m,v)$ is a smooth $\Gamma$--invariant function on $D$. Therefore,
$$\hat\Cal Z(m,v) = (0,\Xi(m,v)), \qquad m<0\tag5.12$$
again defines an element of $\CH^1(\M)$.

For $m=0$, the definition of $\hat\Cal Z(0,v)$ is more speculative.
Using the canonical map from $\Pich(\M)$ to $\CH^1(\M)$,
we let
$$\hat\Cal Z(0,v) = - \bigg(\,\hat\o + (0,\log v)\,\bigg).\tag5.13$$

We now have defined elements $\hat\Cal Z(m,v)\in \CH^1(\M)$ for all $m\in\Z$
and $v\in \R^\times_+$.
We define the following two generating series, which are formal Laurent series
in a parameter $q$. Later we will take $q = e(\tau)= e^{2\pi i\tau}$, where $\tau = u+iv\in\H$.

The first generating function involves only the orbifold $\M(\C)= [\Gamma\back D]$. Let
$$\o(m,v) = dd^c\Xi(m,v) + \delta_{D_{\Cal Z(m)}}\tag5.14$$
be the right hand side of the Green's equation for $\Xi(m,v)$. Then let
$$\deg(\hat\Cal Z(m,v)) = \int_{[\Gamma\back D]} \o(m,v).\tag5.15$$
If $m>0$, then $\deg(\hat\Cal Z(m,v))$
is just the usual degree of the $0$--cycle $\Cal Z(m)_\C$ (in the stack sense).
If $m<0$, then $\deg(\hat\Cal Z(m,v))=0$, since $\Xi(m,v)$ is smooth is this case so that $\o(m,v)$ is exact.
For $m=0$, we take $\o(0,v)$ to be the Chern form of $-\hat\o$, i.e., $-\mu$, and hence
$$\deg(\hat\Cal Z(0,v)) := \int_{[\Gamma\back D]} \o(0,v) = -\,\vol(\M(\C)).\tag5.16$$
The generating function for degrees is then
$$\align
\phi_{\deg}(\tau) :&= \sum_{m} \deg(\hat\Cal Z(m,v))\,q^m\tag5.17\\
\nass
{}&= -\,\vol(\M(\C)) + \sum_{m>0} \deg(\Cal Z(m)_\C)\,q^m.
\endalign
$$

For the second generating function, we use the height pairing (4.15) of our cycles  with the
class $\hat\o\in  \Pich(\M)$, and let
$$\phi_{\text{\rm height}}(\tau) = \sum_m \langle\ \hat\Cal Z(m,v),\hat\o\ \rangle\,q^m.\tag5.18$$
At the moment, we regard $\phi_{\deg}(\tau)$ (resp. $\phi_{\text{\rm height}}(\tau)$ )
as a formal generating series, but
our main theorem will identify it as a bona fide  holomorphic (resp. non-holomorphic) function of the variable
$\tau$ by identifying it with the Fourier expansion of a special value of
an Eisenstein series (resp. of the derivative of an Eisenstein series).

\subheading{Part II. Eisenstein series}

\subheading{6. Eisenstein series of weight 3/2}

In this section, we introduce the Eisenstein series of half--integral weight which will be connected
with the arithmetic geometry discussed in Part I.  A more general discussion of
such series from an adelic point of view can be found in \cite{\annals}. The series we consider
are, of course, rather familiar from a classical point of view, and an expression for them
in this language will emerge in section 8 and 16 below. Thus, one purpose of the present section is to
explain how such classical series are associated to indefinite quaternion algebras in a
natural way, via the Weil representation. A second advantage of the adelic viewpoint is that
it allows one to assemble the Fourier coefficients out of local quantities. This construction
shows in a very
clear way the dependence of these coefficients, and more importantly, of their derivatives
on the choice of local data.

Let $G'_{\A}$  be the metaplectic extension of $Sp_1(\A) = SL_2(\A)$ by
$\Bbb C^1$, and let $P'_{\A}$ be the preimage of the subgroup $N(\A)M(\A)$ of $\SL_2(\A)$ where
$$
N(\A) = \{ n(b) = \pmatrix 1 & b \\ 0 & 1 \endpmatrix\mid\, b \in \A\}\tag6.1
$$
and
$$
M(\A) = \{ m(a) = \pmatrix a & 0 \\ 0 & a^{-1} \endpmatrix\mid \, a \in \A^\times\}.\tag6.2
$$
As in \cite{\waldspurger} we have an identification $G'_\A\simeq \SL_2(\A)\times\C^1$
where the multiplication on the right is given by
$[g_1,z_1][g_2,z_2] = [g_1g_2,c(g_1,g_2)z_1z_2]$ with cocycle $c(g_1,g_2)$
as in \cite{\waldspurger} or \cite{\gelbart}.
Let $G'_{\Bbb Q}= \SL_2(\Bbb Q)$, identified with a subgroup of
$G'_\A$ via the canonical splitting homomorphism $G'_\Q\rightarrow G'_\A$,
and let $P'_\Q= P'_\A\cap G'_\Q$.
An id\`ele character $\chi$ of  $\Bbb Q^\times \backslash \A^\times$, determines
a character $\chi^{\psi}$ of $P'_{\Bbb Q} \backslash P'_{\A}$ via
$$
\chi^\psi([n(b)m(a), z]) = z\, \chi(a)\, \gamma(a, \psi)^{-1},\tag6.3
$$
where $\psi$ is our fixed additive character of $\Q\back \A$ and
$\gamma(a, \psi)$ is the Weil index (\cite{\weilI} or \cite{\rao}, appendix). For $s\in \C$, let
$$I(s, \chi) =\hbox{Ind}_{P'_{\A}}^{G'_{\A}} \chi^\psi |\, |^s\tag6.4$$
be the principal series representation of $G'_\A$ determined by $\chi^\psi$. A section
$\P(s)\in I(s,\chi)$ is thus a smooth function on $G'_{\A}$ such that
$$
\Phi(p' g',s)
 =\chi^\psi(p')\, |a|^{s+1} \Phi(g',s).\tag6.5
$$
where $p'=[n(b)m(a),z]$. Such a section is called {\it standard} if its restriction
to the maximal compact subgroup $K'\subset G'_\A$ is independent of $s$ and {\it factorizable} if
$\P(s) = \tt_p\P_p(s)$ for the decomposition of the induced representation
$I(s,\chi) = \otimes_p' I_p(s, \chi_p)$. Here, for each prime $p$, $I_p(s,\chi_p)$ is
the corresponding induced representation of $G'_p$, the metaplectic extension
of $\SL_2(\Q_p)$.
The Eisenstein series associated to a standard section $\Phi(s) \in I(s, \chi)$ is given by
$$
E(g, s, \Phi)
 =\sum_{\gamma \in P'_{\Bbb Q}\backslash G'_{\Bbb Q}} \Phi(\gamma g, s).\tag6.6
$$
This series is absolutely convergent for $\Re(s) > 1$ and has a meromorphic continuation to the whole
complex $s$-plane. Note that this series is normalized so that it has a functional equation
$$E(g',-s,M(s)\P) = E(g',s,\P),\tag6.7$$
where $M(s): I(s,\chi) \rightarrow I(-s,\chi^{-1})$ is the intertwining operator.  It has a Fourier
expansion
$$
E(g', s, \Phi)
 =\sum_{m \in \Bbb Q} E_m(g', s, \Phi) \tag6.8
$$
where, in the half-plane of absolute convergence,
$$
E_m(g', s, \Phi)
 =\int_{\Bbb Q \backslash \A} E(n(b) g', s, \Phi)\, \psi(-mb)\, db,\tag6.9
$$
for $db$ the  self-dual measure on $\A$ with respect to $\psi$.
When $m \ne 0$ and $\P(s) = \tt_p\P_p(s)$ is factorizable, the $m$th Fourier coefficient has a product expansion
$$
E_m(g', s, \Phi)
 = \prod_{p \le \infty} W_{m, p}(g'_p, s, \Phi_p),\tag6.10
$$
where
$$
 W_{m, p}(g'_p, s, \Phi_p) =\int_{\Bbb Q_p} \Phi_p(w n(b) g'_p, s)\, \psi(-mb)\, db\tag6.11
$$
is the local Whittaker function, and
$w = \pmatrix 0 & -1 \\ 1 & 0 \endpmatrix\in G'_\Q$. Here $db$ is the self dual measure on
$\Q_p$ for $\psi_p$. On the other hand, the constant
term is
$$
E_0(g', s, \Phi) = \Phi(g', s) + \prod_{p \le \infty} W_{0, p}(g'_p, s, \Phi).\tag6.12
$$
Recall that the poles of the Eisenstein series are precisely those of its constant term.

In this paper, we will only be concerned with the case of a quadratic character $\chi$
given by  $\chi(x) = (x,\kappa)_\A$ for $\kappa\in \Q^\times$,
where $(\ ,\ )_\A$ denotes the global quadratic Hilbert symbol, so we
will omit $\chi$ from the notation and write $I(s) = \tt_p' I_p(s)$ for the induced representation, etc.
We now begin to make specific choices of the local sections $\P_p(s)$.

As before, let $B$ be an indefinite quaternion algebra over $\Q$ with a fixed
maximal order $O_B$. Once again, the case $B=M_2(\Q)$ and $O_B=M_2(\Z)$ will be allowed.
Let
$$V=\{\ x\in B\mid \tr^o(x) = 0\ \}\tag6.13$$
with quadratic form defined by $Q(x) = -x^2$, and
let $L = O_B\cap V$. Note that the determinant of the quadratic space
$(V,Q)$, i.e., $\det(S)$ where $S$ is the matrix for the quadratic form, is a square.
Therefore, the discriminant $-\det(S)$ is $-1$ and the quadratic character $\chi_V$ associated to $V$ is
given by $\chi_V(x) = (x,-1)_\A$.
We therefore take $\chi=\chi_V$ and $\kappa=-1$ in this case.

The group $G'_\A$ (resp. $G'_p$) acts on the Schwartz space $S(V(\A))$ (resp. $S(V_p)$ ) via the Weil representation $\o$
(resp. $\o_p$) determined by $\psi$ (resp. $\psi_p$).

For a finite prime $p$, let $\P_p(s)\in I_p(s)$ be the standard section extending
$\l_p(\ph_p)$, where
$$\l_p:S(V_p)\rightarrow I_p(\frac12),\qquad \l_p(\ph_p)(g') = \big(\o(g')\ph_p\big)(0)\tag6.14$$
is the usual map
and $\ph_p\in S(V_p)$ is the characteristic function of $L_p= L\tt_\Z\Z_p$.

Let $K'_\infty$  be the inverse image in
$G'_\A$ of $\text{\rm SO}(2)\subset \SL_2(\R)$.
For $\ell\in\frac12\Z$, there is a character $\nu_\ell$ of $K'_\infty$ such that
$$\nu_\ell([k_\theta,1])^2 = e^{2i\theta\ell}.\tag6.15$$
For $\ell\in \frac{3}2+2\Z$,
there is a unique standard section $\P_\infty^\ell(s)\in I_\infty(s)$
with
$$\P_\infty^\ell(k,s) = \nu_\ell(k),\tag6.16$$
for $k\in K'_\infty$.

Let
$$\P^{\ell,D(B)}(s) = \P_\infty^{\ell}(s)\tt\bigg(\tt_p \P_p(s)\bigg)\tag6.17$$
be the associated global standard section. A little more generally,
for a finite prime $p$, let $\P_p^+(s)$ be the standard section \
arising from the maximal order $M_2(\Z_p)$ in $M_2(\Q_p)$ and
let $\P_p^-(s)$ be the standard section arising from the maximal order in the
division quaternion algebra over $\Q_p$. Then for any square free positive integer $D$, we have
a global section
$$\P^{\ell,D}(s) = \P_\infty^{\ell}(s)\tt\bigg(\tt_p \P_p^{\e(D)}(s)\bigg),\tag6.18$$
where $\e(D) = (-1)^{\ord_p(D)}$.

Since, by strong approximation,
$G'_\A= G'_\Q G'_\R K^0$
for any open subgroup $K^0$ of $G'_{\A_f}$, we loose no information
by restricting automorphic forms to $G'_\R$, the inverse image of $\SL_2(\R)$
in $G'_\A$.
For $\tau=u+iv\in \H$, let
$$g'_\tau = [\,n(u)m(v^{\frac12}),1\,]\in G'_\infty \subset G'_\A.\tag6.19$$
Then, if $\P(s)$ is a standard factorizable section with $\P_\infty(s) =\P^\ell_\infty(s)$,  we set
$$E(\tau,s,\P) = v^{-\frac{\ell}2}\,E(g'_\tau,s,\P),\tag6.20$$
and, by (6.10), we have
$$E_m(\tau,s,\P) = v^{-\frac{\ell}2}\,W_{m,\infty}(g'_\tau,s,\P^\ell_\infty)\cdot \prod_p W_{m,p}(s,\P_p),\tag6.21$$
for $m\ne0$, and
$$E_0(\tau,s,\P) = v^{\frac12(s+1-\ell)}\cdot \P_f(e) + v^{-\frac{\ell}2}\,W_{0,\infty}(g'_\tau,s,\P^\ell_\infty)\prod_p W_{0,p}(s,\P_p).\tag6.22$$

The main series of interest to us will be
$E(\tau,s,\P^{\ell,D})$, associated to the standard section $\P^{\ell,D}(s)$ of (6.18).
This series has weight $\ell$, where $\ell = \frac32$, $\frac72$, $\frac{11}2, \dots$.
Note that the character $\chi$ is given by $\chi(x) = (x,-1)_\A$ in this case.
A second family $E(\tau,s,\P^{\ell,D})$, with $\ell = \frac12$, $\frac52$, $\frac92, \dots$ etc. is obtained by the same construction
applied to the quadratic space $(V,Q_-)$ where $Q_-(x) = x^2$. In this case, $\kappa=1$ and
$\chi$ is trivial. These cases will be discussed in more detail
in \cite{\ky}. In the present paper, we will only be concerned with the case $\ell=\frac32$,
and so, from now on, we take $\kappa=-1$.

In the next section, we will give a geometric interpretation of the first two terms of the
Laurent expansion of the series $E(\tau,s,\P^{\frac32,D(B)})$ at the point $s=\frac12$. For this it will
be convenient to normalize the series as follows.
For any square free positive integer $D$, let
$$
\E(\tau,s;D) := (s+\frac12)\,c(D)\,\L_D(2s+1)\,E(\tau,s,\P^{\frac32,D}).\tag6.23
$$
where
$$\L_D(2s+1) = \left(\frac{D}{\pi}\right)^{s+\frac12}\,\Gamma(s+\frac12)\,\zeta(2s+1)\cdot \prod_{p|D} (1-p^{-2s-1}),\tag6.24$$
and
$$c(D) = - (-1)^{\ord(D)}\,\frac{1}{2\pi}\,D \prod_{p|D} (p+1)^{-1},\tag6.25$$
where $\ord(D) = \sum_p \ord_p(D)$.
Note that at the point $s=\frac12$, of interest to us, the normalizing factor has value
$$
c(D)\,\L_D(2) = -(-1)^{\ord(D)}\,\frac1{12}\prod_{p|D}(p-1).\tag6.26
$$
Then, in the case $D=D(B)$, and recalling (2.7),
$$c(D)\,\L_D(2) = - \,\vol(\M(\C)).\tag6.27$$
This expression explains the choice of $c(D)$.  In addition, the normalized Eisenstein
series satisfies the simple functional equation, c.f. section 16,
$$\E(\tau,s;D) = \E(\tau,-s;D).\tag6.28$$

Finally, we restrict to the case $D=D(B) >1$ and introduce the {\it modified}
Eisenstein series
$$\Cal E(\tau, s; D(B)) := \E(\tau,s;D(B)) + \sum_{p|D} c_p(s)\,\E(\tau,s;D(B)/p),\tag6.29$$
where $c_p(s)$ is any rational function of $p^{-s}$ satisfying
$$c_p(\frac12)=0,\qquad\qquad\text{and}\qquad\qquad c'_p(\frac12)= -\frac{p-1}{p+1}\,\log(p).\tag6.30$$
To retain the functional equation (6.28) one should also require that $c_p(s)=c_p(-s)$, although
we will not use this.
The motivation for the definition of $\Cal E(\tau,s,D(B))$ comes from geometric
considerations which will emerge below. Note that
$$\Cal E(\tau,\frac12;D(B)) = \E(\tau,\frac12;D(B)),\tag6.31$$
and
$$\Cal E'(\tau,\frac12;D(B)) = \E'(\tau,\frac12;D(B)) + \sum_{p|D} c'_p(\frac12)\cdot\E(\tau,\frac12;D(B)/p).\tag6.32$$

\subheading{7. The main identities }

In this section, we state our main results on the generating functions
$$\phi_{\deg}(\tau) = -\vol(\M(\C)) + \sum_{m>0} \deg(\Cal Z(m)_\C)\,q^m\tag7.1$$
and
$$\phi_{\text{\rm height}}(\tau) = \sum_m \langle\ \hat\Cal Z(m,v),\hat\o\ \rangle\,q^m.\tag7.2$$
introduced in
section 5.

The following result is actually well known, cf., for example, \cite{\funkecompo}.
We state it here to bring out the analogy with Theorem~7.2.

\proclaim{Proposition 7.1} For any indefinite quaternion division algebra $B$ over $\Q$
with associated moduli stack $\M$, as in section 1--5 above, the generating function for the
degrees of the special cycles coincides with the value at $s=\frac12$ of the
Eisenstein series $\Cal E(\tau,s;D(B))$ of weight $\frac32$:
$$
\phi_{\deg}(\tau) = \Cal E(\tau, \frac12;D(B)).
$$
\endproclaim

\proclaim{Theorem 7.2} Under the same assumptions, the generating series for heights
of the special cycles coincides, up to an additive constant, with the derivative at $s=\frac12$
of the
Eisenstein series $\Cal E(\tau,s;D(B))$ of weight $\frac32$:
$$\phi_{\text{\rm height}}(\tau) = \Cal E'(\tau, \frac12;D(B)) +\bold c.$$
for some constant $\bold c$.
\endproclaim

These identities are to be understood as follows. We write the Fourier expansion of
the modified Eisenstein series as
$$\Cal E(\tau,s;D(B)) = \sum_{m} A_m(s,v)\,q^m,\tag7.3$$
so that the Fourier expansions of the value and derivative at $s=\frac12$ are
$$\Cal E(\tau,\frac12;D(B)) = \sum_{m} A_m(\frac12,v)\,q^m,\tag7.4$$
and
$$\Cal E'(\tau,\frac12;D(B)) = \sum_{m}A_m'(\frac12,v)\,q^m.\tag7.5$$
Proposition~7.1 then says that
$$A_m(\frac12,v) = \cases
\deg(\Cal Z(m)_\C) &\text{ if $m>0$,}\\
\nass
-\,\vol(\M(\C))&\text{ if $m=0$,}\\
\nass
0&\text{ if $m<0$.}
\endcases\tag7.6
$$
Analogously, Theorem~7.2 says that, for $m\ne0$,
$$A_m'(\frac12,v) =
\lan \hat\Cal Z(m,v),\hat\o\ran.\tag7.7
$$
The ambiguity $\bold c$ in Theorem~7.2 thus arises from the fact that we do not have an explicit expression
for the quantity $\lan\hat\o,\hat\o\ran$.
If we knew, a priori, that $\phi_{\text{\rm height}}(\tau)$
was (the Fourier expansion of) a modular form of weight $\frac32$, then we could
conclude that $\bold c=0$ and, by formula (5.11) of Bost \cite{\bost}, that
$$A'_0(\frac12,v)\  \overset{??}\to{=}\  \lan \hat\Cal Z(0,v),\hat\o\ran = -\lan\hat\o,\hat\o\ran -\frac12\log(v)\deg(\hat\o).\tag7.8$$
For further discussion of this point, see section 13.

As already explained in the introduction, Theorem~7.2 is proved by an explicit computation of both sides of (7.7).
For the left hand side, this will be done in the next section. The right hand side will be computed in sections 9--12.

\subheading{8. Fourier expansions and derivatives}

In this section, we describe the first two terms of the Laurent expansion of the
Eisenstein series $\Cal E(\tau, s;D(B))$ at the point $s = \frac12$. By (6.21) and (6.22),
the essential point is to describe the behaviour of the local Whittaker functions $W_{m,p}(s,\P_p^D)$
and
$$W_{m,\infty}(\tau,s,\P_\infty^{\frac32}) := v^{-\frac34}\,W_{m,\infty}(g'_\tau,s,\P^\frac32_\infty).\tag8.1$$
The calculations of this section will be elementary manipulations based on results about these Whittaker
functions quoted from Part IV below.

In what follows, for a nonzero integer $m$, we write
$$4 m = n^2 d\tag8.2$$
where $-d$ is a fundamental discriminant, i.e., discriminant of the field $\kay =\kay_d=\Q(\sqrt{-m})$.
Note that if $4m=-n^2$, then $\kay=\Q\oplus\Q$. Let $\chi_d$ be the corresponding Dirichlet character, so that
$$\chi_d(p) = \cases 1&\text{ if $p$ is split in $\kay_d$,}\\
-1&\text{ if $p$ is inert in $\kay_d$,}\\
0&\text{ if $p$ is ramified in $\kay_d$.}
\endcases\tag8.3
$$

For a given $m$ and for a square free positive integer $D$, define a modification of the
standard Dirichlet L-series $L(s,\chi_d)$ by
$$
L(s, \chi_m; D)  :=  L(s, \chi_d) \prod_{p| n D} b_p(n, s; D) \tag8.4
$$
where $b_p(n, s; D)$ is defined as follows. Set
$$k = k_p(n)= \ord_p(n)\tag8.5$$
and
$X=p^{-s}$.  Then for $p\nmid D$
$$
b_p(n, s; D) =\frac{1-\chi_d(p)\, X  + \chi_d(p)\, p^k X^{(1+2k)} - (pX^2)^{k+1}}{1-p X^2},\tag8.6
$$
and, for $p\mid D$,
$$
\align
&b_p(n, s; D)\tag8.7
\\
 &=\frac{(1-\chi\, X)(1-p^2 X^2) - \chi\, p^{k+1}X^{2k+1} + p^{k+2}X^{2k+2}
          +\chi\, p^{k+1}X^{2k+3} - p^{2k+2}X^{2k+4}}
       {1- p X^2}.
\endalign
$$
Here, for a moment, we write $\chi$ for $\chi_d(p)$.
Depending on whether or not $p\mid d$ we can rewrite (8.7) as
$$
b_p(n, s; D) =\frac{  1- p^2 X^2 + p^{k+2} X^{2k+2}(1- X^2)}{1-p X^2}\qquad
\text{ if $ p \mid  d$ and $p\mid D$,}\tag8.8
$$
and
$$
b_p(n, s; D) =\frac{ (1-\chi\, X)(1-p^2 X^2)-\chi\, p^{k+1}X^{2k+1}(1-\chi\, p X)(1-X^2)}{1-pX^2}\qquad
\text{ if $p \nmid d$ and  $p \mid D$.}\tag8.9
$$
In all cases the local factor $b_p(n, s; D)$ is, in fact, a polynomial in $X= p^{-s}$
and is, hence, entire in $s$.
It satisfies the functional equation
$$
|nD|_p^{-s} b_p(n, s; D) = |nD|_p^{s-1} b_p(n, 1-s; D).\tag8.10
$$

One of the main results of section 14 is the following.
\proclaim{Proposition 8.1} For a fixed prime $p$,\hfill\break
(i) if $m\ne 0$, then
$$W_{m,p}(s+\frac12,\P^D_p) = L_p(s+1,\chi_d)\,b_p(n,s+1;D)\cdot\cases
C^+_p\cdot \frac{1}{\zeta_p(2s+2)}&\text{ for $p\nmid D$,}\\
\nass
C^-_p&\text{ if $p\mid D$.}
\endcases
$$
where the constants $C_p^\pm$ are given by
$$
C_p^+=\cases 1 &\text{ if $ p \ne 2$,}\\
\frac{1}{\sqrt 2} \zeta_8^{-1} &\text{ if $p=2$,}
\endcases
$$
and $C_p^-= -p^{-1} C_p^+$. Here $\zeta_8 = e(\frac18)$.
\hfill\break
(ii) If $m=0$, then
$$W_{0,p}(s+\frac12,\P^D_p) = \zeta_p(2s)\cdot\cases
C^+_p\cdot \frac{1}{\zeta_p(2s+1)}&\text{ for $p\nmid D$,}\\
\nass
C^-_p\cdot\frac{1}{\zeta_p(2s-1)}&\text{ if $p\mid D$.}
\endcases
$$
\endproclaim

From (6.21), (6.22), and these formulas, we obtain a nice description of
the Fourier expansion of $\E(\tau,s;D)$. \proclaim{Corollary 8.2} Let
$C_f(D) = \prod_p C_p^{\e(D)}$.\hfill\break (i) For $m\ne0$,
$$E_m(\tau,s;D) = C_f(D)\cdot
W_{m,\infty}(\tau,s,\P^{\frac32}_\infty)\cdot
\frac{L(s+\frac12,\chi_d)}{\zeta_{D}(2s+1)} \cdot (nD)^{-2s}\,\prod_p
b_p(n,\frac12-s;D),$$ and $$\align \E_m(\tau,s;D) &=
c(D)\,C_f(D)\,\left(\frac{D}{\pi}\right)^{s+\frac12}
\Gamma(s+\frac32)\cdot W_{m,\infty}(\tau,s,\P^{\frac32}_\infty)\\ \nass
&\qquad\qquad\qquad\times L(s+\frac12,\chi_d) \cdot (nD)^{-2s}\,\prod_p
b_p(n,\frac12-s;D).
\endalign$$
(ii) For $m=0$,
$$E_0(\tau,s;D) = v^{\frac12(s-\frac12)} + W_{0,\infty}(\tau,s,\P^{\frac32}_\infty)\,C_f(D)\cdot
\frac{\zeta(2s)}{\zeta_D(2s+1)}\cdot\prod_{p\mid D} \frac{1}{\zeta_p(2s-1)},$$
and
$$\align
\E_0(\tau,s;D) &= v^{\frac12(s-\frac12)}\,(s+\frac12)\,c(D)\,\L_D(2s+1)\\
\nass
&{} + W_{0,\infty}(\tau,s,\P^{\frac32}_\infty)\,c(D)\,C_f(D)\,\left(\frac{D}{\pi}\right)^{s+\frac12}
\Gamma(s+\frac32)\cdot
\zeta(2s)\cdot\prod_{p\mid D} \frac{1}{\zeta_p(2s-1)},
\endalign$$
\endproclaim
Here
$$c(D)\,C_f(D) = - \frac{1}{\sqrt 2}\, \zeta_8^{-1}\,\frac{1}{2\pi} \prod_{p|D} (p+1)^{-1}.\tag8.11$$

Using Corollary~8.2, we now compute the value of $\E(\tau,s;D)$ at $s=\frac12$.
We start with the constant term.

The following result is a special case of (iii) of Proposition~15.1 below.
\proclaim{Lemma 8.3}
$$W_{0,\infty}(\tau,s,\P_\infty^{\frac32}) = 2\pi\,(-i)^{\frac32}\,v^{-\frac12(s+\frac12)}\,2^{-s}\,\frac{\Gamma(s)}{\Gamma(\a)
\Gamma(\b)},$$
for $\a=\frac12(s+\frac52)$ and $\b=\frac12(s-\frac12)$. Here $(-i)^{\frac32}=e(-\frac38)$.
\endproclaim

Since the zero of $\Gamma(\beta)^{-1}$ at $s=\frac12$ cancels the pole of $\zeta(2s)$ there,  the second term in $\E_0(\tau,s;D)$
has a zero of order equal to the number of primes dividing $D$, and we obtain
\proclaim{Corollary 8.4} For $D>1$, the constant term at $s=\frac12$ is
$$\E_0(\tau,\frac12;D) = c(D)\,\L_D(2) = -(-1)^{\ord(D)}\,\frac1{12}\prod_{p|D}(p-1) = \zeta_D(-1).$$
\endproclaim

Next we consider the coefficients of $\E(\tau,\frac12;D)$ for $m\ne 0$.

If $p\nmid D$, then by (8.6),
$$
b_p(n, 0; D) =\frac{1-\chi_d(p)\, + \chi_d(p)\, p^k - p^{k+1}}{1-p}.\tag8.13
$$
Note that, when $\chi_d(p)=1$, this simplifies to $p^k = |n|_p^{-1}$.

If $p\mid D$, then by (8.7)--(8.9),
$$
b_p(n, 0; D)= (1-\chi_d(p)) (1+p).\tag8.14
$$
Note that this quantity is actually independent of $n$, and that
$b_p(n, 0; D) =0$ if and only if $p \mid D$ and $\chi_d(p)=1$.

The proof of the following identity is a simple combinatorial
exercise, which we omit.
\proclaim{Lemma 8.5} (i) For $p\nmid D$, and $k=\ord_p(n)$,
$$b_p(n, 0; D) =\sum_{c | p^{k}} c \prod_{\ell|c} (1- \chi_d(\ell) \ell^{-1}),$$
where $\ell$ runs over the prime factors of $c$ and the product is taken to be $1$ when $c=1$.
\hfill\break
(ii)
$$\prod_{p\nmid D} b_p(n,0;D) = \sum_{ c | n \atop (c,D)=1 } c \prod_{\ell|c} (1- \chi_d(\ell) \ell^{-1}).
$$
Here, again, $\ell$ runs over the prime factors of $c$ and the product is taken to be $1$ when $c=1$.
\qed
\endproclaim

On the other hand, the following fact is a special case of (iv) of Proposition~15.1 below.
\proclaim{Lemma 8.6}
$$W_{m,\infty}(\tau,\frac12,\P^{\frac32}_\infty) =
\cases 2\,C_\infty\cdot m^{\frac12}\,q^m&\text{ if $m>0$, and}\\
\nass
0&\text{ otherwise.}
\endcases
$$
where $C_\infty = (-2i)^{\frac32}\,\pi$.
\endproclaim

Combining these facts, we obtain the following results.\hfill\break
For $m<0$, the vanishing of the archimedean factor yields:
$$\E_m(\tau,\frac12;D)=0,\qquad \text{when $\chi_d\ne1$, or $D>1$.}
\tag8.16
$$
For $m>0$, (i) of Corollary~8.2 gives
$$\align
\E_m(\tau,\frac12;D) &= c(D)\,C_f(D)\,C_\infty\cdot\frac{D}{\pi}\cdot 2\,m^{\frac12}\,q^m\cdot L(1,\chi_d)\,(nD)^{-1}
\prod_p b_p(n,0;D)\tag8.17\\
\nass
{}&=c(D)\,C_f(D)\,C_\infty\cdot q^m\cdot 2\, \frac{h(d)}{w(d)}\cdot
\bigg( \sum_{ c | n \atop (c,D)=1 } c \prod_{\ell|c} (1- \chi_d(\ell) \ell^{-1})\bigg)\\
\nass
&\qquad\qquad\qquad\times \bigg(\prod_{p|D} (1-\chi_d(p)) (1+p)\bigg)\\
\nass
\nass
{}&= q^m\cdot 2\,\frac{h(d)}{w(d)}\cdot
\bigg( \sum_{ c | n \atop (c,D)=1 } c \prod_{\ell|c} (1- \chi_d(\ell) \ell^{-1})\bigg)\cdot\bigg(\prod_{p|D} (1-\chi_d(p))\bigg).
\endalign
$$
Here, $w(d) = |O_{\smallkay}^\times|$ is the number of roots of unity in the maximal order $O_{\smallkay}$  of
$\kay_d$, $h(d)$ is the class number, and
$$c(D)\,C_f(D)\,C_\infty =  \prod_{p|D} (p+1)^{-1}.\tag8.18$$

For $m>0$, let
$$
H_0(m;D) = \frac{h(d)}{w(d)}\cdot
\bigg( \sum_{ c | n \atop (c,D)=1 } c \prod_{\ell|c} (1- \chi_d(\ell) \ell^{-1})\bigg),\tag8.19
$$
where, as before, in the product, $\ell$ runs over the prime factors of $c$ and the product is taken to be $1$
when $c=1$,
and
$$\delta(d;D) = \prod_{p|D} (1-\chi_d(p)).\tag8.20$$
Thus, we obtain the Fourier expansion of $\Cal E(\tau,\frac12;D)$.
\proclaim{Proposition 8.7}
For $D>1$,
$$\Cal E(\tau,\frac12;D) = c(D)\,\L_D(2)+ \sum_{m>0} 2\,\delta(d;D)\,H_0(m;D)\, q^m.$$
Here
$$c(D)\,\L_D(2)= -(-1)^{\ord(D)}\,\frac1{12}\prod_{p|D}(p-1).$$
\endproclaim

This Eisenstein series of weight $\frac32$ is a familiar object.
Recall that,
if $O_{c^2d}$ is
the order in $O_{\smallkay}$ of conductor $c$, with class number $h(c^2d)$ and with $w(c^2d) = |O_{c^2d}^\times|$, then
\cite{\BorevichShafarevic}, p.250\footnote{The quantity $e_c$ there is $|\Cal O_d^\times:\Cal O_{c^2d}^\times| = w(d)/w(c^2d)$},
$$\frac{h(c^2d)}{w(c^2d)} =\frac{h(d)}{w(d)}\cdot c\, \prod_{\ell|c} (1- \chi_d(\ell) \ell^{-1}).\tag8.21$$
Thus,
$$
H_0(m;D)=\sum_{c|n\atop (c,D)=1} \frac{h(c^2d)}{w(c^2d)}.\tag8.22
$$
For example, if $D=1$, i.e., in the case of $B=M_2(\Q)$,
$$H_0(m;1) = \sum_{c|n} \frac{h(c^2d)}{w(c^2d)}\tag8.23$$
is quite close to\footnote{Precisely, $2\,H_0(m;1) = H(4m)$.} the `class number' $H(m)$ which appears in
the Fourier expansion
$$\Cal F(\tau)= -\frac1{12} + \sum_{m>0} H(m)\,q^m + \sum_{n\in \Z} \frac1{16\pi}\, v^{-\frac12}
\int_1^\infty e^{-4\pi n^2 v r}\, r^{-\frac32}\,dr\, q^{-n^2},\tag8.24$$
of Zagier's nonholomorphic Eisenstein series of weight $\frac32$, \cite{\cohen}, \cite{\zagier}. In fact,
when $D=1$, we have
$$\Cal E(\tau,\frac12;1)= -\frac1{12} + \sum_{m>0} 2\,H_0(m;1)\,q^m+
\sum_{n\in\Z}\frac{1}{8\pi}\,v^{-\frac12}\,\int_1^\infty e^{-4\pi n^2 vr}\,r^{-\frac32}\,dr \cdot q^{-n^2}.\tag8.25$$
This case will be discussed in detail in the sequel \cite{\kryIII}.

Next we consider the derivative $\Cal E'(\tau,\frac12;D)$ in the case $D= D(B)>1$.
In this case, the only terms which contribute are the following:
\roster
\item"{(i)}" $m>0$ and $\delta(d;D)\ne0$,
\item"{(ii)}" $m>0$ and there is a unique $p\mid D(B)$ such that $\chi_d(p)=1$,
\item"{(iii)}" $m<0$ and $\delta(d;D)\ne0$, and
\item"{(iv)}" $m=0$.
\endroster
In cases (i) and (iv),  $\Cal E_m(\tau,\frac12;D)\ne0$. In cases (ii) and (iii),
$\Cal E_m(\tau,s;D)$ has a simple zero at $s=\frac12$ due to the vanishing of the local
factor $b_p(n,0;D)$ in case (ii) and the archimedean factor $W_{m,\infty}(\tau,\frac12;\P_\infty^{\frac32})$
in case (iii). In all other cases, $\Cal E_m(\tau,\frac12;D)$ has a zero of order at least $2$ at $s=\frac12$.

\proclaim{Theorem 8.8} Assume that $D=D(B)>1$. \hfill\break
(i) If $m>0$ and there is no prime $p\mid D$ for which $\chi_d(p)=1$, then
$$\align
&\Cal E'_m(\tau,\frac12;D)\\
\nass
{}&= 2\,\delta(d;D)\,H_0(m;D)\cdot q^m\cdot\,\bigg[ \,
\frac12\, \log(d)
 + \frac{L'(1,\chi_d)}{L(1,\chi_d)} -\frac12\log(\pi) -\frac12\gamma\\
\nass
\nass
{}&\hbox to .7in{}
+\frac12 J(4\pi mv)
+\sum_{p \atop p\nmid D} \bigg(\ \log|n|_p - \frac{b_p'(n,0;D)}{b_p(n,0;D)}\ \bigg)
+\sum_{p \atop p\mid D} K_p\,\log(p)
\ \bigg].\\
\endalign
$$
Here
$$K_p = \cases -k + \frac{(p+1)(p^k-1)}{2(p-1)}&\text{ if $\chi_d(p)=-1$, and}\\
\nass
-1-k + \frac{p^{k+1}-1}{p-1} &\text{ if $\chi_d(p)=0$,}
\endcases
$$
with $k=\ord_p(n)$, and
$$J(t) = \int_0^\infty e^{- t r}\big[\, (1+r)^{\frac12} - 1\,\big]\,r^{-1}\,dr.$$
An explicit expression for the logarithmic derivative of $b_p(n,s;D)$ is given by (i) of Lemma~8.10, and
$H_0(m;D)$ and $\delta(d;D)$ are given by (8.19) and (8.20) respectively. \hfill\break
(ii) If there is a unique prime $p\mid D$ such that $\chi_d(p)=1$, then
$$\Cal E'_m(\tau,\frac12;D) = 2\,\delta(d;D/p)\,H_0(m;D)\cdot (p^k-1)\,\log(p)\cdot  q^m.$$
(iii) If $m<0$, then
$$
\Cal E'_m(\tau,\frac12;D) =2\,\delta(d;D)\,H_0(m;D)\cdot q^m
\cdot
\frac{1}{4\pi}\,|m|^{-\frac12}\,v^{-\frac12}\,\int_1^\infty e^{-4\pi |m|v r} r^{-\frac32}\,dr,
$$
where, for $m<0$, $H_0(m;D)$ is defined by (8.32) below.\hfill\break
(iv)
$$\Cal E'_0(\tau,\frac12;D) =
c(D)\,\L_D(2)\,\bigg[\frac12\,\log(v) - 2\frac{\zeta'(-1)}{\zeta(-1)} -1 + 2C+ \sum_{p\mid D}\frac{p\log(p)}{p-1}\,\bigg].$$
(v) All other Fourier coefficients of $\Cal E'(\tau,\frac12;D)$ vanish.
\endproclaim

\demo{Proof} We begin with the Eisenstein series $\E(\tau,s;D)$ for any $D>1$.

First consider case (i), so that $m>0$ and that there are no primes $p\mid D$ with $\chi_d(p)=1$.
Then, using (i) of Corollary~8.2, we have
$$\align
\E'_m(\tau,\frac12;D)=\E_m(\tau,\frac12;D)\,\bigg[ \, \log(D)-&\log(\pi) +1-\gamma +
\frac{W'_{m,\infty}(\tau,\frac12,\P^{\frac32})}{W_{m,\infty}(\tau,\frac12,\P^{\frac32})}\\
\nass
&{}
+ \frac{L'(1,\chi_d)}{L(1,\chi_d)} -2 \log(nD)
-\sum_p \frac{b_p'(n,0;D)}{b_p(n,0;D)}\, \bigg]
\endalign
$$

The following fact is proved in section 15.

\proclaim{Lemma 8.9}
For $m>0$,
$$
\frac{W'_{m,\infty}(\tau,\frac12,\P^{\frac32})}{W_{m,\infty}(\tau,\frac12,\P^{\frac32})}
 =
\frac12 \bigg[\, \log(\pi m) -\frac{\Gamma'(\frac32)}{\Gamma(\frac32)} + J(4\pi mv)\,\bigg].
$$
\endproclaim

Using this result and the fact that
$$\frac{\Gamma'(\frac32)}{\Gamma(\frac32)} = 2-\gamma-2\log(2),$$
and recalling that $4m=n^2d$,
we obtain
$$\align
&\E'_m(\tau,\frac12;D)\tag8.27\\
\nass
{}= \ &\E_m(\tau,\frac12;D)\,\bigg[ \,
\frac12\, \log(d)
 + \frac{L'(1,\chi_d)}{L(1,\chi_d)} -\frac12\log(\pi) -\frac12\gamma+ \frac12 J(4\pi mv)\\
\nass
{}&\hbox to 2in{} \quad
 +\sum_{p \atop p\nmid D} \bigg(\ \log|n|_p - \frac{b_p'(n,0;D)}{b_p(n,0;D)}\ \bigg)\\
\nass
{}&\hbox to 1.5in{}
-\log(D)  +\sum_{p \atop p\mid D} \bigg(\ \log|n|_p-\frac{b_p'(n,0;D)}{b_p(n,0;D)}\ \bigg)
\ \bigg].\\
\endalign
$$
Next we note the explicit expressions for the logarithmic derivatives of the $b_p(n,s;D)$'s
which will be useful later.

\proclaim{Lemma 8.10}
(i) For a prime $p\nmid D$,
$$
\align
\frac{1}{ \log p}\cdot\frac{b_p'(n, 0; D)}{b_p(n, 0; D)}
&=\frac{\chi_d(p) -\chi_d(p) \,(2k+1) p^{k} +(2k+2) p^{k+1}}
       {1-\chi_d(p)+\chi_d(p)\, p^{k} -p^{k+1}}
  -\frac{2p}{1-p}
\\
\nass
\nass
&=\cases
 \frac{p^k-1}{p^k(p-1)}-2k &\text{ if  $\chi_d(p)=1$,}
\\
\nass
 -\frac{2p(1-(k+1)p^k+ kp^{k+1})}{(p-1)(p^{k+1}-1)} &\text{ if  $\chi_d(p)=0$,}
\\
\nass
 -\frac{1+3p -(2k+1)p^k-3p^{k+1} +2 k p^{k+2}}{(p-1)(p^{k+1}+p^k-2)}
       &\text{ if  $\chi_d(p)=-1$.}
  \endcases
\endalign
$$
(ii)
For a prime $p\mid D$ with $\chi_d(p)\ne1$,
$$
\frac{1}{\log p}\cdot\frac{b_p'(n, 0; D)}{b_p(n, 0; D)}
= \cases
 -\frac{2p(p^{k+1}-1)}{p^2-1}  &\text{ if $\chi_d(p)=0$,}\\
\nass
 -\frac{2(1+p)p^{k+1} + p^2 -4p -1}{2(p^2-1)} &\text{ if $\chi_d(p)=-1$.}
 \endcases
$$
Here $k=\ord_p(n)$.
\endproclaim

In case (ii),  $m>0$ and there is a unique prime $p\mid D$ such that
$\chi_d(p)=1$.
In
this case, it is easy to verify
$$
b_p'(n, 0; D) = (1+p - 2 p^{k+1})\log p.\tag8.28
$$
Then, with the notation introduced above and using (8.14), we have
$$
\E'_m(\tau,\frac12;D) = - 2\,\delta(d;D/p)\,H_0(m;D)\cdot  q^m\cdot (p+1)^{-1}\cdot
(1+p - 2 p^{k+1})\log p.\tag8.29
$$
Recall that $k=\ord_p(n)$.

Finally, in case (iii), we need another result to be proved in section 15.
\proclaim{Lemma 8.11}
For $m<0$:
$$\align
W'_{m,\infty}(\tau,\frac12,\P^{\frac32}_\infty) &= C_\infty\,|m|^{\frac12}\,  q^m \,e^{-4\pi |m| v}
\,\int_0^\infty e^{-4\pi |m| v r} (r+1)^{-1}\,r^{\frac12}\,dr\\
\nass
{}&=C_\infty\cdot \frac14\,q^m \, v^{-\frac12}\,
\,\int_1^\infty e^{-4\pi |m|vr} \,r^{-\frac32}\,dr.
\endalign
$$
\endproclaim
Using the second expression of this Lemma and (i) of Corollary~8.2,
we have, for $m<0$,
$$\align
\E'_m(\tau,\frac12;D) &= c(D)\,C_f(D)\cdot \frac{D}{\pi}\cdot W'_{m,\infty}(\tau,\frac12,\P^{\frac32})\cdot L(1,\chi_d)\cdot (nD)^{-1}
\prod_p b_p(n,0;D)\tag8.30\\
\nass
{}&=c(D)\,C_f(D)\,C_\infty\cdot \frac{4\,h(d)\,\log|\e(d)|}{w(d)\,|d|^{\frac12}}\cdot n^{-1}\,
\bigg( \sum_{ c | n \atop (c,D)=1 } c \prod_{\ell|c} (1- \chi_d(\ell) \ell^{-1})\bigg)\\
\nass
&\qquad\qquad\qquad\times \bigg(\prod_{p|D} (1-\chi_d(p)) (1+p)\bigg)\\
\nass
{}&\qquad\qquad\qquad\qquad\times q^m\,\frac1{4\pi}\,v^{-\frac12}\,\int_1^\infty e^{-4\pi|m|vr}\,r^{-\frac32}\,dr.
\endalign
$$
where $\e(d)$ is the fundamental unit of the real quadratic field $\kay_d=\Q(\sqrt{|d|})$.
Using the value (8.18), this can rewritten as
$$
\E'_m(\tau,\frac12;D) = 2\,H_0(m;D)\,\delta(d;D)
\cdot q^m\cdot
\frac{1}{4\pi}\,|m|^{-\frac12}\,v^{-\frac12}\,\int_1^\infty e^{-4\pi |m|v r} r^{-\frac32}\,dr\tag8.31
$$
where
$$\align
H_0(m;D) &= \frac{h(d)\,\log|\e(d)|}{w(d)}\cdot
\bigg( \sum_{ c | n \atop (c,D)=1 } c \prod_{\ell|c} (1- \chi_d(\ell) \ell^{-1})\bigg)\tag8.32\\
\nass
\nass
{}&=\sum_{c | n \atop (c,D)=1} h(c^2d)\cdot \frac{\log|\e(c^2d)|}{w(c^2d)}.
\endalign
$$
is the analogue of (8.19) and (8.22) in the case of a
real quadratic field, i.e., for $m<0$.

Finally, we consider the constant term using (ii) of Corollary~8.2 and
noting that for $D=D(B)$, the second term there
has a zero of order at least $2$. This gives
$$\align
\E'_0(\tau,\frac12;D) &= c(D)\,\L_D(2)\,\bigg[\frac12\,\log(v) + 1+ 2\frac{\L'_D(2)}{\L_D(2)}\,\bigg]\tag8.33\\
\nass
{}&=c(D)\,\L_D(2)\,\bigg[\frac12\,\log(v) + 1+ \log(D)-\log(\pi)-\gamma+2\frac{\zeta'(2)}{\zeta(2)}
+2\sum_{p\mid D} \frac{\log(p)}{p^2-1}\,\bigg].
\endalign$$

Now we return to the modified Eisenstein series
$$\Cal E(\tau, s ;D) = \E(\tau,s;D) + \sum_{p|D} c_p(s)\,\E(\tau,s;D/p)$$
of (6.29) for $D=D(B)>1$.
By (6.32), the Fourier coefficients of $\Cal E'(\tau,\frac12;D)$ for $m<0$ agree with those
of $\E'(\tau,\frac12;D)$, so that (8.31) gives part (iii) of Theorem~8.8.

If $m>0$ and for all $p\mid D$, $\chi_d(p)\ne1$, note that by (8.13) and Lemma~8.5,
$$\align
\E_m(\tau,\frac12;D/p) &= 2\,\delta(d;D/p)\,H_0(m;D/p)\,q^m\tag8.34\\
\nass
{}& = (1-\chi_d(p))^{-1}\cdot
\frac{1-\chi_d(p)+\chi_d(p)p^k- p^{k+1}}{1-p}\cdot 2\,\delta(d;D)\,  H_0(m;D)\cdot q^m\\
\nass
{}&= (1-\chi_d(p))^{-1}\cdot
\frac{1-\chi_d(p)+\chi_d(p)p^k- p^{k+1}}{1-p}\cdot \E_m(\tau,\frac12;D).
\endalign$$
Therefore,
$$\Cal E'_m(\tau,\frac12;D) = \E_m(\tau,\frac12;D)\,\bigg[ \dots +
\sum_{p|D} c'_p(\frac12)\cdot (1-\chi_d(p))^{-1}\cdot
\frac{1-\chi_d(p)+\chi_d(p)p^k- p^{k+1}}{1-p} \ \bigg]\tag8.35$$
where the dots indicate the expression in (8.27) for $\E'_m(\tau,\frac12;D)$.
For a prime $p\mid D$, we write $c'$ for $c'_p(\frac12)/\log(p)$ and $k=\ord_p(n)$.
Then, the coefficient of $\log(p)$ inside the bracket is
$$-1-k -\frac1{\log(p)}\cdot \frac{b'_p(n,0;D)}{b_p(n,0;D)} +  c'\cdot (1-\chi_d(p))^{-1}\cdot
\frac{1-\chi_d(p)+\chi_d(p)p^k- p^{k+1}}{1-p} .\tag8.36$$

We now use (ii) of Lemma~8.10.
If $\chi_d(p)=-1$, (8.36) gives
$$\align
K_p:=-1-k +&\frac{2(p+1)p^{k+1} +p^2-4p-1}{2(p^2-1)}
+  c'\cdot\frac12\cdot
\frac{p^{k+1}+p^k -2}{p-1}\tag8.37\\
\nass
{}&= -1-k + \frac12+ \frac{p^{k+1}+p^k-2}{2(p-1)}\\
\nass
{}&= -k + \frac{(p+1)(p^k-1)}{2(p-1)}.
\endalign
$$
If $\chi_d(p)=0$, (8.36) gives
$$\align
K_p:=-1-k + &\frac{2p(p^{k+1}-1)}{p^2-1}
+  c'\cdot
\frac{p^{k+1}-1}{p-1}\tag8.38\\
\nass
{}&= -1-k + \frac{p^{k+1}-1}{p-1}.
\endalign
$$
Thus, (8.27), (8.35), and these expressions for the coefficients $K_p$ of $\log(p)$ for $p\mid D$ yield (i) of
Theorem~8.8.

To prove (ii),  suppose that $m>0$ and that there is a unique prime $p\mid D$ for which $\chi_d(p)=1$. Then,
using (8.29) and (6.32),
we have
$$\align
\Cal E'_m(\tau,\frac12;D)&= \E'_m(\tau,\frac12;D) + \sum_{\ell|D} c'_\ell(\frac12)\cdot \E_m(\tau,\frac12;D/\ell)\tag8.39\\
\nass
\nass
{}&= - 2\,\delta(d;D/p)\,H_0(m;D)\cdot  q^m\cdot (p+1)^{-1}\cdot
(1+p - 2 p^{k+1})\log p\\
\nass
{}&\qquad\qquad\qquad + c'_p(\frac12)\cdot 2\,\delta(d;D/p)\,H_0(m;D/p)\,q^m\\
\nass
\nass
{}&= 2\,\delta(d;D/p)\,H_0(m;D)\cdot  q^m\cdot\bigg[-(p+1)^{-1}\cdot
(1+p - 2 p^{k+1})  + c'\cdot
p^k\, \bigg]\,\log(p)\\
\nass
\nass
{}&=2\,\delta(d;D/p)\,H_0(m;D)\cdot  q^m\cdot (p^k-1)\,\log(p),
\endalign
$$
as claimed.

Finally, we consider the constant term.
By (8.33), Corollary~8.4 and (6.32), we have
$$\align
&\Cal E'_0(\tau,\frac12;D)\tag8.40\\
\nass
{}&= \E'_0(\tau,\frac12;D) - \sum_{p\mid D} \frac{p-1}{p+1}\,\log(p)\cdot c(D/p)\,\L_{D/p}(2)\\
\nass
{}&= c(D)\L_D(2)\bigg[\,\frac12\,\log(v)
+ 1+ \log(D)-\log(\pi)-\gamma+ 2\frac{\zeta'(2)}{\zeta(2)} +\sum_{p\mid D}\bigg( \frac{2}{p^2-1} + \frac{1}{p+1}\bigg)\log(p)\bigg]\\
\nass
{}&=c(D)\L_D(2)\bigg[\,\frac12\,\log(v)
+ 1 -\log(\pi)-\gamma+ 2\frac{\zeta'(2)}{\zeta(2)} +\sum_{p\mid D} \frac{p\log(p)}{p-1}\bigg]\\
\nass
{}&=c(D)\,\L_D(2)\,\bigg[\frac12\,\log(v) - 2\frac{\zeta'(-1)}{\zeta(-1)} -1 + 2\log(2) + \log(\pi)
+\gamma + \sum_{p\mid D}\frac{p\log(p)}{p-1}\,\bigg]\tag13.1\\
{}&=c(D)\,\L_D(2)\,\bigg[\frac12\,\log(v) - 2\frac{\zeta'(-1)}{\zeta(-1)} -1 + 2C+ \sum_{p\mid D}\frac{p\log(p)}{p-1}\,\bigg],
\endalign
$$
where $C$ is as in Definition~3.4.
Here we use the fact that
$$c(D/p)\,\L_{D/p}(2) = -c(D)\,\L_D(2)\,\frac{p+1}{p^2-1}.\tag8.41$$
\qed\enddemo

For later comparison, we note that the coefficient $A'_m(\frac12,v)$ in the term
$$\Cal E'_m(\tau,\frac12;D) = A'_m(\frac12,v)\,q^m\tag8.42$$
in (i) of Theorem~8.8 can be written as a sum of  four quantities:
$$\align
&2\,\delta(d;D)\,H_0(m;D)\cdot\bigg[ \,
\frac12\, \log(d)  + \frac{L'(1,\chi_d)}{L(1,\chi_d)} -\frac12\log(\pi) -\frac12\gamma\,\bigg],\tag8.43\\
\nass
&2\,\delta(d;D)\,H_0(m;D)\cdot \frac12 J(4\pi mv) ,\tag8.44\\
\nass
&2\,\delta(d;D)\,H_0(m;D)\cdot\sum_{p \atop p\nmid D} \bigg(\ \log|n|_p - \frac{b_p'(n,0;D)}{b_p(n,0;D)}\ \bigg),\tag8.45\\
\nass
\noalign{and}
\nass
&2\,\delta(d;D)\,H_0(m;D)\cdot \,\sum_{p \atop p\mid D} K_p\,\log(p).\tag8.46
\endalign$$

\subheading{Part III. Computations: geometric}

\subheading{\Sec9. The geometry of $\Cal Z(m)$'s}

In this section we will prepare the calculation of the
coefficients of the generating series $\phi_{\deg}(\tau)$ and $\phi_{\text{\rm height}}(\tau)$ by
describing some of the geometry of the special cycles $\Cal Z(m)$.

It turns out that the primes of bad reduction (i.e.\ $p\mid D(B)$) play a
very special role. Namely,
$$\Cal Z(m)\times_{\Spec\, \Z}\Spec\, \Z[D(B)^{-1}]\ \text{is reduced and is finite
and flat over}\ \Spec\, \Z[D(B)^{-1}].\leqno(9.1)$$
We denote by $\Cal Z(m)^{\roman{horiz}}$ the closure of
$\Cal Z(m)\times_{\Spec\, \Z}\Spec\, \Z[D(B)^{-1}]$ in $\Cal Z(m)$ and call it the {\it
horizontal part of the special cycle.}

We first describe the generic fiber of $\Cal Z(m)$. As in section 5, let
$L=V\cap O_B$, let $L(m)$ be as in (5.7), and let
$$D_{\Cal Z(m)}=\prod\limits_{x\in L(m)}D_x\ \ ,\leqno(9.2)$$
as in (5.11).
Then $\Gamma$ acts on $D_{\Cal Z(m)}$ compatibly with its action on $D$, and we may
represent $\Cal Z(m)_{\C}$ as an orbifold mapping to $[\Gamma\setminus D]$
$$\Cal Z(m)_{\C}=[\Gamma\setminus D_{\Cal Z(m)}]\ \ .\leqno(9.3)$$
We next give the degree of this orbifold. We recall the following notation from (8.2). For $x\in
L$ with $Q(x)= m>0$, let
$$\kay=\Q[x]\simeq \Q[X]/(X^2+m)=\Q(\sqrt{-m})\ \ .\leqno(9.4)$$
Let
$O_{\smallkay}$ be its ring of integers and let $-d={\roman{disc}}(O_{\smallkay})$
be its discriminant. Then the discriminant of the order
$\Z[x]=\Z[X]/(X^2+m)$ is equal to $-4m$. Write $4m=n^2d$, as in (8.2).
We note that there is a map of discrete orbifolds
$$[\Gamma\back D_{\Cal Z(m)}] \lra [\Gamma\back L(m)],\tag9.5$$
which is $2$ to $1$.
\proclaim{Proposition 9.1} For $m>0$ and $\kay$ as above, \hfill\break
(i) if $\kay$ cannot be embedded into $B$, then
$\Cal Z(m)_{\Q}=\emptyset$,\hfill\break
(ii) otherwise,
$$\deg\, \Cal Z(m)_{\Q}=2\,\delta (m,D)\,H_0(m,D).$$
\endproclaim

Here the degree of $\Cal Z(m)_{\Q}$ is taken in the stack sense, i.e.\ each
geometric point $\eta$ of $\Cal Z(m)_{\Q}$ counts with multiplicity
$1/\vert \Aut(\eta)\vert$.

\demo{Proof} (i) For any $\C$-valued point $(A,\iota)$ of ${\Cal M}$ we
have an injection
$$\End(A,\iota)\hookrightarrow B\ \ .$$
Hence, if $(A,\iota, x)$ is a $\C$-valued point of $\Cal Z(m)$ we obtain an
injection $\kay=\Q[x]\hookrightarrow B$.

(ii) Using (9.8),  we obtain
$$\deg\, \Cal Z(m)_{\Q}=2\sum\limits_{\matrix \scr x\in L(m)\\\scr \mod \Gamma \endmatrix} {1\over \vert \Gamma_x\vert}\ \ .\tag9.6$$
The following result finishes the proof of the Proposition.
\proclaim{Lemma 9.2} If $m>0$, then
$$ \sum_{\matrix \scr x\in L(m)\\ \scr \mod \Gamma\endmatrix}
\frac{1}{|\Gamma_x|}  = \delta(d;D)\cdot H_0(m;D),$$
where
$
H_0(m;D)$ and $\delta(d;D)$ are given by (8.19) and (8.20).
\endproclaim
\demo{Proof}
For any $x\in V(\Q)\cap O_B$ with $Q(x)=m$, there is an associated embedding
$i_x:\Q(\sqrt{-m}) \rightarrow B$, taking $\sqrt{-m}$ to $x$. The order
$O_{c^2d} = i_x^{-1}(O_B)$ is an invariant of the $\Gamma=O_B^\times$ conjugacy class
of $i_x$ and $i_x$ is an optimal embedding of $O_{c^2d}$,
in the terminology of Eichler, \cite{\Eichler}.
Recall that the order $\Z[\sqrt{-m}]$ has discriminant $-4m$ and hence, conductor $n$.
\define\Optoc{\text{\rm Opt}(O_{c^2d},O_B)}
Let
$$\Optoc = \{\, i:\Q(\sqrt{-m})\rightarrow B\mid i^{-1}(O_B) = O_{c^2d}\,\}/\Gamma\tag9.7$$
be the set of $\Gamma$ orbits of optimal embeddings.
Recall  that the order $i^{-1}(O_B)$ is maximal at all primes $p|D(B)$.
The following fact is classical, \cite{\Eichler}:
$$|\Optoc| = \delta(d;D)\cdot h(c^2d).\tag9.8$$
Using (8.21), we have
$$\align
\bigg(\sum_{\matrix \scr x\in L(m)\\\scr Q(x)=m\\ \scr \mod \Gamma\endmatrix}
|\Gamma_x|^{-1} \bigg)& = \sum_{\scr c|n \atop \scr (c,D) =1} |\Optoc|\cdot|O_{c^2d}^\times|^{-1}\\
\nass
{}&= \,\delta(d;D)\,\frac{h(d)}{w(d)}\sum_{\scr c|n \atop \scr (c,D) =1}
c\prod_{\ell|c}(1-\chi_d(\ell)\ell^{-1})\tag9.9\\
\nass
{}&= \delta(d;D)\,H_0(m;D).
\endalign
$$
\qed\enddemo
\enddemo

\demo{Proof of Proposition~7.1} Comparing the expression just found for $\deg \Cal Z(m)_\Q$ together with (2.7), we have
$$\phi_{\deg}(\tau) = \zeta_D(-1) + \sum_{m>0}2\,\delta(d;D(B))\,H_0(m;D(B))\,q^m,$$
which coincides with $\Cal E(\tau,\frac12;D(B))$, via Proposition~8.5, so that Proposition~7.1 is proved.
\qed\enddemo

\medskip\noindent
{\bf Remark 9.3.} The map $\Cal Z(m)_{\Q}\to {\Cal M}_{\Q}$ is not a closed
immersion, hence $\Cal Z(m)_{\Q}$ is not a divisor on ${\Cal M}_{\Q}$. In fact,
the morphism is of degree 2 over its image. To see this, note that if
$(A,\iota, x)\in \Cal Z(m)(\C)$, then $\End (A,\iota)_{\Q}=\kay=\Q(\sqrt{-m})$.
Hence the only other point of $\Cal Z(m)(\C)$ mapping to $(A,\iota)\in {\Cal
M}(\C)$ is $(A,\iota, -x)$. That the degree is 2 even in the stack sense
follows from

\centerline{\hfill $\Aut(A,\iota)=\Aut(A,\iota, x).$\hfill\hfill\qed }

The stack $\Cal Z(m)$ can have some pathological features in characteristic $p$
for $p\mid  D(B)$. Namely, as already mentioned in section 5, it can
happen that $\Cal Z(m)$ has dimension 0 (only if $\kay=\Q(\sqrt{-m})$ does not
embed into $B$, cf.\ (i) of Proposition 9.1) and also that $\Cal Z(m)$ has embedded
components. This leads us to introduce the Cohen-Macauleyfication
$\Cal Z(m)^{\roman{pure}}$, \cite{\krinvent}. In the case that $\Cal Z(m)$ has
dimension 0, this is empty. In all other cases $\Cal Z(m)$ may be considered a
divisor on ${\Cal M}$ (but note that since $\Cal Z(m)$ is not a closed substack
of ${\Cal M}$, the degree  of $\Cal Z(m)$ over its image must be taken into
account). But even after $\Cal Z(m)$ is replaced by $\Cal Z(m)^{\roman{pure}}$, one
interesting feature remains, namely the existence of vertical components
in characteristic $p\mid  D(B)$. We write
$$\Cal Z(m)^{\roman{pure}}=\Cal Z(m)^{\roman{horiz}}+\Cal Z(m)^{\roman{vert}}\leqno(9.10)$$
(equality of ``divisors'' on ${\Cal M}$), where $\Cal Z(m)^{\roman{vert}}$ is
the sum with multiplicities of the irreducible vertical components in
characteristic $p$ as $p$ runs over primes dividing $D(B)$. We note that if we redefine
$$\hat \Cal Z(m,v)=(\,\Cal Z(m)^{\roman{pure}}, \Xi(m,v)\,)\leqno(9.11)$$
then the expression $\langle\hat \Cal Z(m,v),\hat\omega\rangle$ appearing in
the definition of $\phi_{\text{\rm height}}(\tau)$ remains unchanged,
cf. \cite{\krinvent}, section 4. We remark in passing
that if $\Cal Z(m)_{\Q}=\emptyset$, then $\Xi(m,v)=0$, cf.\ (5.8)

Taking (9.10)
into account and using formula (5.11) of Bost \cite{\bost}, we may write
$$\langle\,\hat \Cal Z(m,v),\hat\omega\,\rangle =
h_{\hat\omega}(\Cal Z(m)^{\roman{horiz}})+h_{\hat\omega}(\Cal Z(m)^{\roman{vert}})
+{1\over 2}\int_{[\Gamma\setminus D]}\Xi(m,v)c_1(\hat\omega)\ .\leqno(9.12)$$
Also note that the second summand on the right hand side may be written as the sum
over contributions of the bad fibers ,
$$h_{\hat\omega}(\Cal Z(m)^{\roman{vert}})=\sum_{p\mid
D(B)}h_{\hat\omega}(\Cal Z(m)_p^{\roman{vert}})\ \ ,\leqno(9.13)$$
where, since $\Cal Z(m)^{\roman{vert}}$ has empty generic fiber,
$$h_{\hat\omega}(\Cal Z(m)_p^{\roman{vert}})=\deg (\omega\vert
\Cal Z(m)_p^{\roman{vert}})\,\log(p) \ .\leqno(9.14)$$
Here $\Cal Z(m)_p^{\roman{vert}}=\Cal Z(m)^{\roman{vert}}\times_{\Spec\, \Z}\Spec\,
\Z_p$.

In the next three sections we will evaluate explicitly each summand on the
right hand side of (9.12).

\subheading{10. Contributions of horizontal components}

In this section, we compute the quantity $h_{\hat\o}(\Cal Z(m)^\hor)$. This will be done in
two steps. We first express this quantity in terms of the Faltings heights of
certain abelian surfaces which are isogenous to products of CM--elliptic curves.
We then determine the effect of the isogeny on the Faltings height.

Observe that $\Cal Z(m)^\hor$ is a union of horizontal integral substacks
$$\Cal Z(m)^\hor = \sum_{\xi\in \Cal Z(m)^\hor_\Q} \Cal Z_\xi,\tag10.1$$
where $\xi$ is the generic point of $\Cal Z_\xi$. Let $\tilde \Cal Z_\xi$ be the normalization of $\Cal Z_\xi$ and
let $j_\xi:\tilde\Cal Z_\xi \rightarrow \M$ be the composition of
the  normalization map $\tilde \Cal Z_\xi\rightarrow\Cal  Z_\xi$ with the
morphism $\Cal Z_\xi \rightarrow \M$.
By linearity and the definition of $h_{\hat\o}$ for horizontal cycles, cf. (4.4) section 4 and \cite{\bost}, we have
$$\align
h_{\hat\o}(\Cal Z(m)^\hor) & = \sum_{\xi\in \Cal Z(m)^\hor_\Q} h_{\hat\o}(\Cal Z_\xi)\\
\nass
{}&= \sum_{\xi\in\Cal  Z(m)^\hor_\Q}\degh j_\xi^* \hat\o\cdot \frac1{|\Aut(\bar\xi)|}\tag10.2\\
\nass
{}&= \frac{1}{|L:\Q|}\cdot\sum_{\eta\in\Cal  Z(m)^\hor(L)} \degh j_\eta^* \hat\o\cdot \frac1{|\Aut(\eta)|},
\endalign
$$
for any sufficiently large number field $L\subset \bar\Q$, where $\eta$ runs over the
$L$ points of $\Cal Z(m)$ and
where the factors involving $\Aut({\bar \xi}) = \Aut((A,\iota,x)_{\bar \xi})$
and $\Aut(\eta) = \Aut((A,\iota,x)_\eta)$ come in
due to the stack.  Here we also write $\eta$ for the extension to
$\eta:\Spec(O_L)\rightarrow \M$. We may assume that $A_\eta$,
the abelian variety over $L$ determined by $\eta$,  has
semistable reduction over $L$. Then,
by definition, the Faltings height $\hfal^*(A_\eta)$ is given by
$$\hfal^*(A_\eta) = |L:\Q|^{-1} \degh j_\eta^* \hat\o.\tag10.3$$
Here the notation $\hfal^*$ indicates that we have used the metric $||\ ||$ given in
Definition~3.4, rather than the standard metric of (3.11). This point is discussed further below.
Then:
$$\align
h_{\hat\o}(\Cal Z(m)^\hor) &= \sum_{\eta\in \Cal Z(m)^\hor(L)} \hfal^*(A_\eta)\cdot \frac{1}{|\Aut(\eta)|}\tag10.4\\
\nass
{}&= 2\sum_{\matrix\scr x\in L(m)\\\scr \mod \Gamma\endmatrix} \hfal^*(A_x)\cdot\frac1{|\Gamma_x|}.
\endalign
$$
In the last expression,
we have used the description of $\Cal Z(m)(\C) = \Cal Z(m)^\hor(\C)$
as the orbifold $[\Gamma\back D_{\Cal Z(m)}]$, as in section 9, where
the map $[\Gamma\back D_{\Cal Z(m)}] \rightarrow [\Gamma\back L(m)]$ is $2$ to $1$,
together with the fact that
the abelian varieties associated to the two points in $D_x$ (i.e., having
opposite complex structures) have the same Faltings height.

We next turn to the computation of the Faltings height $\hfal^*(A)$ of an abelian surface $A$
occuring in a triple $(A,\iota,x)$ where $\iota:O_B\rightarrow \End(A)$ and $x\in \End(A,\iota)$
is a special endomorphism with $Q(x)=m$, all
defined over a number field $L$, where they all have good reduction.
Let
$\phi_x:\kay\rightarrow \End^0(A)$ be the embedding determined by $\phi_x(\sqrt{-m})= x$ and let
$$O_{c^2d} = \phi_x^{-1}(\Q[x]\cap \End(A)),\tag10.5$$
where $O_{c^2d}$ is the order in $O_{\smallkay}$ of conductor $c$.
In this case, we will say that the triple $(A,\iota,x)$ is of type $c$.
Recall that the order $\Z[\sqrt{-m}\,]\subset  O_{c^2d}\subset O_{\smallkay}$ has discriminant $4m$,
and that we have written $4m = n^2 \,d$,  where $-d$ is the discriminant of $O_\smallkay$. Then the order $\Z[\sqrt{-m}\,]$
has conductor $n$ and $c\mid n$.  Note, in addition, that the order $O_{c^2d}$ must be maximal
at all primes $p\mid D(B)$, and hence $(c,D(B))=1$.

Here the key point is that the special endomorphism
$x$ with $x^2=-Q(x)\cdot 1_A$ forces $A$ to be isogenous to a product of elliptic
curves with CM by $\kay = \Q(\sqrt{-m})$. The isogeny of interest
is constructed as follows.
Since $\kay$ splits $B$, we can choose an embedding
of $\psi:\kay\hookrightarrow B$ such that
$$O_\smallkay = \psi^{-1}(\OB),\tag10.6$$
i.e., an optimal embedding, in the sense of Eichler.
Then the endomorphisms
$$\a_x^{\pm} := x\pm \iota\circ\psi(\sqrt{-m})\in \End(A)\tag10.7$$
satisfy
$$\a_x^++\a_x^- = 2x\qquad\text{and}\qquad \a_x^+\cdot \a_x^- = 0.\tag10.8$$
Let
$$E^\pm=(\ker(\a_x^\pm))^0\tag10.9$$
be the identity component of the kernel of the endomorphism $\a_x^\pm$.
Choose an element $\eta\in O_B$ with $\tr(\eta)=0$ and such that
conjugation by $\eta$ induces the Galois automorphism on $\psi(\kay)$.
Note that
$$\iota(\eta)\,\a_x^+ = \a_x^-\,\iota(\eta),\tag10.10$$
and so
$$\iota(\eta)E^\pm = E^\mp.\tag10.11$$
Thus
$E^\pm$ are both elliptic curves, and we obtain an isogeny
$$u_L:E^+\times E^- = (\ker(\a_x^+))^0\times (\ker(\a_z^-))^0 \lra A,\tag10.12$$
rational over $L$.
Moreover, the elliptic curves $E^\pm=(\ker(\a_x^\pm))^0$ are stable under $\psi(O_\smallkay)$,
i.e., have complex multiplication by the full ring of integers of $\kay$.
Note that the kernel of $u_L$ is the subgroup
$$M_L=(\ker(\a_x^+))^0\cap (\ker(\a_x^-))^0\tag10.13$$
embedded antidiagonally in $E^+\times E^-$.

The behavior of the Faltings height under isogeny is nicely described in the
article of Raynaud, \cite{\raynaud}. Our normalization of the Faltings height is the following.
For an abelian variety $B$ of dimension $g$ over a number field $L$, let $N(B)$ be the N\'eron model of $B$ over $S=\Spec(O_L)$, and
let $\o_B = \e^*(\wedge^g \O_{N(B)/S})$ be the pullback by the zero section $\e$ of
the top power of the sheaf of relative differentials on $N(B)$. This
invertible sheaf on $\Spec(O_L)$ has natural metrics; if $\s:L\hookrightarrow \C$ is an embedding of
$L$, then a section $\b$ of $\o_B$ determines a holomorphic $g$--form on $B_\s(\C)$,
and\footnote{ Note that we use the factor $\left(\frac{i}{2\pi}\right)^g$ rather than $\left(\frac{i}{2}\right)^g$.
This is the normalization used in Bost, \cite{\bost}, for example.}
$$||\b||^2_{\s,\text{\rm nat}} = \bigg| \left(\frac{i}{2\pi}\right)^g \int_{B_\s(\C)} \b\wedge \bar\b \ \bigg|.\tag10.14$$
If $B$ has semi--stable reduction, then the Faltings height of $B$ is given by
$$\hfal(B) = |L:\Q|^{-1}\,\degh(\o_B)\tag10.15$$
where $\degH(\o_B)$ denotes that Arakelov degree of $\o_B$. The
quantity $\hfal(B)$ does not depend on the choice of $L$ over which $B$ has semi--stable
reduction.

In view of the normalization used in Definition~3.4, we introduce the metrics
$$||\b||^2_{\s} = \bigg| \left(e^{-C}\frac{i}{2\pi}\right)^g \int_{B_\s(\C)} \b\wedge \bar\b \ \bigg|,\tag10.16$$
where, as before,
$C= \frac12\,\big(\,\log(4\pi)+\gamma\,\big)$.
We denote the resulting height by $\hfal^*(B)$. The two heights are related by
$$\hfal^*(B) = \hfal(B) +\frac12 g C.\tag10.17$$

Assume that $A$ has good reduction over $L$ and let $u_L: A\rightarrow B$ be an isogeny defined over $L$.
Let $u:N(A)\rightarrow N(B)$
be the resulting homomorphism of N\'eron models with $M:=\ker(u)$.
Then, as a special case of \cite{\raynaud}, p.205,
$$\hfal(B) = \hfal(A) + \frac12\log(\deg(u_L)) - |L:\Q|^{-1}\, \log|\e^*(\O^1_{M/R})|.\tag10.18$$
The quantity $\delta(u):= \log|\e^*(\O^1_{M/S})|$ is a sum of local contributions as follows.
For each prime $v$ of $L$ with \hbox{$v\mid\deg(u_L)$}, let $R_v$ be the completion of $O_L$ at $v$
and let
$M_v = M\tt_{O_L}R_v$. Then
$$\delta(u) = \log|\e^*(\O^1_{M/R})| = \sum_{v\mid\deg(u_L)} \log|\e^*(\O^1_{M_v/R_v})|.\tag10.19$$
For convenience, we set
$$\delta_v(u)=\log|\e^*(\O^1_{M_v/R_v})|.\tag10.20$$
This quantity is invariant under base change in the sense that
if $L'$ is a finite extension of $L$ and if $u'$ is the base change of $u$ to
$\Spec(O_{L'})$, then
$$\delta_v(u) = \sum_{w\mid v} \delta_w(u'),\tag10.21$$
where $w$ runs over the primes of $L'$ dividing $v$.

We now return to our isogeny $u_L$, noting that $A$ and $E^\pm$
all have good reduction over $L$.  Let
$$u:N(E^+)\times N(E^-)\lra N(A)\tag10.22$$
be the homomorphism induced by $u_L$ and let $M$ be its kernel.
To calculate
$\delta_v(u)$ for a prime $v\mid \deg(u_L)$, we pass to the $p$-divisible groups, where $p$ is the
residue characteristic for $v$.

Let
$G= G^+\times G^-$ (resp. $A(p)$)
be the p-divisible group over $R_v$ associated to $E^+\times E^-$ (resp. $A$) so that
we have an exact sequence
$$0 \lra C \lra G \lra A(p) \lra 0\tag10.23$$
determined by the isogeny $u$. Since the prime to $p$ part of the kernel of
$u$ is automatically \'etale over $R_v$, the invariant $\delta_v(u)$ depends only
on the $p$--divisible groups and hence on $C$.
The isogeny (10.23) corresponds to a submodule $T'$ of $V(G)$, the rational Tate module of $G$,
$$T(G)\subset T'\subset V(G) = V(G^+)\oplus V(G^-).\tag10.24$$
The fact that $E^\pm$ maps injectively into $A$ implies that $G^\pm \hookrightarrow A(p)$
and hence
$$T'\cap V(G^\pm) = T(G^\pm).\tag10.25$$
Thus there are isomorphisms
$$\pr^+(T')/T(G^+) \lisoarrow T'/T(G) \isoarrow \pr^-(T')/T(G^-).\tag10.26$$

\proclaim{Proposition 10.1} Suppose that $p$ splits in $\kay$. Then $\ord_p(\deg(u_L)) = 2\,\ord_p(c)$.
Moreover, for any place $v$ of $L$ with $v\mid p$, the group $C$ is \'etale over $R_v$, and hence
$$\delta_v(u)=\log|\e^*(\O^1_{M_v/R_v})|=\log|\e^*(\O^1_{C/R_v})|=0.$$
\endproclaim
\demo{Proof}
We write
$$\kay_p\simeq \Q_p\oplus \Q_p, \qquad \a\mapsto (\a_1,\a_2),\tag10.27$$
and let $\lambda_1$ and $\lambda_2$ be the corresponding algebra homomorphisms from $\kay_p$ to $\Q_p$.
Let $G_0$ (resp. $G_{\text{\'et}}$) be the connected (resp.  \'etale) part of $G$, and note that, for example
$$G_0 = G_0^+\times G_0^-.\tag10.28$$
Since the action of $\Cal O:= O_{\smallkay}\tt\Z_p\simeq \Z_p\oplus \Z_p$ preserves $G^\pm_0$, the action of $\Cal O$
on the Tate module $T(G^\pm)$ must have the form
$$\align
T(G^+) &= T(G^+_0)\oplus T(G^+_{\text{\'et}})\simeq \Z_p\oplus \Z_p, \qquad \l_1\oplus\l_2\tag10.29\\
\nass
T(G^-) &= T(G^-_0)\oplus T(G^-_{\text{\'et}})\simeq \Z_p\oplus \Z_p, \qquad \l_2\oplus\l_1,
\endalign
$$
where the switch of characters is due to the fact that the isogeny induced by
$\eta$ is $\Cal O$--antilinear
but must preserve the connected--\'etale decomposition.
Thus there is a canonical decomposition
$$T(G) = T(G_{0}^+)\oplus T(G_0^-)\oplus T(G_{\text{\'et}}^+) \oplus T(G_{\text{\'et}}^-),
\qquad \l_1\oplus\l_2\oplus \l_2\oplus\l_1.\tag10.30$$
Since the $\Z_p$--lattice $T'\subset V(G)$ is stable under $\Cal O$, it must be generated by
coset representatives of the form $(x,0,0,w)$ and $(0,y,z,0)$. Condition (10.24) implies, for example,
that if $x=0$ for such a representative, then $w=0$ (i.e., lies in $\Z_p$). It follows that
$$T' = \Z_p\cdot (p^{-r},0,0,\e_1 p^{-r}) + \Z_p\cdot (0,p^{-s},\e_2 p^{-s},0) + T(G).\tag10.31$$
for units $\e_1$ and $\e_2$ and non--negative integers $r$ and $s$. In addition, we can choose the element $\eta\in O_B$ above so that
$\eta^2$ is prime to $p$ and so $\eta$ induces an automorphism on $A(p)$. Since $T'$ is stable under
this automorphism, we must have $r=s$ and $\e_1\equiv \e_2 \mod p^r$.
On the other hand, the $\phi_x$--action of $\Cal O$ on $V(G)$ is given by
$$\l_2 \oplus \l_2\oplus \l_1\oplus \l_1.\tag10.32$$
The $\phi_x$ action of $\a\in \Cal O$ preserves $T'$ if and only if
$$\a\in O_{c^2d}\tt \Z_p =\{ (a_1,a_2)\Cal O\mid a_1\equiv a_2 \mod p^{\ord_p(c)}\}.\tag10.33$$
Thus, we conclude that $r = \ord_p(c)$. This proves that $\ord_p(\deg(u_L)) = 2\,\ord_p(c)$.

To finish the proof of the Proposition, it suffices to show that
$$T'\cap V(G_0) = T(G_0),\tag10.34$$
so that the projection to the \'etale part $G\rightarrow G_{\text{\'et}}$
induces an isomorphism on $C$. But this is clear from our description of the
coset representatives.
\qed\enddemo

\proclaim{Proposition 10.2}
Suppose that $p$ is inert or ramified in $\kay$. \hfill\break
(i) If $p\nmid D(B)$, then
for any place $v$ of $L$ with $v\mid p$, there is a factorization
$u = u^\dagger\circ u^o$ with isogenies $u^\dagger$ and $u^o$
such that
$$\delta_v(u^o) = \frac12\cdot|L_v:\Q_p|\cdot \big(\,\ord_p(\deg(u_L))-2\,\ord_p(c)\,\big),$$
and $\ord_p(\deg(u^\dagger)) = 2\ord_p(c) = 2s$. Moreover,
$$\delta_v(u^\dagger) = |L_v:\Q_p|\, \frac{(1-p^{-s})\cdot (1-\chi(p))}{(1-p^{-1})\cdot(p-\chi(p))}\,\log(p).$$
Here $\chi(p)=-1$ if $p$ is inert and $\chi(p)=0$ if $p$ is ramified in $\kay$. \hfill\break
(ii) If $p\mid D(B)$, then $\ord_p(\deg(u_L))=0$ and $\delta_v(u)=0$.
\endproclaim

\demo{Proof}
Let $\Bbb F_q$ be the residue
field $\Cal O/\pi\Cal O$, where $\pi$ is a fixed prime element of $\Cal O$.
Also write $\Cal O_s =  O_{c^2d}\tt\Z_p$, where $s=\ord_p(c)$.  For convenience, we temporarily
write $L$ in place of $L_v$ and $\Cal O_L$ in place of $O_{L_v}$.

Now $G^\pm=E^\pm(p)$ is a formal group of dimension $1$ and height $2$ over $R_v=\Cal O_{L_v}$ with an action
of $\Cal O$, i.e., a special formal $\Cal O$--module in the sense of Drinfeld.
We consider the sequence
$$C\lra G^+\times G^- \overset{u}\to{\lra}\ A(p)\tag10.35$$
and note that $O_B\tt \Cal O_s$ acts on $A(p)$.  After replacing $L$ by a finite
extension and using the invariance property (10.21), we may assume that $G_0:=G^+\simeq G^-$.

First suppose that $p\nmid D(B)$. Then, fixing an isomorphism
$$O_B\tt_\Z\Cal O_s \simeq M_2(\Z_p)\tt_{\Z_p}\!\Cal O_s \simeq M_2(\Cal O_s),\tag10.36$$
we may write $A(p)\simeq G_s\times G_s$, where $G_s$ is a $1$--dimensional
formal group of height $2$. Since $(A,\iota,x)$ was supposed to be of type $c$,
we have $\End(G_s) = \Cal O_s$, hence the notation for $G_s$, consistent with that
for $G_0$. The isogeny $u$ corresponds to an inclusion of Tate modules
$$T(G_0)\oplus T(G_0) = T(G_0\times G_0) \subset T(A(p)) = T(G_s\times G_s) = T(G_s)\oplus T(G_s).\tag10.37$$
Note that the two direct sums here are not necessarily compatible.
Let $T(G_s)^\dagger$ be the largest $\Cal O$--module contained in the $\Cal O_s$--module $T(G_s)$.
Note that
$$T(G_s)^\dagger\oplus T(G_s)^\dagger\tag10.38$$
is the largest $\Cal O$--module contained in $T(G_s)\oplus T(G_s)$. Hence the
inclusion (10.37) gives rise to a chain of inclusions
$$T(G_0)\oplus T(G_0)\subset T(G_s)^\dagger\oplus T(G_s)^\dagger\subset T(G_s)\oplus T(G_s).\tag10.39$$
Hence the isogeny $u$ factors as
$$G_0\times G_0 \overset{u^o}\to{\lra}\ G_s^\dagger\times G_s^\dagger \overset{u^\dagger}\to{\lra}\ G_s\times G_s.\tag10.40$$
By the elementary divisor theorem, we can then find automorphisms
$\a$ of $G_0\times G_0$ and $\b$ of $G_s^\dagger\times G_s^\dagger$ such that
$\b\circ u^o\circ \a$ is of the form
$u^o_1\times u^o_2$ for isogenies $u^o_i:G_0\rightarrow G_s^\dagger$, $i=1$, $2$.
Each of the isogenies $u^o_i$ corresponds to an inclusion of Tate modules
$T(G_0)\subset T(G_s^\dagger)\subset V(G_0)$. Since both Tate modules are free
$\Cal O$--modules of rank $1$, there exists an isomorphism $G_0\simeq G_s^\dagger$
such that $u_i^o$ is given by multiplication by an element $\nu_i\in \Cal O$.
On the other hand, the isogeny $u^\dagger$ is of the form $u^\dagger = u^\dagger_1\times u^\dagger_1$,
where the isogeny $u^\dagger_1$ is determined by the inclusion $T(G_s)^\dagger\subset T(G_s)$.
If we choose an isomorphism $T(G_s)^\dagger \simeq \Cal O$, then
$$T(G_s) \simeq p^{-s}\e\,\Z_p + \Cal O,\tag10.41$$
where $\e$ is a unit in $\Cal O$. This implies that the degree of the isogeny $u_1^\dagger$ is $p^s$,
as claimed.

The contribution $\delta_v(u^\dagger) = 2\delta_v(u^\dagger_1)$ of the isogeny $u^\dagger$
to the invariant $\delta_v(u) = \delta_v(u^0) + \delta_v(u^\dagger)$
can now be obtained from the following result, whose proof we include, for the sake of completeness.
\proclaim{Proposition 10.3} {\rm (Nakkajima--Taguchi, \cite{\naktag})} Let $\kay_p/\Q_p$ be a
quadratic extension and let $L/\kay_p$ be a finite extension. Let $G_0$ be a
one-dimensional formal ${\Cal O}$-module over $\Cal O_L$. For $s\geq 0$
suppose that $\lambda:G_0\to G_s$ is an isogeny of degree $p^s$ over $\Cal O_L$ such that
$\End(G_s)={\Cal O}_s$. Let $D={\roman{Ker}}\, \lambda$. Then
$${\roman{log}}\vert \varepsilon^*\Omega^1_{D/\Cal O_L}\vert = \frac12\,{\vert
L:\Q_p\vert}\, {(1-p^{-s})\cdot (1-\chi(p))\over (1-p^{-1})\cdot
(p-\chi(p))}\,{\roman{log}}\, p\ \ .$$
Here $\chi(p)$ is $-1$ {\rm (}resp. $0${\rm )} depending on whether $\kay_p/\Q_p$ is
unramified or ramified.
\endproclaim

\demo{Proof} (Sketch) Let $G_0$ be defined by the formal group law
$g(X,Y)\in \Cal O_L[[X,Y]]$. Then by Serre's isogeny formula, \cite{\grossqc},
$$D=\Spec\, \Cal O_L[[X]]/\prod_{d\in D}g(X,d).\tag10.42$$
It follows that
$$\align
\varepsilon^*\Omega^1_{D/\Cal O_L}
&
=\Cal O_L/\bigg(\prod_{d\in D}g(X,d)\bigg)'_{X=0}
\tag10.43
\\
&
=\Cal O_L/(\prod_{d\in D\setminus \{ 0\}}d)\ \ .
\endalign$$
A consideration of the Newton polygon of $[\pi^r]_{G_0}=\pi^rX+\ldots$ shows
that if $d\in D\setminus\{ 0\}$ is of precise $\pi$-order $r$ then
$$\ord(d)= {1\over q^r-q^{r-1}}.\tag10.44$$
Here, as before, $\pi$ denotes a prime element of $\kay_p$ and the $\ord$
function is normalized by $\ord(\pi)=1$. It follows that
$$\align
{\roman{lg}}_{\Cal O_L}(\varepsilon^*\Omega^1_{D/\Cal O_L})
&
=e_{L/\smallkay_p}\cdot\sum_{d\in D\setminus\{ 0\}}\ord(d)\tag10.45
\\
&
=e_{L/\smallkay_p}\cdot \sum_{r=1}^\infty\cdot {\ell(r)\over q^r-q^{r-1}}\ \ .
\endalign$$
Here $\ell(r)$ denotes the number of elements of $D(\bar L)$ of exact
$\pi$-order $\ell(r)$. But the conditions $\deg(\lambda)=p^s$ and
$\End(G_s)=\Cal O_s$ imply that
$$D(\bar L)\ \cong\ T(G_s)/T(G_0)\ \cong\ \big(\,p^{-s}\varepsilon\,\Z_p+{\Cal O}\,\big)/{\Cal O}\ \cong\ p^{-s}\,\Z_p/\Z_p,\tag10.46$$
where $\varepsilon$ is a unit in $\Cal O$.
Hence, for $r\ge1$.
$$\ell(r) = \cases p^r-p^{r-1}&\text{ $1\leq r\leq s$, \quad\  if $\kay_p/\Q_p$ is unramified}\\
\nass
p^{r\over 2}-p^{{r\over 2}-1}&\text{ $1\leq r\leq 2s$,  $2\mid r$, if $\kay_p/\Q_p$ is ramified}
\endcases\tag10.47
$$
and is zero in all other cases. It follows that
$$\align
{\roman{log}}\vert\varepsilon^*\Omega^1_{D/\Cal O_L}\vert
&
= f_{L/\Q_p}\cdot
{\roman{lg}}_{\Cal O_L}(\varepsilon^*\Omega^1_{D/\Cal O_L})\cdot{\roman{log}}\, p\tag10.48
\\
&
= f_{L/\Q_p}\cdot e_{L/\smallkay_p}\cdot\sum_{r=1}^s {\ell(r)\over
q^r-q^{r-1}}{\roman{log}}\, p\ \ ,
\endalign$$
which yields the assertion.\qed
\enddemo

Finally, the contribution $\delta_v(u^o) = \delta_v(u^o_1)+\delta_v(u^o_2)$
is given by the following result.

\proclaim{Lemma 10.4}
For extensions $L/\kay_p/\Q_p$ as in the previous Proposition, suppose that $G_0$ is a
one-dimensional formal $\Cal O$-module over $\Cal O_L$. Let $\lambda:G_0\to
G_0$ be the isogeny given by multiplication by a non-zero element
$\nu\in{\Cal O}$. Let $D={\roman{ker}}(\lambda)$. Then
$${\roman{log}}\vert\varepsilon^*\Omega^1_{D/\Cal O_L}\vert = \frac12\,{\vert
L:\Q_p\vert}\, {\roman{log}}(\deg\, \lambda)\ \ .$$
\endproclaim

\demo{Proof} Let $s=\ord(\nu)$. Then, there are $q^r-q^{r-1}$ elements in
$D(\bar L)\simeq \pi^{-s}{\Cal O}/{\Cal O}$ of exact order $\pi^r$, for
$1\leq r\leq s$, and none for all other $r$. Hence
$${\roman{lg}}_{\Cal O_L}(\varepsilon^*\Omega^1_{D/\Cal O_L})=
e_{L/k_p}\cdot\sum_{r=1}^s {q^r-q^{r-1}\over q^r-q^{r-1}} =e_{L/k_p}\cdot
s.\tag10.49$$
It follows that
$${\roman{log}}\vert\varepsilon^*\Omega^1_{D/O_L}\vert = e_{L/k_p}\cdot
f_{L/\Q_p}\cdot s\cdot {\roman{log}}\, p= \frac12\,|L:\Q_p|\cdot s\cdot \log(q).\tag10.50$$
The assertion follows since the isogeny $\lambda$ has degree $q^s$.\qed
\enddemo

\enddemo

Finally, we consider the case $p\mid D(B)$, and we recall that $\ord_p(c)=0$,
so that $\Cal O_s=\Cal O_0 =\Cal O$. Once again, we consider the action
of $O_B\tt\Cal O$ on $A(p)$. We fix an isomorphism
$$B\tt\kay_p \simeq M_2(\kay_p), \qquad\text{\rm with}\qquad O_B\tt\Cal O \hookrightarrow M_2(\Cal O).\tag10.51$$
and such that, for $\a\in \kay_p$,
$$1\tt\a \mapsto \pmatrix \a&{}\\{}&\a\endpmatrix\qquad\text{\rm and}\qquad \psi(\a)\tt1\mapsto \pmatrix \a&{}\\{}&\a^\s\endpmatrix,\tag10.52$$
where $\psi$ is the embedding of $\kay$ into $B$ chosen above.
\proclaim{Lemma 10.5} (i) If $\kay_p/\Q_p$ is unramified, then the image of $O_B\tt\Cal O$
in $M_2(\Cal O)$ is the order
$$\pmatrix \Cal O&\Cal O\\p\Cal O&\Cal O\endpmatrix.$$
(ii) If $\kay_p/\Q_p$ is ramified, then there is an element $\lambda\in (O_B\tt\Cal O)^\times$
whose image in $GL_2(\F_p)$ under the composition of maps
$$O_B\tt\Cal O \lra M_2(\Cal O) \lra M_2(\F_p)$$
has eigenvalues which are not rational over $\F_p$. \hfill\break
(iii) In the ramified case, suppose that $\Lambda$ is an $\Cal O$--lattice
contained in $\kay^2$ which is stable under $O_B\tt\Cal O$. Then $\Lambda$
is homothetic to $\Lambda_0= \Cal O^2$.
\endproclaim
{\bf Remark 10.6.} In fact, if $p\ne2$ and $\kay_p/\Q_p$ is ramified, then
the image of $O_B\tt\Cal O$
in $M_2(\Cal O)$ is conjugate to $\text{\rm red}^{-1}(\F_{p^2})$, where
$$\text{\rm red}: M_2(\Cal O) \lra M_2(\F_p)$$
is the reduction modulo $\pi$, and $\F_{p^2}$
is the nontrivial quadratic extension of $\F_p$, viewed as a subalgebra of
$M_2(\F_p)$.

\demo{Proof} In the unramified case, we may write
$$O_B\tt\Z_p = \Cal O\lan\Pi\ran,\tag10.53$$
where $\Pi\in O_B$ is an element with $\Pi^2 = p$ which normalizes $\psi(\kay_p)$
and acts on it by the Galois automorphism $\s$, i.e., $\Pi \a = \a^\s \Pi$.
The image of $O_B\tt\Cal O$ in $M_2(\Cal O)$ is then generated by the elements
of the form
$$\pmatrix \a&{}\\{}&\a\endpmatrix, \qquad\pmatrix \a&{}\\{}&\a^\s\endpmatrix,
\qquad\text{and}\qquad \pmatrix {}&1\\p&{}\endpmatrix,\tag10.54$$
for $\a\in\Cal O$. Since $\kay_p$ is unramified, there is an $\a\in \Cal O$
such that $\a-\a^\s$ is a unit, and hence such elements generate the Eichler order
as claimed in (i). \hfill\break
Next suppose that $\kay_p$ is ramified. Let $\kay_o$ be the unramified quadratic
extention of $\Q_p$ and let $\Cal O_o$ be its ring of integers, with generator
$\lambda$ having unit norm. Again, we can write
$O_B\tt\Z_p = \Cal O_o\lan\Pi\ran$.
But now, the image of $\lambda\tt1\in (O_B\tt\Cal O)^\times$ gives the element required in (ii).
\hfill\break
Finally, to prove (iii), observe that the lattice $\Lambda_0=\Cal O^2$ is preserved by $O_B\tt\Cal O$,
and hence is fixed by the element $\lambda$ of part (ii).  If $\Lambda$ is another $\Cal O$--lattice,
preserved by $O_B\tt\Cal O$, then $\Lambda$ must also be fixed by $\lambda$. The whole geodesic
joining the vertices $[\Lambda_0]$ and $[\Lambda]$ in the building of $PGL_2(\kay_p)$
is then fixed by $\lambda$. In particular, $\lambda$ must then fix a vertex at distance
$1$ from $[\Lambda_0]$. But this implies that the image of $\lambda$ in $PGL_2(\F_p)$
has a fixed point on $\Bbb P^1(\F_p) = \Bbb P(\Lambda_0/\pi\Lambda_0)$, and hence
an $\F_p$--rational eigenvector/eigenvalue, which has been excluded.
\qed\enddemo

Returning to $A(p)$ and our isogeny, the isomorphism (10.51) determines an
isomorphism
$$V(A(p)) \simeq \kay^2,\tag10.55$$
under which
$$V(G^+) = \kay\cdot\pmatrix 0\\1\endpmatrix\qquad\text{\rm and}\qquad
V(G^-) = \kay\cdot\pmatrix 0\\1\endpmatrix.\tag10.56$$
The image of the Tate module $T(A(p))$ in $\kay^2$ is an $\Cal O$--lattice
which is stable under the action of $O_B\tt\Cal O$.

If $\kay_p$ is unramified, then $T(A(p))\subset \kay^2$ is
an $\Cal O$--lattice stable under the Eichler
order $O'\subset M_2(\Cal O)$ in (i) of Lemma. Let $1_2 = e_++e_- \in O'$
where
$$e_+ = \pmatrix 0&0\\0&1\endpmatrix \qquad\text{\rm and}\qquad e_- = \pmatrix 1&0\\0&0\endpmatrix\in O'.\tag10.57$$
If $y = y_++y_-\in T(A(p))$ with $y_\pm\in V(G^\pm)$, then $y_\pm = e_\pm y \in V(G^\pm)\cap T(A(p)) = T(G^\pm)$,
and hence $T(A(p))= T(G^+\times G^-)$. Thus, our isogeny has degree 1 and $\delta_v(u) = 0$.

If $\kay_p/\Q_p$ is ramified, then, by (iii) of Lemma~10.5,  $T(A(p))$ is homothetic to $\Lambda_0=\Cal O^2$.
But then, since $T(G^\pm) = T(A(p))\cap V(G^\pm)$, we have simply $T(G^+\times G^-) = T(A(p))$,
so again our isogeny has degree 1 and $\delta_v(u) = 0$. This finishes the proof of (ii) of Proposition~10.2.\qed

We can now compute the Faltings height.

\proclaim{Theorem 10.7} Suppose that the triple $(A,\iota,x)$, defined over a number field $L$, is of type $c$.
Write $4m= n^2 d$, with $-d$ a fundamental discriminant, and let $E=E_d$ be an elliptic curve over $L$
with complex multiplication by $O_\smallkay$, the ring of integers in $\kay=\Q(\sqrt{-d})$.
Then
$$\hfal^*(A) = 2\,\hfal^*(E) +\log(c) - \sum_p\frac{(1-p^{-\ord_p(c)})\cdot (1-\chi(p))}{(1-p^{-1})\cdot(p-\chi(p))}\cdot \log(p).$$
Here $\chi=\chi_d$ is as in (8.3).
\endproclaim

\demo{Proof}
We apply formula (10.18) to the isogeny $E^+\times E^-\rightarrow A$ defined above.
The change in the Faltings height due to the isogeny has the form
$$\frac12\,\log(\deg(u_L)) - \frac1{|L:\Q|}\sum_v \delta_v(u).\tag10.58$$
We write
$$\frac12\,\log(\deg(u_L)) = \sum_p \frac12\ord_p(\deg(u_L))\cdot \log(p),\tag10.59$$
so that (10.58) can be written as a sum of local contributions (10.58)${}_p$ which we now describe case by case.

If $p$ is  split in $\kay$, then by Proposition~10.1,
$$(10.58)_p = \ord_p(c)\cdot \log(p).\tag10.60$$

If $p\nmid D(B)$ is inert or ramified, by Proposition~10.2 , (10.58)${}_p$ is equal to:
$$
\align
&\frac12\ord_p(\deg(u_L))\cdot \log(p)\tag10.61
\\
\nass
&-\frac1{|L:\Q|}\sum_{v\mid p} |L_v:\Q_p|
\bigg( \frac12\big[\ord_p(\deg(u_L)) - 2r\big] +
\frac{(1-p^{-r})\cdot (1-\chi(p))}{(1-p^{-1})\cdot(p-\chi(p))} \bigg)\, \log(p)\\
\nass
\nass
{}&= \bigg(r - \frac{(1-p^{-r})\cdot (1-\chi(p))}{(1-p^{-1})\cdot(p-\chi(p))} \bigg)\, \log(p)
\endalign
$$
where we have set $r=\ord_p(c)$.

If $p\mid D(B)$, then by the considerations after Lemma~10.5,  (10.58)${}_p=0$.
\qed\enddemo

\proclaim{Theorem 10.8} The contribution of the `horizontal' part to
the pairing $\langle\hat\Cal Z(m,v),\hat\omega\rangle$ is
$$\align
h_{\hat\o}(\Cal Z(m)^\hor) &=
2\, \delta(d,D)\,H_0(m,D(B))\,
2\,\hfal^*(E)\\
\nass
\nass
{}&\qquad+2\delta(d,D)\,\frac{h(d)}{w(d)}\sum_{\matrix \scr c|n\\\scr (c,D(B))=1\endmatrix}
c\prod_{\ell|c}(1-\chi(\ell)\ell^{-1})
\cdot
\sum_p\eta_p(\ord_p(c))\log(p),
\endalign$$
where, for $r\in \Z_{\ge0}$,
$$\eta_p(r)=r
-\frac{(1-p^{-r})\cdot (1-\chi(p))}{(1-p^{-1})\cdot(p-\chi(p))}.
$$
\endproclaim
\demo{Proof}
Continuing (10.4) above, we have
$$\align
&h_{\hat\o}(\Cal Z(m)^\hor)\\
\nass
{} &= 2 \sum_{c|n}\sum_{\matrix\scr x\in L(m)\\\scr \mod \Gamma\\ \scr type \ c\endmatrix} \hfal^*(A_x)\cdot\frac1{|\Gamma_x|}\\
\nass
&= 2 \sum_{c|n}\bigg(\sum_{\matrix\scr x\in L(m)\\\scr \mod \Gamma\\ \scr type \ c\endmatrix} \frac1{|\Gamma_x|}\bigg)
\bigg(2\,\hfal^*(E) +\log(c) - \sum_p\frac{(1-p^{-\ord_p(c)})\cdot (1-\chi(p))}{(1-p^{-1})\cdot(p-\chi(p))} \log(p)\bigg)\\
\nass
{}&=2\, \delta(d,D)\,H_0(m,D(B))\,
2\,\hfal^*(E)\\
\nass
\nass
{}&\qquad\qquad\qquad +
2\delta(d,D)\,\frac{h(d)}{w(d)}\sum_{\matrix \scr c|n\\\scr (c,D(B))=1\endmatrix}
c\prod_{\ell|c}(1-\chi(\ell)\ell^{-1})
\cdot
\sum_p\eta_p(\ord_p(c))\log(p),
\endalign
$$
as claimed.
\qed\enddemo

We now make the comparison of this expression with terms arising in the
the positive Fourier coefficients of the derivative of the modified Eisenstein series.
To do this, we need a better expression for the sum on $c$ in the second term in
Theorem~10.8.
For convenience in the calculations, we let
$$\beta_p(k)  = -2k +
\cases
\frac{p^k-1}{p^k(p-1)},
&\text{if $\chi_d(p)=1$,}\\
\nass
\frac{(3p+1)(p^k-1)-4k(p-1)}{(p-1) (p^{k+1}+p^k-2)},
&\text{if $\chi_d(p)=-1$,}\\
\nass
\frac{2}{p-1}-\frac{2k+2}{p^{k+1}-1}&
\text{if $\chi_d(p)=0$.}\endcases\tag10.62
$$
Note that when $k=\ord_p(n)$, then, by (i) of Lemma~8.7,
$$\beta_p(k)= \frac{1}{\log(p)}\cdot\frac{b'_p(n,0;D)}{b_p(n,0;D)}.\tag10.63$$

\proclaim{Lemma 10.9} Let $4m=n^2d$, as before. Then the following
identity holds for any square-free $D>0$:
$$\align
&{h(d)\over w(d)}\sum\limits_{c\vert n\atop (c,D)=1} c\,
\prod\limits_{p\vert c} (1-\chi_d(p)\, p^{-1})\,\sum\limits_{p\vert
c} \eta_p(\ord_p(c))\,\log(p)
\\
\nass
\nass
{}&=
H_0(m;D)\cdot\sum\limits_{p\atop (p,D)=1} \big(\,-\ord_p(n)-\beta_p(\ord_p(n))\,\big)\,\log(p)\\
\nass
\nass
{}&= H_0(m;D)\cdot \sum\limits_{p\atop (p,D)=1} \bigg(\log|n|_p-\frac{b'_p(n,0;D)}{b_p(n,0;D)}\bigg).
\endalign$$
\endproclaim

\demo{Proof}
We note that the sum on the last expression of the Lemma (right hand side) is finite since only summands for $p$ with
$p\vert n$ are non-zero. We proceed by induction on the number of prime factors
of $n$. To start the induction, let $n=p^t$. Then the first expression of the Lemma (left hand side) is equal to
$${h(d)\over w(d)}\cdot\sum\limits_{r=1}^t p^r\, (1-\chi_d(p)p^{-1})\,\eta_p(r)\cdot {\roman{log}}\, p\tag10.64$$
(note that the contribution of $c=1$ is trivial). By (8.19) and (10.62), the right hand side is equal
to
$${h(d)\over w(d)}\cdot\left( \sum_{r=1}^tp^r (1-\chi_d(p)p^{-1})+1\right)
\cdot (-t-\beta_p(t))\cdot {\roman{log}}\, p.\tag10.65$$
Case by case,  one can check that these two expressions coincide.
\par\noindent
{\bf Case $\chi_d(p)=1$:} Then (10.64) without the factor ${h(d)\over
w(d)}\cdot\log(p)$ is equal to
$$\sum_{r=1}^t p^r\, (1-p^{-1})\, r=(p-1)\,\sum_{r=1}^t r\,
p^{r-1} =t\, p^t -{p^t-1\over p-1}.\tag10.66$$
On the other hand, (10.65) without the factor ${h(d)\over w(d)}\cdot\log(p)$ is equal
to
$$\left( \sum_{r=1}^t p^r (1-p^{-1})+1\right)\cdot \left(t-p^{-t}\,
{p^t-1\over p-1}\right) =p^t\, \left(t -p^{-t}\,{p^t-1\over
p-1}\right).\tag10.67$$
\par\noindent
{\bf Case $\chi_d(p)=-1$:} Then (10.64) without the factor ${h(d)\over
w(d)}\cdot\log(p)$ is equal to
$$\align
&\sum_{r=1}^t p^r(1+p^{-1})\, \left( r-2\, p^{-r+1}\, {p^r-1\over
p^2-1}\right)\\
\nass
{} &=(1+p^{-1})\, \sum_{r=1}^t \left( rp^r-2p\,
{p^r-1\over p^2-1}\right)\tag10.68
\\
\nass
\nass
&{}=
{(p-1)\, t\,(p^t(p+1)+2)-(3p+1)\, (p^t-1)\over (p-1)^2}\ \ .
\endalign$$
On the other hand, (10.65) without the factor ${h(d)\over
w(d)}\cdot\log(p)$ is equal to
$$\align
&
\left( \sum_{r=1}^t p^r(1+p^{-1})+1\right) \left( t -{(3p+1)\,
(p^t-1)-4t(p-1)\over (p-1)\,(p^{t+1}+p^t-2)}\right)
\\
\nass
\nass
{}&=
{p^{t+1}+p^t-2\over p-1} \left( t - {(3p+1)(p^t-1)-4t(p-1)\over (p-1)\,
(p^{t+1}+p^t-2)}\right)\tag10.69
\\
\nass
\nass
{}&=
{(p-1)\,t\,(p^t(p+1)+2)-(3p+1)\, (p^t-1)\over (p-1)^2}\ \ .
\endalign$$
\par\noindent
{\bf Case $\chi_d(p)=0$:} Then (10.64) without the factor ${h(d)\over
w(d)}\cdot\log(p)$ is equal to
$$\align
&
\sum_{r=1}^tp^r\, \left( r-p^{-r}\, {p^r-1\over p-1}\right)\\
\nass
{}& =
\sum_{r=1}^t r\, p^r-\sum_{r=1}^t {p^r-1\over p-1}\tag10.70
\\
\nass
{}&=
{p^{t+1}-1\over p-1}\left( t-{2\over p-1} +{2t+2\over p^{t+1}-1}\right)\ \
.
\endalign$$
On the other hand, (10.65) without the factor ${h(d)\over
w(d)}\cdot\log(p)$ is equal to
$$\left( \sum_{r=0}^t p^r\right)\, \left( t-{2\over p-1}+ {2t+2\over
p^{t+1}-1}\right) = {p^{t+1}-1\over p-1} \left( t-{2\over p-1}+
{2t+2\over p^{t+1}-1}\right)\ \ .\tag10.71$$
We therefore have checked the beginning of the induction. Let us now perform the
induction step. Let $n=p^t\cdot n_0$ where $p\nmid n_0$.
Let us put $m_0=m/p^{2t}$, so that $4m_0=n_0^2d$. We may assume that
$p\nmid D$ because otherwise both sides of the identity for $n$
coincide with the corresponding sides of the identity for $n_0$, so that
we may apply the induction hypothesis. We write ${\Cal L}(m)$ resp.\ ${\Cal
R}(m)$ for the left hand side resp.\ right hand side of our identity corresponding to $m$. Then
${\Cal L}(m)$ is equal to
$$\align
{h(d)\over w(d)}\sum_{r=1}^t p^r
(1-\chi_d(p)p^{-1})&\sum_{c_0\vert n_0\atop (c_0,D)=1}
c_0\cdot\bigg(\prod_{\ell \vert c_0} (1-\chi_d(\ell)\ell^{-1})\bigg)\,
\bigg[\,\eta_p(r)\log(p) +
\sum_{\ell\mid c_0}\eta_{\ell} (r_{\ell}(c_0))\log(\ell)\,\bigg]
\tag10.72
\\
\nass
{}+{h(d)\over w(d)} &\sum_{c_0\vert n_0\atop (c_0, D)=1} c_0\cdot
\bigg(\prod_{\ell\vert c_0} (1-\chi_d(\ell)\ell^{-1})\bigg)\cdot
\sum_{\ell\vert c_0}\eta_{\ell} (r_{\ell} (c_0))\,\log(\ell)
\endalign$$
We recall that
$${h(d)\over w(d)}\cdot\sum_{\matrix\scr c_0\vert n_0\\\scr (c_0,D)=1\endmatrix} c_0\prod_{\ell\vert
c_0}(1-\chi_d(\ell)\ell^{-1})= H_0(m_0;D).\tag10.73$$
Hence we can write the above expression as a sum of three terms, the first
one being
$$\align
&
{h(d)\over w(d)}\sum_{c_0\vert n_0\atop (c_0, D)=1} c_0\,
\bigg(\prod_{\ell\vert c_0} (1-\chi_d(\ell)\ell^{-1})\bigg)\cdot
(1-\chi_d(p)p^{-1})\,\sum_{r=1}^t p^r\, \eta_p(r)\log(p)
\tag10.74
\\
=
&
H_0(m_0, D)\cdot (1-\chi_d(p)p^{-1})\,\sum_{r=1}^t
p^r\,\eta_p(r)\,\log(p)\ \ .
\endalign$$
The second and the third term are respectively equal to
$${\Cal L}(m_0)\cdot (1-\chi_d(p)\cdot p^{-1})\,\sum_{r=1}^t p^r =
{\Cal L} \left( m_0\right)\cdot (1-\chi_d(p)p^{-1})\, p\,
{p^t-1\over p-1}\tag10.75$$
and ${\Cal L}(m_0)$.
\par\noindent
We thus obtain
$${\Cal L}(m)= {p^{t+1}-\chi_d(p)p^t+\chi_d(p)-1\over
p-1}\cdot{\Cal L}(m_0)
+H_0(m_0,D)\cdot (1-\chi_d(p)p^{-1})\, \sum_{r=1}^t
p^r\,\eta_p(r)\log(p)\ \ .\tag10.76$$
By induction we have for the last summand
$$\align
&H_0(m_0;D)\cdot
(1-\chi_d(p)p^{-1})\,\sum_{r=1}^tp^r\eta_p(r)\,\log(p)\\
\nass
{}& =H_0(m_0;D)\cdot \bigg( \sum_{r=1}^tp^r\, (1-\chi_d(p)p^{-1})+1\bigg)\big(\,-t-\beta_p(t)\,\big)\,\log(p)\tag10.77\\
\nass
\nass
{}&= H_0(m;D)\cdot\big(\,-t-\beta_p(t)\,\big)\,\,\log(p)\ \ .
\endalign$$
Hence
$${\Cal L}(m)= {p^{t+1}-\chi_d(p)p^t+\chi_d(p)-1\over p-1}\cdot {\Cal
L}(m_0)+ H_0(m;D)\cdot (-t-\beta_p(t))\,\log(p).
\tag10.78$$
Now recall from (8.13) and Lemma~8.5 that
$$H_0(m;D)= {h(d)\over w(d)} \prod_{q\nmid D}
{q^{t+1}-\chi_d(q)q^t+\chi_d(q)-1\over q-1}\ \ .\tag10.79$$
 It follows that
$${H_0(m;D)\over H_0(m_0;D)} = {p^{t+1}-\chi_d(p)p^t+\chi_d(p)-1\over
p-1}\ \ .\tag10.80$$
From the definition of ${\Cal R}(m)$ we have $$\align
{\Cal R}(m) & ={H_0(m,D)\over H_0(m_0; D)}\cdot {\Cal
R}(m_0)+H_0(m,D)\cdot (-t-\beta_p(t))\,\log(p) \tag10.81\\ \nass
{}&={p^{t+1}-\chi_d(p)p^t+\chi_d(p)-1\over p-1}\cdot {\Cal R}
(m_0)+H_0(m;D)\cdot (-t-\beta_p(t))\,\log(p)\ \ .
\endalign$$
Comparing (10.78) with (10.81), the induction hypothesis ${\Cal
L}(m_0)={\Cal R}(m_0)$ implies the assertion.
\qed
\enddemo

The following result is well known, cf., for example, \cite{\colmez}.
\proclaim{Proposition 10.10} With the normalization given by (10.14) above,
the Faltings height $\hfal(E)$ of an elliptic curve $E$ with CM by $O_{\smallkay}$ is given by
$$\align
2\,\hfal(E) &= -\frac12\log(d)-\frac{L'(0,\chi_d)}{L(0,\chi_d)}\\
\nass\phantom{-}
{}&= \phantom{-}\frac12\,\log(d) - \frac{w(d)}{2 h(d)}\sum_{a=1}^{d-1} \chi_d(a)\,\log\Gamma\left(\frac{a}{d}\right)\\
\nass
{}&=\phantom{-}\frac12\log(d)+\frac{L'(1,\chi_d)}{L(1,\chi_d)} - \log(2\pi) -\gamma.
\endalign
$$
\endproclaim

{\bf Remark 10.11.} The value for $2\hfal(E)$ in Colmez \cite{\colmez},
p.633 is our $2\hfal(E) - \log(2\pi)$ due to a difference in the normalization of the metric on the Hodge bundle.

Our {\it renormalized} Faltings height is then given by
$$\align
2\,\hfal^*(E) &= 2\,\hfal(E) + \frac12\log(\pi) +\frac12\gamma+\log(2)\tag10.82\\
\nass
{}&= \frac12\log(d)+\frac{L'(1,\chi_d)}{L(1,\chi_d)}- \frac12\log(\pi) -\frac12\gamma\\
\endalign$$

Combining these facts, we have
\proclaim{Corollary 10.12} The contribution of the `horizontal' part to
the pairing $\langle\,\hat\Cal  Z(m,v),\hat\omega\,\rangle$ is
$$\align
h_{\hat\o}(\Cal Z(m)^\hor) =
2\, \delta(d;D)\,H_0(m;D)\,&
\bigg[\,\frac12\log(d)+\frac{L'(1,\chi_d)}{L(1,\chi_d)} - \frac12\log(\pi) -\frac12\gamma
\\
\nass
\nass
\nass
{}&\qquad\qquad\qquad+ \sum\limits_{p\atop p\nmid D=1} \bigg(\log|n|_p-\frac{b'_p(n,0;D)}{b_p(n,0;D)}\bigg)\ \bigg].
\endalign$$
\endproclaim

\demo{Proof of Theorem~7.2} Looking back to the end of section 8, we see that the expression
of Corollary~10.12 for $h_{\hat\o}(\Cal Z(m)^\hor)$ coincides exactly with the
sum of (8.43) and (8.45). The remaining terms will be considered in the next two sections.
\enddemo

\subheading{\Sec11. Contributions of vertical components}

In this section we fix a prime number $p$ with $p\mid D(B)$. We wish to
determine the quantity $\deg(\omega\vert\Cal Z(m)_p^{\roman{vert}})$, cf.\
(9.11), using the results of \cite{\krinvent}.

We describe $\Cal Z(m)\times_{\Spec\, \Z}\Spec\, W(\bar\F_p)$ in terms of the
$p$-adic uniformization of ${\Cal M}\times_{\Spec\, \Z}\Spec(\F_p)$,
comp.\ section 2. To this end we fix $x\in O_{B'}$ with
${\roman{tr}}^\circ(x)=0$ and $x^2=-m$. As in section 2 we identify
$B'\otimes \A_f^p$ with $B\otimes \A_f^p$ and $H(\A_f^p)$ with
$H'(\A_f^p)$ and $K^p$ with $K^{\prime p}$. Put
$$I(x)=\{ gK^p\in H'(\A_f^p)/K^{\prime,p}\mid\ g^{-1}xg\in\hat O_{B'}^p\}\ \
.\leqno(11.1)$$
We also use the abbreviation $\hat\Omega_{W(\bar\F_p)}$ for
$\hat\Omega^2\times_{{\roman{Spf}}\,\Z_p}{\roman{Spf}}\, W(\bar\F_p)$. Let
$\tilde x = x$, if $\ord_p(m)=0$ (resp. $\tilde x = 1+x$, if $\ord_p(m)>0$).
Let
$$Z(x)=(\hat\Omega_{W(\bar\F_p)}\times \Z)^{\tilde x}\leqno(11.2)$$
be the fixed point set of $\tilde x\in H'(\Q_p)$. Denoting by $H'_x$
the stabilizer of $x$ in $H'$, we have, \cite{\krinvent},
$$\Cal Z(m)\times_{\Spec\, \Z}\Spec\, W(\bar\F_p)=[H'_x(\Q)\setminus I(x)\times
Z(x)]\tag11.3$$
(quotient in the sense of stacks). Since $\ord_p\det(\tilde x)=0$,
we have
$$Z(x)=\hat\Omega^{\tilde x}_{W(\bar\F_p)}\times\Z\ \ .\leqno(11.4)$$
Since the set
$$H'_x(\A_f^p)\setminus\{ g\in H'(\A_f^p)\mid\ g^{-1}xg\in\hat O_{B'}^p\}$$
is compact, the group $H'_x(\A_f^p)$ has only finitely many orbits on
$I(x)$. Let $g_1,\ldots, g_r\in H'(\A_f^p)$ such that
$$I(x)= \coprod^r_{i=1}H'_x(\A_f^p)\, g_i\, K^{\prime, p}\ \
.\leqno(11.5)$$
Then we may rewrite (11.3) as
$$\coprod^r_{i=1}\left[H'_x(\Q)\setminus \left(\,H'_x(\A_f^p)/ (K_i^{\prime, p}\cap
H'_x(\A_f^p))\times \Z\times \hat\Omega^{\tilde x}_{W(\bar\F_p)}\,\right)\right]\ \ ,$$
where $K_i^{\prime p}=g_iK^{\prime p}g_i$.
\par
Note that $H'_x(\Q)\cong\kay^\times$, where $\kay=\Q(\sqrt{-m})$ is the imaginary
quadratic field associated to $m$. Let us first consider the case where
$p$ does not split in $\kay$. Then
$$\ord_p\det(\kay_p^\times)=\delta_p\cdot\Z\ \ ,\leqno(11.6)$$
where $\delta_p=2$ if $p$ is unramified in $\kay$ and $\delta_p=1$ if $p$ is
ramified in $\kay$. Let
$$H'_x(\Q)^1=\{ g\in H'_x(\Q)\mid\ \ord_p(\det(g)) =0\}\ \ .\leqno(11.7)$$
Then $H'_x(\Q)^1$ acts with finite stabilizer groups on
$H'_x(\A_f^p)/(K_i^{\prime, p}\cap H'_x(\A_f^p))$. Hence we may rewrite
(11.3) as $\delta_p$ copies of
$$\coprod^r_{i=1}\left[H'_x(\Q)^1\setminus \left(\,H'(\A_f^p)/(K_i^{\prime, p}\cap
H'_x(\A_f^p))\,\right)\right]\times \hat\Omega^{\tilde x}_{W(\bar\F_p)}\leqno(11.8)$$
(here the first factor is taken in the sense of stacks).
\par
Appealing now to \cite{\krinvent}, Proposition~3.2, we obtain the following expression
for the vertical
components of $\Cal Z(m)$:
$$\align
&\Cal Z(m)^{\roman{vert}}\times_{\Spec\,\Z}\Spec\, W(\bar\F_p)
\tag11.9\\
\nass
&
=\Z/\delta_p\Z
\times \coprod^r_{i=1}\left[\,H'_x(\Q)^1 \setminus \left(\,H'(\A_f^p)/
(K_i^{\prime, p}\cap H'_x(\A_f^p))\,\right)\right]
\times \bigg(\sum_{[\Lambda]\in{\Cal
B}}{\roman{mult}}_{[\Lambda]}(x)\cdot{\Bbb P}_{[\Lambda]}\bigg).
\endalign$$
Here $[\Lambda]$ ranges over the vertices of the Bruhat-Tits tree of
$PGL_2(\Q_p)$ and the multiplicity with which the prime divisor ${\Bbb
P}_{[\Lambda]}$ occurs is given by loc.cit., (3.9) for $p\ne2$ and by Proposition~A.1
in the Appendix below for $p=2$.

\proclaim{Proposition 11.1}
For any $[\Lambda]\in{\Cal B}$ we have
$$\deg(\omega\vert{\Bbb P}_{[\Lambda]})=p-1\ \ .$$
\endproclaim

\demo{Proof}
Of course, here $\deg(\omega\vert{\Bbb P}_{[\Lambda]})$ is shorthand for
$$\deg\, i^*_{[\Lambda]}(\omega\otimes_{\Z}W(\bar\F_p))\ \ ,$$
where $i_{[\Lambda]}:{\Bbb P}_{[\Lambda]}\to {\Cal M}\times_{\Spec\,
\Z}\Spec\, W(\bar\F_p)$ is the natural morphism. We write $O_{B_p}$ as
$$O_{B_p}=\Z_{p^2} [\Pi]/(\Pi^2=p,\ \Pi a=a^{\sigma}\Pi,\ \ \forall a\in
\Z_{p^2})\ \ .$$
For the inverse image of the universal abelian scheme $(\Cal A,\iota)$ on
${\Cal M}\times_{\Spec\, \Z}\Spec \, W(\bar\F_p)$ we have
$$\Lie\,\Cal A={\Cal L}_0\oplus {\Cal L}_1\ \ ,\leqno(11.10)$$
where ${\Cal L}_i=\{ x\in \Lie\,\Cal A\mid\ \iota(a)x=a^{\sigma^{-i}}x,\ \
\forall a\in\Z_{p^2}\}$.
\par
Due to the determinant condition (1.1), both ${\Cal L}_0$ and ${\Cal
L}_1$ are line bundles on ${\Cal M}\times_{\Spec\, \Z}\Spec\, W(\bar\F_p)$ and
$$\omega\otimes_{\Z}W(\bar\F_p)={\Cal L}_0^{-1}\otimes {\Cal L}_1^{-1}\ \
.\leqno(11.11)$$
The fiber of ${\Cal L}_i$ at a $\bar\F_p$-valued point of ${\Cal M}$ is
expressed as follows in terms of the Dieudonn\'e module $M$ of the
corresponding abelian variety,
$${\Cal L}_0=M_0/VM_1\ \ ,\ \ {\Cal L}_1= M_1/VM_0\ \ .\leqno(11.12)$$
Here $M=M_0\oplus M_1$ is the eigenspace decomposition under the action of
$\Z_{p^2}$ analogous to (11.10).
\par
To fix ideas assume that $\Lambda$ is even (\cite{\krinvent}). Then for every
$x\in {\Bbb P}_{[\Lambda]}(\bar\F_p)$ we have
$$M_0= \Lambda\otimes_{\Z_p}W(\bar\F_p),\ VM_0= \Pi M_0;\ {\Cal
L}_{0x}=M_0/\ell_x\ \ ,\leqno(11.13)$$
where $\ell_x$ is the line in $\Lambda\otimes_{\Z_p}\bar\F_p$ corresponding
to $x$. It follows that
$$i^*_{[\Lambda]}({\Cal L}_0)= {\Cal O}_{{\Bbb P}_{[\Lambda]}}(1)\ \
.\leqno(11.14)$$
(It is ${\Cal O}_{{\Bbb P}_{[\Lambda]}}(1)$ rather than ${\Cal O}_{{\Bbb
P}_{[\Lambda]}}(-1)$ since ${\Cal L}_0$ obviously has global sections.) To
calculate $i^*_{[\Lambda]}({\Cal L}_1)$, we use the exact sequence
of coherent sheaves on $\hat\Omega_{W(\bar\F_p)}$
$$0\lra {\Cal L}_1\buildrel\Pi\over\lra {\Cal L}_0\lra {\Cal O}_{{\Bbb
P}^{\roman{odd}}}\lra 0\ \ .\leqno(11.15)$$
Here ${\Cal O}_{{\Bbb P}^{\roman{odd}}}$ denotes the structure sheaf of
the closed subscheme of the special fiber (\cite{\krinvent}, section 2),
$$\bigcup_{[\Lambda]\ {\roman{odd}}}{\Bbb P}_{[\Lambda]}\subset
\hat\Omega_{W(\bar\F_p)}\otimes_{W(\bar\F_p)}\bar\F_p\ \ .\leqno(11.16)$$
This sequence remains exact after pulling back and yields
$$0\lra i^*_{[\Lambda]}({\Cal L}_1)\lra i^*_{[\Lambda]}({\Cal L}_0)\lra
{\Cal O}_{{\Bbb P}_{[\Lambda]}(\F_p)}\lra 0\leqno(11.17)$$
(note that $i^*_{[\Lambda]}({\Cal L}_1)$ is torsion-free). Here we used
that $ {\Bbb P}_{[\Lambda]}(\F_p)={\Bbb P}_{[\Lambda]}\cap
(\bigcup\limits_{[\Lambda]\ {\roman{odd}}}{\Bbb P}_{[\Lambda]})$. From
(11.14) we obtain
$$\align
\deg\, i^*_{[\Lambda]}({\Cal L}_1)
&
=\deg\, i^*_{[\Lambda]}({\Cal L}_0)-(p+1)\tag11.18\\
&
=-p\ \ .
\endalign$$
Now the identity (11.11) yields the assertion. The case where $[\Lambda]$
is odd is similar.
\qed
\enddemo
\medskip\noindent
{\bf Remark 11.2.} Another proof of Proposition 11.1 may be obtained by
using Proposition 3.2. Indeed, by that proposition we may identify $\omega$ and the relative dualizing sheaf $\omega_{{\Cal
M}/\Z}$. It follows that
$$\deg(\omega\vert{\Bbb P}_{[\Lambda]})=\deg (\omega_{{\Cal M}/\Z}\vert\,
{\Bbb P}_{[\Lambda]})=\deg(\omega_{{\Cal M}\otimes_{\Z}\bar\F_p/\bar\F_p}\vert\,
{\Bbb P}_{[\Lambda]})\ \ .$$
By expressing the dualizing sheaf $\omega_{{\Cal
M}\otimes_{\Z}\bar\F_p/\bar\F_p}$ explicitly, it is easy to calculate the last term.
\qed

\proclaim{Corollary 11.3} Let $k=\ord_p n$, where, as usual, $4m=n^2d$. Then
$$\sum_{[\Lambda]\in{\Cal
B}}\text{\rm mult}_{[\Lambda]}(x)\cdot\deg(\omega\vert{\Bbb
P}_{[\Lambda]}) =
\cases
-2k+(p+1)\cdot {p^{k}-1\over p-1} &\text{ if $p$ is unramified in $\kay$}\\
\nass
-2k-2+2\cdot {p^{k+1}-1\over p-1} &\text{  if $p$ is ramified in $\kay$.}
\endcases
$$
\endproclaim

\demo{Proof}
It remains to calculate $\sum_{[\Lambda]}\text{\rm mult}_{[\Lambda]}(x)$. If $p$ is odd, this is an
easy exercise using the results of \cite{\krinvent}, section 6. For instance, let  $p$ be odd and
unramified and put $\a=\ord_p(m)$ so that $\a=2k$. Then
$$\sum_{[\Lambda]}\text{\rm mult}_{[\Lambda]}(x)= {\alpha\over 2}
+(p+1)\sum_{r=1}^{{\alpha\over 2}-1} p^{r-1}\left( {\alpha\over
2}-r\right) ={-\alpha\over p-1}+{p+1\over p-1} \cdot {p^{\alpha\over
2}-1\over p-1}\ .\eqno$$
The case $p=2$ is handled in the Appendix to this section. \qed
\enddemo
\medskip\noindent
It now remains to determine the degree of the discrete stack
$$\coprod\limits_{i=1}^r\left[H'_x(\Q)^1\setminus \left( H'(\A_f^p)/(K_i^{\prime, p}\cap
H'_x(\A_f^p))\right)\right].$$
 Now the elements $g_1,\ldots, g_r$ are in one-to-one
correspondence with the $\hat{O}_B^{p,\times}$--conjugacy classes of
embeddings of rings
$$j^p:\kay\otimes \A_f^p\lra B\otimes \A_f^p,$$
with $x\in (j^p)^{-1}(\hat{\Cal O}_B^p)$. To each embedding $j^p$ there is
associated an order of $\kay$,
$$O(j^p)=((j^p)^{-1}(\hat O_B^p).O_{\smallkay_p})\cap \kay.$$
The conductor $c = c(j^p)$ of this order satisfies $(c,D(B))=1$ and $c\mid n$, because $x\in
(j^p)^{-1}(O_B^p)$. Conversely, given $c$ with those two properties,
there are precisely
$$\prod_{\ell\vert D(B), \ell\neq p}(1-\chi_d(\ell))$$
classes of embeddings $j^p$ yielding the order of conductor $c$. Finally,
if $g_i$ yields the order of conductor $c$, the stack $[H'_x(\Q)^1\setminus
(H'(\A_f^p)/ K'_i\cap H'_x(\A_f^p))]$ may be identified with
$$[\kay^{\times,1}\setminus (\kay\otimes \A_f^p/\hat{O}_{c^2d}^{p,\times})]$$
 which has degree
$h(c^2d)/w(c^2d)$. Summarizing these arguments, we therefore obtain

\proclaim{Lemma 11.4}
$$\align
&
\delta_p\cdot\deg \left( \coprod\limits_{i=1}[H'_x(\Q)^1\setminus
H'(\A_f^p)/K_i^{\prime, p}\cap H'_x(\A_f^p)]\right)
\\
\nass
&
=\prod\limits_{\ell\vert D(B)} (1-\chi_d(\ell))\cdot\sum\limits_{c\mid n\atop
(c,D(B))=1} h(c^2d)/w(c^2d)
\\
\nass
&
=
\delta(d;D(B))\cdot H_0(m;D(B))\ \ .
\endalign$$
\rightline{\hfill\qed}
\endproclaim

Now let us consider the case when $p$ splits in $\kay$. In this case
$\ord_p(\det(\kay_p^{\times}))=\Z$, but $H'_x(\Q)^1$ does not act with finite
stabilizer groups on $H'_x(\A_f^p)/(K_i^{\prime p}\cap H'_x(\A_f^p))$. Let
$\epsilon(x)\in H'_x(\Q)=\kay^{\times}$ be an element whose localization in
$\kay_p=\Q_p\oplus \Q_p$ has valuation $(1,-1)$. Let $H'_x(\Q)^{1,1}$ be the subgroup of
elements of $H'_x(\Q)$ which are units at $p$. Then
$H'_x(\Q)^1=H'_x(\Q)^{1,1}\times \langle\epsilon(x)^{\Z}\rangle$, and
$H'_x(\Q)^{1,1}$ acts with finite stabilizer groups on
$H'_x(A_f^p)/(K_i^{\prime p}\cap H'_x(\A_f^p))$, whereas $\epsilon(x)$
acts freely on $\hat\Omega^{\tilde x}_{W(\bar\F_p)}$ by translations by 2 on
the `apartment of central components' (\cite{\krinvent}). We obtain therefore
the following expression for (11.3) in this case
$$\left( \coprod\limits_{i=1}^r \left[H'_x(\Q)^{1,1}\setminus \left(\,H'_x(\A_f^p)
/(K_i^{\prime, p}\cap H'_x(\A_f^p))\, \right)\right]\right) \times \left( \langle\epsilon
(x)^{\Z}\rangle\setminus \hat\Omega^{\tilde x}_{W(\bar\F_p)}\right)
.\leqno(11.19)$$
The same analysis as before yields
$$\align
&\deg(\,\coprod\limits^r_{i=1}\left[H'_x(\Q)^{1,1}\setminus\left(\, H'_x(\A_f^p)/(K_i^{\prime,p} \cap H'_x(\A_f^p))\, \right)\right]\,)
\tag11.20
\\
\nass
&
=\prod\limits_{\ell\vert D(B)\atop \ell\neq p}
(1-\chi_d(\ell))\cdot\sum\limits_{c\vert n\atop (c,D(B))=1}
h(c^2d)/w(c^2d)
\\
\nass
&
=
\delta(d;D(B)/p)\cdot H_0(m;D(B))\ \ .
\endalign$$
Using Proposition 11.1 we have
$$
\deg(\o\mid\big(\langle\epsilon (x)^{\Z}\rangle\setminus \hat\Omega^{\tilde
x}_{W(\bar\F_p)}\big)\,)
=2\, (p-1)\,\sum_{[\Lambda]} \text{\rm mult}_{[\Lambda]}(x).
\tag11.21
$$
The sum on the right hand side runs over all vertices $[\Lambda]$ such that the
closest vertex on the apartment corresponding to $\kay^\times_p$ is a fixed vertex.
This sum can again be evaluated using \cite{\krinvent}, 3.9, for $p$ odd (resp. the
appendix to this section for $p=2$).

We summarize our findings in the following theorem.

\proclaim{Theorem 11.5}
Let $k=\ord_p(n)$, where, as usual, $4m= n^2d$. \hfill\break
(i) If $p$ splits in $\kay$, then
$$\deg(\,\omega\mid \Cal Z(m)_p^{\roman{vert}}\,)= 2\,H_0(m;D(B))\,\delta(d;D(B)/p)\cdot (p^{k}-1).$$
(ii)
$$
\deg(\,\omega\mid \Cal Z(m)_p^{\roman{vert}}\,)=2\,H_0(m;D(B))\,\delta(d;D(B))\cdot\cases
\, -k+ \frac{(p+1)(p^{k}-1)}{2(p-1)} &\text{if $\chi_d(p)=-1$,}\\
\nass
-1-k+\,\frac{p^{k+1}-1}{p-1}&\text{if $\chi_d(p)=0$.}
\endcases
$$
\endproclaim

\demo{Proof of Theorem~7.2 (continued)} In the case $\chi_d(p)=1$, the quantity
$$\deg(\,\omega\mid \Cal  Z(m)_p^{\roman{vert}}\,)\,\log(p)\cdot q^m$$
coincides exactly with the term (ii) of Theorem~8.8.
On the other hand, in the cases $\chi_d(p)=-1$ and $\chi_d(p)=0$, we find that
$$\deg(\,\omega\mid \Cal Z(m)_p^{\roman{vert}}\,)\,\log(p)\cdot q^m
=2\,\delta(d;D)\,H_0(m;D)\,K_p\,\log(p)\cdot q^m,$$
where $K_p$ is as in Theorem~8.8.  Thus, summing on $p\mid D$, we obtain (8.45).
\enddemo

\define\notmid{\mkern-5mu\not\mkern5mu\mid}

\subheading{Appendix to section 11: The case $p=2$}

In \cite{\krinvent} we made the blanket assumption $p\neq 2$. In this
appendix we indicate the modifications needed to arrive at the formulas
given in Theorem 11.5 in the case $p=2$.
\par
We will use the same notation as in \cite{\krinvent}. We fix a special
endomorphism $j\in V$ with $q(j)=j^2\in\Z_2\setminus\{ 0\}$. We denote by
$Z(j)$ the associated closed formal subscheme of the Drinfeld moduli space
${\Cal M}\simeq \hat\Omega\times_{\Spf\ \Z_2}\Spf\ W(\bar\F_2)$. We will
content ourselves with giving the structure of the divisor
$Z(j)^{\roman{pure}}$ associated to $Z(j)$, loc.cit., section 4. Our
discussion will proceed by distinguishing cases. Let $\kay=\Q_2(j)$ (hence in
the global case $\kay$ is the localization at 2 of the imaginary quadratic
field). Let ${\Cal O}={\Cal O}_{\smallkay}$ be the ring of integers in $\kay$. We write
as usual
$$q(j)=\varepsilon\cdot 2^{\alpha}\ \ ,\ \ \varepsilon\in\Z_2^{\times}\ \
,\ \ \alpha \ge 0\ \ .\leqno(A.1)$$
We define $k\geq 0$ by
$$\alpha+2= 2k+\ord_2(d)\ \ ,\leqno(A.2)$$
where $d$ denotes the discriminant of ${\Cal O}/\Z_2$. Note that in the
global context, when $\Cal Z(m)$ is $p$-adically uniformized by $Z(j)$ (cf.\ (11.3)
above), then
$\alpha=\ord_p(m)$. If we write as usual $4m=n^2d$, then $k=\ord_p(n)$.
\par
We have then the following cases

\catcode`\_=11%
\def\my_self#1{\def#1{\noexpand#1}}%
\my_self\Uebersetze_Ende
\my_self\Uebersetzung_Ende
\newtoks\Spalten_Format
\newtoks\Spalteneintrag_Format
\newtoks\Spaltenzwischenraum_Format
\newtoks\Spaltenende_Format
\newtoks\Vielfach_Format
\newif\if_vor_erster_Spalte_
\def\_vor_erster_Spalte_true{\global\let\if_vor_erster_Spalte_=\iftrue}%
\def\_vor_erster_Spalte_false{%
 \global\let\if_vor_erster_Spalte_=\iffalse
}%
\newif\if_Mathe_
\newif\if_Mathe_normal_
\newcount\ne_Zahl
\newcount\noch_ne_Zahl
\newdimen\tabrulewidth
\newdimen\hfree
\newdimen\vfree
\newskip\Init_tabskip
\def\Frei_raum#1{%
 \vbox{%
  \vskip\vfree
  \vtop{%
   \hbox{%
    \hskip\hfree
    \ignorespaces
    #1%
    \unskip
    \hskip\hfree
   }%
   \vskip\vfree
  }%
 }%
}%
\def\wider#1{%
 \relax
 \ifmmode
  \def\next{\Frei_raum{$#1$}}%
 \else
  \def\next{\Frei_raum{#1}}%
 \fi
 \next
}%
\tabrulewidth=.4pt
\hfree=2pt
\vfree=2pt
\newskip\par_width
\def\par_box#1#2{%
 \par_width=#1%
 \hbox to \par_width{%
  \vtop{%
   \normalbaselines
   \noindent
   \leftskip=0pt%
   \rightskip=-\par_width
   \advance\rightskip by\hsize
   #2%
  }%
  \hss
 }%
}%
\def\verknuepfe_Tokens_zur_ner_Liste#1#2#3{\global #3={#1#2}}%
\def\fuege_Tokens_an#1#2#3{%
 \expandafter\verknuepfe_Tokens_zur_ner_Liste
 \expandafter{\the#2}{#1}#3%
}%
\def\verknuepfe_Token_Listen#1#2#3{%
 \expandafter\fuege_Tokens_an
 \expandafter{\the#2}#1#3%
}%
\def\fuege_Tokens_ans_Spaltenformat#1{%
 \fuege_Tokens_an
 {#1}%
 \Spalten_Format
 \Spalten_Format
}%
\def\fuege_Tokens_ans_Vielfachformat#1{%
 \fuege_Tokens_an
 {#1}%
 \Vielfach_Format
 \Vielfach_Format
}%
\def\Spalteneintrag_hfil{\fuege_Tokens_ans_Spaltenformat\hfil}%
\def\Mathe_zu_dollar{%
 \if_Mathe_
  \fuege_Tokens_ans_Spaltenformat{\relax$\relax}%
 \fi
}%
\def\Spalten_Eintrag{%
 \Mathe_zu_dollar
 \verknuepfe_Token_Listen
  \Spalten_Format
  \Spalteneintrag_Format
  \Spalten_Format
 \Mathe_zu_dollar
}%
\def\Uebersetze_#1{%
 \def\Uebersetzung_##1#1##2##3\Uebersetzung_Ende{##2}%
 \Uebersetzung_
  +\Uebersetze_plus
  |\Uebersetze_vrule
  *\Uebersetze_star
  @\Uebersetze_at
  c\Uebersetze_c
  d\Uebersetze_d
  h\Uebersetze_h
  l\Uebersetze_l
  m\Uebersetze_m
  p\Uebersetze_p
  r\Uebersetze_r
  \Uebersetze_Ende\empty
  0\Uebersetze_zero
  \Uebersetzung_Ende
}%
\def\Uebersetze_Eintrag#1{%
 \if_vor_erster_Spalte_
  \_vor_erster_Spalte_false
 \else
  \verknuepfe_Token_Listen
   \Spalten_Format
   \Spaltenzwischenraum_Format
   \Spalten_Format
 \fi
 #1%
 \Uebersetze_
}%
\def\Uebersetze_l{%
 \Uebersetze_Eintrag{%
  \Spalten_Eintrag
  \Spalteneintrag_hfil
 }%
}%
\def\Uebersetze_r{%
 \Uebersetze_Eintrag{%
  \Spalteneintrag_hfil
  \Spalten_Eintrag
 }%
}%
\def\Uebersetze_c{%
 \Uebersetze_Eintrag{%
  \Spalteneintrag_hfil
  \Spalten_Eintrag
  \Spalteneintrag_hfil
 }%
}%
\def\Uebersetze_vrule{%
 \fuege_Tokens_ans_Spaltenformat\vline
 \Uebersetze_
}%
\def\vline{%
 \vrule width\tabrulewidth\relax
 \futurelet\next\v_line
}%
\def\v_line{%
 \ifx\next\vline
  \hskip\hfree\relax
 \fi
}%
\def\Uebersetze_at#1{%
 \fuege_Tokens_ans_Spaltenformat{#1}%
 \Uebersetze_
}%
\def\Uebersetze_p#1{%
 \Uebersetze_Eintrag{%
  \expandafter\Uebersetze_p_hilf
  \expandafter{\the\Spalteneintrag_Format}%
  {#1}%
 }%
}%
\def\Uebersetze_p_hilf#1#2{%
 \fuege_Tokens_ans_Spaltenformat
  {\wider{\par_box{#2}{\un_wider#1}}\hfil}%
}%
\def\un_wider{\futurelet\next\skip_wider}%
\def\skip_wider{%
 \ifx\next\wider
  \let\next\eat_wider
 \else
  \let\next\relax
 \fi
 \next
}%
\def\eat_wider#1#2{\un_wider#2}%
\def\Uebersetze_d{%
 {%
  \Uebersetze_ r0@{\hskip-2\hfree}l\Uebersetze_Ende
 }%
 \Uebersetze_
}%
\def\Uebersetze_star#1#2{%
 {%
  \ne_Zahl=#1%
  \loop
  \ifnum\ne_Zahl>0
   \advance \ne_Zahl by-1
   \Uebersetze_ #2\Uebersetze_Ende
  \repeat
 }%
 \Uebersetze_
}%
\def\Uebersetze_zero{\Uebersetze_tabskip{0pt}}%
\def\Uebersetze_plus{\Uebersetze_tabskip{0pt plus 1fil}}%
\def\Uebersetze_tabskip#1{%
 \if_vor_erster_Spalte_
  \Init_tabskip=#1\relax
 \else
  \fuege_Tokens_ans_Spaltenformat{\tabskip=#1\relax}%
 \fi\Uebersetze_
}%
\def\Uebersetze_h{\expandafter\global\_Mathe_false\Uebersetze_}%
\def\Uebersetze_m{\expandafter\global\_Mathe_true\Uebersetze_}%
\def\Beginn_der_Tabelle#1#2#3#4{%
 \_vor_erster_Spalte_true
 \global\Spalten_Format={#1}%
 \if_Mathe_normal_
  \Uebersetze_ 0m#4\Uebersetze_Ende
 \else
  \Uebersetze_ 0h#4\Uebersetze_Ende
 \fi
 \verknuepfe_Token_Listen
  \Spalten_Format
  \Spaltenende_Format
  \Spalten_Format
 \v_box\bgroup
  \offinterlineskip
  \tabskip=\Init_tabskip
  #2%
  \halign #3\bgroup
   \span\the\Spalten_Format
}%
\def\Ende_der_Tabelle{\crcr\egroup\egroup}%
\def\beginmatrix{%
 \global\let\v_box=\vcenter
 \expandafter\global\_Mathe_true
 \global\Spalteneintrag_Format={\wider{##}}%
 \global\Spaltenzwischenraum_Format={&}%
 \global\Spaltenende_Format={\cr}%
 \null
 \,%
 \_Mathe_normal_true
 \Beginn_der_Tabelle{}{\mathsurround=0pt}{}%
}%
\def\endmatrix{\Ende_der_Tabelle\,}%
\def\begintab{%
  \global\let\v_box=\vtop
  \expandafter\global\_Mathe_false
  \global\Spalteneintrag_Format={\wider{##}}%
  \global\Spaltenzwischenraum_Format={&}%
  \global\Spaltenende_Format={\cr}%
  \_Mathe_normal_false
  \Beginn_der_Tabelle{}{}{}%
}%
\let\endtab=\Ende_der_Tabelle
\def\beginfixtab#1{%
 \global\let\v_box=\vtop
 \expandafter\global\_Mathe_false
 \global\Spalteneintrag_Format={\wider{##}}%
 \global\Spaltenzwischenraum_Format={&}%
 \global\Spaltenende_Format={\cr}%
 \_Mathe_normal_false
 \Beginn_der_Tabelle{}{}{to #1}%
}%
_der_Tabelle
\def\hline{%
 \noalign{%
  \unskip
  \hrule height\tabrulewidth
  \vskip\vfree
  \vskip-\vfree
 }%
}%
\def\span_omit{\fuege_Tokens_ans_Vielfachformat{\span\omit}}%
\def\und_omit{\fuege_Tokens_ans_Vielfachformat{&\omit\unskip}}%
\def\c_line[#1-#2]{%
 \noalign{\vskip-\tabrulewidth}%
 \omit
 \ne_Zahl=#1%
 \noch_ne_Zahl=#2%
 \global\Vielfach_Format={}%
 \advance \ne_Zahl by -1%
 \advance \noch_ne_Zahl by -\ne_Zahl%
 \loop
 \ifnum \ne_Zahl>1
  \advance\ne_Zahl by-1
  \span_omit
 \repeat
 \ifnum \ne_Zahl>0
  \und_omit
 \fi
 \fuege_Tokens_ans_Vielfachformat{%
  \leaders\hrule height\tabrulewidth\hfill}%
 \ifnum \noch_ne_Zahl>1
  \advance \noch_ne_Zahl by-1
  \span_omit
 \repeat
 \the\Vielfach_Format
 \crcr
}%
\def\multicolumn#1#2#3{%
 \multispan{#1}%
 {\_vor_erster_Spalte_true
  \global\Spaltenzwischenraum_Format={}%
  \global\Spalteneintrag_Format={#3}%
  \global\Spalten_Format={}%
  \if_Mathe_normal_
   \Uebersetze_ m#2\Uebersetze_Ende
  \else
   \Uebersetze_ h#2\Uebersetze_Ende
  \fi
  \the\Spalten_Format
 }%
 \ignorespaces
}%
\catcode`\_=8

$$
\vcenter{
\begintab{|c|c|c|c|c|}
\hline
Case&$q(j)$&2 in $\kay$ &value of $k$&${\Cal B}^{{\Cal O}^{\times}}$
\cr
\hline
1&$2\vert\alpha,\ \varepsilon\equiv 1(8)$&split&$k={\alpha\over
2}+1$&${\Cal A}$
\cr
\hline
2&$2\vert\alpha,\ \varepsilon\equiv 5(8)$&unramified&$k={\alpha\over
2}+1$&$\{ [\Lambda_0]\}$
\cr
\hline
3&$2\vert\alpha,\ \varepsilon\equiv -1(4)$&ramified&$k={\alpha\over
2}$&$\{ [\Lambda_0], \, [\Lambda_1]\}$
\cr
\hline
4&$2\nmid\alpha$&ramified&$k={\alpha-1\over 2}$&$\{[\Lambda_0],\,
[\Lambda_1]\}$
\cr
\hline
\endtab
}
$$
\medskip\noindent
We explain the last column in this table. In cases 1 and 2, writing
$j=2^{\alpha/2}\cdot
\bar j$, the index of $\Z_2[\bar j]$ in ${\Cal O}$ is 2. In case 1, the fixed point
set of ${\Cal O}^{\times}$ is the apartment $\Cal A$ in ${\Cal B}$ corresponding to
the split Cartan subgroup $\kay^{\times}$ of $GL_2(\Q_2)$. In case 2, the
fixed point set of ${\Cal O}^{\times}$ is the vertex corresponding to the
lattice  $\Lambda_0={\Cal O}$ in $\Q_2^2$. Note that in case 1 the fixed
point set of $j$ is
$${\Cal B}^j=\{[\Lambda];\ d([\Lambda],\ {\Cal A})\leq 1\}\ \
.\leqno(A.3)$$
In case 2, denoting by $[\Lambda_1]$ the vertex corresponding to the
lattice $\Lambda_1=\Z_2[\bar j]$, we have
$${\Cal B}^j=\{ [\Lambda_0],[\Lambda_1]\}\ \ .\leqno(A.4)$$
In cases 3 and 4 we write $j=2^{[\alpha/2]}\bar j$. Then we have ${\Cal
O}=\Z_2[\bar j]$. The fixed point set of ${\Cal O}^{\times}$ consists of
the vertices corresponding to the lattices $\Lambda_0={\Cal O}$ and
$\Lambda_1 =\pi {\Cal O}$, where $\pi$ denotes a uniformizer in ${\Cal
O}$. In case 3 this coincides with the fixed point set of $j$, whereas, in
case 4, $j$ permutes the two vertices $[\Lambda_0]$ and $[\Lambda_1]$ so
that ${\Cal B}^j$ consists of the midpoint of the edge formed by
$[\Lambda_0]$ and $[\Lambda_1]$.
\par
To formulate the theorem we write the divisor as usual as a sum of a vertical part
and a horizontal part,
$$Z(j)^{\roman{pure}}=Z(j)^{\roman{vert}}+Z(j)^{\roman{horiz}}\ \ .$$

\define\mult{\text{\rm mult}}

\proclaim{Proposition A.1} (i) Let
$$Z(j)^{\roman{vert}}=\sum_{[\Lambda]\in{\Cal B}}
\mult_{[\Lambda]}(j)\cdot{\Bbb P}_{[\Lambda]}\ \ .$$
Then the multiplicity $\mult_{[\Lambda]}(j)$ is given by
$$\mult_{[\Lambda]}(j)=\max (k-d([\Lambda], {\Cal B}^{{\Cal O}^{\times}}),
0)\ \ .$$
(ii) In case 1, $Z(j)^{\roman{horiz}}=0$. In case 2, $Z(j)^{\roman{horiz}}$
is isomorphic to the disjoint union of two copies of $\Spf\ W(\bar\F_2)$ and
meets the special fiber in two ordinary special points of ${\Bbb
P}_{[\Lambda_0]}$. In cases 3 and 4, $Z(j)^{\roman{horiz}}$ is isomorphic
to $\Spf\ W'$, where $W'$ is the ring of integers in a ramified quadratic
extension of $W(\bar\F_2)$, and meets the special fiber in the superspecial
point corresponding to the midpoint of the edge formed by $[\Lambda_0]$
and $[\Lambda_1]$.
\endproclaim

\demo{Proof}
We first determine $Z(j)\cap (\hat\Omega_{[\Lambda]}\times_{\Spf\
\Z_2}\Spf\ W(\bar\F_2))$ for a vertex $[\Lambda]$ where the intersection is
non-empty. Let $m=\max\{ r;\ j(\Lambda)\subset 2^r\Lambda\}$. Then
$$m=\alpha/2-d([\Lambda],\ {\Cal B}^j)\ \ ,\leqno(A.5)$$
cf.\ loc.cit., Lemma 2.8. After choosing a basis of $\Lambda$ we may write
$$j=2^m\cdot\left(\matrix\bar a&\bar b\cr \bar c&-\bar
a\cr\endmatrix\right)
=2^m\cdot\bar j\ \ ,\leqno(A.6)$$
where $\bar a, \bar b, \bar c$ are not simultaneously divisible by $p$.
The equation of $Z(j)$ on $\hat\Omega_{[\Lambda]}\times_{\Spf\ \Z_2}\Spf\
W(\bar\F_2)= \Spf\ W(\bar\F_2)[\,T, (T^2-T)^{-1}]\widehat{\phantom{R}}$ is given by
$$2^m\cdot (\bar b T^2-2\bar aT-\bar c)=0\ \ ,\leqno(A.7)$$
cf.\ loc.cit., (3.5). We now distinguish 2 cases.
\medskip\noindent
{\bf Case a:} $2\vert\bar b$ and $2\vert\bar c$. Then $2\notmid \bar a$
and we may write (A.7) in the form
$$2^{m+1}\cdot (\bar b_0T^2-\bar a T-\bar c_0)=0\ \ ,\leqno(A.8)$$
where $\bar b=2\bar b_0$ and $\bar c= 2\bar c_0$.
\par
Hence in this case the multiplicity $\mult_{[\Lambda]}(j)$ equals $m+1$. However, in this case
$\bar j\in GL_2(\Z_2)$ and hence $[\Lambda]$ is fixed by $j$. We check now
case by case when alternative a) can occur. Case 4 can be excluded right
away since in this case no vertex is fixed by $j$. In case 3 let
$[\Lambda]=[\Lambda_0]$ with the notation introduced in this case, i.e.\
$\Lambda_0={\Cal O}=\Z_2[\bar j]$.  Choosing as basis $1,\bar j$ we see
that $j$ is given by the matrix
$$j=2^{\alpha/2}\cdot\left(\matrix 0&\varepsilon\cr
1&0\cr\endmatrix\right) =2^m\,\bar j\ \ ,\leqno(A.9)$$
hence alternative a) does not occur for $[\Lambda_0]$. The case when
$[\Lambda]=[\Lambda_1]$ with $\Lambda_1=\pi{\Cal O}$ is identical and
hence alternative a) does not occur in case 3.
\par
In case 2, the vertex $[\Lambda_1]$ with $\Lambda_1=\Z_2[\bar j]$ is
excluded for the same reason. Now let us consider the vertex
$[\Lambda_0]$, with $\Lambda_0={\Cal O}$. We may choose the basis $1,
{1+\bar j\over 2}$ of $\Lambda_0$ and then $j$ is given by the matrix
$$j=2^{\alpha\over 2}\cdot\left(\matrix -1&2\lambda\cr
2&1\cr\endmatrix\right)\ \ ,\ \ \hbox{where}\ \varepsilon-1=4\lambda\ \
.\leqno(A.10)$$
Hence alternative a) applies here. Furthermore, in this case, the second factor in (A.8)
is equal to
$$\lambda T^2+T-1\ \ .\leqno(A.11)$$
Since, in case 2, we have $\lambda\in \Z_2^{\times}$, the ring
$$\Z_2[T]/(\lambda T^2+T-1)$$
is the ring of integers in an unramified quadratic extension of $\Q_2$ and
the zero's of the polynomial $T^2+T-1\in \F_2[T]$ lie in
$\F_4\setminus\F_2$ and define 2 ordinary special points of ${\Bbb
P}_{[\Lambda_0]}$.

In case 1 the analysis is similar to case 2. First one checks that if
$[\Lambda]\not\in {\Cal B}^{{\Cal O}^{\times}}$, then alternative a) does
not occur. If $[\Lambda]\in {\Cal B}^{{\Cal O}^{\times}}$, then after
replacing $[\Lambda]$ by $[g\Lambda]$ for some $g\in \kay^{\times}$ we may
assume that either $[\Lambda]=[\Lambda_0]$ with $\Lambda_0={\Cal
O}=\langle 1, {1+j\over 2}\rangle$ or $[\Lambda]=[\Lambda'_0]$ with
$\Lambda'_0=\langle 2, {1+\bar j\over 2}\rangle$. In the first case the
matrix of $j$ is given by (A.10) and hence we are in alternative a). The
second factor in (A.8) is equal to
$$\lambda T^2+T-1\equiv T-1\ {\roman{mod}}\ 2\ \ ,$$
since $2\vert\lambda$ in case 1. It follows that $Z(j)^{\roman{horiz}}\cap
(\hat\Omega_{[\Lambda_0]}\times_{\Spf\ \Z_2}\Spf\ W(\bar\F_2))=\emptyset$. In the
second case the matrix of $j$ is given by
$$j=2^{\alpha\over 2}\cdot\left( \matrix -1&\lambda\cr
4&1\cr\endmatrix\right)\ \ .\leqno(A.12)$$
Again we are in alternative a), since $2\vert\lambda$. The second factor
in (A.8) is equal to
$$\lambda_0T^2+T-2\equiv \lambda_0T^2+T\ {\roman{mod}}\ 2\ \ ,$$
where we have set $\lambda=2\lambda_0$. Since $T^2-T$ is invertible on
$\hat\Omega_{[\Lambda'_0]}$, we again have $Z(j)^{\roman{horiz}}\cap
(\hat\Omega_{[\Lambda'_0]}\times_{\Spf\ \Z_2}\Spf\ W(\bar\F_2))=\emptyset$. This
concludes our analysis of the alternative a).
\medskip\noindent
{\bf Case b:} $2\notmid\bar b$ or $2\notmid \bar c$. In this case the
second factor of (A.7) is not divisible by 2, and hence $\mult_{[\Lambda]} (j)=m$. At this point we have shown that
for any $[\Lambda]\in{\Cal B}$ the multiplicity $\mult_{[\Lambda]}(j)$ is given by the formula in (i).
Indeed, this follows by listing case by case the fixed point sets of $j$,
and the expressions for $k$ and for $m$, cf.\ (A.5), and comparing them
with the multiplicities calculated above.
\par
Now let us analyze the second factor in (A.7) in the alternative b). Its
image in $\F_2[T]$ is
$$\bar b\,T^2-\bar c\ \ ,$$
hence is equal to either $T^2-1$, 1, or $T^2$. In all cases
$Z(j)^{\roman{horiz}}\cap(\hat\Omega_{[\Lambda]}\times_{\Spf\ \Z_2} \Spf\
W(\bar\F_2))=\emptyset$.
\par
Now we determine $Z(j)^{\roman{horiz}}\cap
(\hat\Omega_{\Delta}\times_{\Spf\ \Z_2}\Spf\ W(\bar\F_2))$ for an edge $\Delta
=\{ [\Lambda], [\Lambda']\}$, where the intersection is non-empty.
As in the proof of Prop.\ 3.3 in \cite{\krinvent} we see that this
intersection is non-empty only when $d([\Lambda], {\Cal B}^j )=
d([\Lambda'], {\Cal B}^j)$. In cases 3 and 4, we therefore must have
$$[\Lambda]=[\Lambda_0]\ \ ,\ \ [\Lambda']=[\Lambda_1]\ \ .\leqno(A.13)$$
In case 3, we take as basis in standard form for $\Lambda_0, \Lambda_1$,
and, noting that $1+\bar j$ is a uniformizer of $\Cal O$,
$$\Lambda_0=\langle 1,1+\bar j\rangle\ \ ,\ \ \Lambda_1=\langle 2,1+\bar
j\rangle\ \ .\leqno(A.14)$$
In terms of the basis $1,1+\bar j$ of $\Lambda_0$ the matrix of $j$ is
$$j=p^{\alpha\over 2}\cdot\left( \matrix -1&-2\cdot(1-2\lambda)\cr
1&1\cr\endmatrix\right)\ \ ,\ \hbox{where}\ \varepsilon =4\lambda-1\ \
.\leqno(A.15)$$
By loc.\ cit.\ the equations for $Z(j)\cap
(\hat\Omega_{\Delta}\times_{\Spf\ \Z_2}\Spf\ W(\bar\F_2))$ in
$$\hat\Omega_{\Delta}\times_{\Spf\ \Z_2}\Spf\ W(\bar\F_2)= \Spf\ W(\bar\F_2) [\,T_0,
T_1, (1-T_0)^{-1}, (1-T_1)^{-1}]\widehat{\phantom{R}}/(T_0T_1-2)$$
are given by
$$\align
p^{\alpha/2}\cdot T_0\big(-(1-2\lambda)T_0+2-T_1\big)
&
=0\\
p^{\alpha/2}\cdot T_1\big(-(1-2\lambda)T_0+2-T_1\big)
&
=0\ .
\endalign$$
Hence, $Z(j)^{\roman{horiz}}$ is defined by the second factor in these
equations. Now putting \hfill\break
$\mu=(-(1-2\lambda))^{-1}\in\Z_2^{\times}$, we obtain
$$Z(j)^{\roman{horiz}}=\Spf\ W(\bar\F_2)[T_0]/(T_0^2+2\mu T_0-2\mu)\ \
.\leqno(A.16)$$
Since $T_0^2+2\mu T_0-2\mu$ is an Eisenstein polynomial, we see that
$Z(j)^{\roman{horiz}}$ is the formal spectrum of the ring of integers in a
ramified quadratic extension of $W(\bar\F_2)$ and it meets the special fiber
of $\hat\Omega_{\{ [\Lambda_0],[\Lambda_1]\} }\times_{\Spf\ \Z_2}\Spf\
W(\bar\F_2)$ in ${\roman{pt}}_{\Delta}$,  which finishes the proof in this case.

\par
The  case 4 is similar to the case of loc.\ cit., p.\ 180.\footnote{We note
that at this point in loc.\ cit.\ there is a slight error. The
equations (3.23) of loc.\ cit.\ do not define the same closed subscheme as
equations (3.22). The correct expression for $Z(j)^h$, replacing (3.24) is
$$\align
Z(j)^h
&
=\Spf\ W[T_0, T_1]\widehat{\phantom{j}} / (\bar b_0T_0-2\bar a-\bar c T_1, T_0T_1-p)\\
&
=\Spf\ W[T_0]/(T_0^2+\alpha T_0+\beta)\ ,
\endalign$$
where $\alpha\in p{\Z}_p$ and $\beta\in p\Z_p^{\times}$. The conclusions
drawn from these corrected equations (pp.\ 181, 182/183 loc.\ cit.) are
unchanged.}
In this case we have
$$j(\Lambda_0)=2^{\alpha-1\over 2}\cdot \Lambda_1\ \ ,\ \ j(\Lambda_1)=
2^{\alpha+1\over 2}\cdot\Lambda_0\ \ .\leqno(A.17)$$
Hence, as in loc.\ cit., we can write $j$ in terms of standard coordinates
for $\Lambda_0, \Lambda_1$
$$j= 2^{\alpha-1\over 2}\cdot\left( \matrix \bar a&\bar b\\ \bar c&-\bar
a\endmatrix\right)\ \ \hbox{with}\ \ \bar b= 2\cdot \bar b_0$$
and where $2\mid\bar a$ and $\bar b_0$ and $\bar c$ are units. Hence
$Z(j)^{\roman{horiz}}$ is isomorphic to
$$\Spf\ W(\bar\F_2)[\,T_0, T_1, (1-T_0)^{-1}, (1-T_1)^{-1}]\widehat{\phantom{R}} /(T_0T_1-2, \bar
b_0 T_0-2\bar a-\bar c T_1)\ \ .$$
Putting $\mu=\left( { b_0\over \bar c}\right)^{-1}\in
\Z_2^{\times}\ \ \hbox{and}\ \ \nu= {\bar a\over \bar c}$, we see that
$$Z(j)^{\roman{horiz}}=\Spf\ W(\bar\F_2) [T_0]/ (T_0^2-2\mu\nu
T_0-2\mu)\leqno(A.18)$$
which yields the assertion as in case 3.
\par
In case 2, we must have $\Delta =\{ [\Lambda_0], [\Lambda_1]\}$, where
$\Lambda_0={\Cal O}$ and $\Lambda_1=\Z_2[\bar j]$. In this case we take as
standard bases
$$\Lambda_0= \left\langle  {1+j\over 2}, 1\right\rangle\ \ ,\ \ \Lambda_1=
\left\langle 2\cdot {1+j\over 2}, 1\right\rangle\ \ .$$
But then, by (A.10), the equation for $Z(j)^{\roman{horiz}}$ is given by
$$T_0-1-\lambda T_1=0\ \ ,\ \ \hbox{where}\ \varepsilon -1=4\lambda\ \ .$$
It follows that ${\roman{pt}}_{\Delta}\not\in Z(j)^{\roman{horiz}}$, which
finishes the proof in this case.
\par
Finally there is case 1. In this case either $\Delta\subset {\Cal A}$ or,
after replacing $\Delta$ by $g\Delta$ where $g\in \kay^{\times}$, we may
assume that $\Delta= \{ [\Lambda_0],[\Lambda_1]\}$ where $\Lambda_0$ and
$\Lambda_1$ are as in case 3. The second alternative is treated as the
case 3 above. If $\Delta\subset {\Cal A}$, we may assume that
$$\Lambda = \left\langle 1,{1+\bar j\over 2}\right\rangle\ \ ,\ \
\Lambda'= \left\langle 2, {1+\bar j\over 2}\right\rangle\ \ .$$
In this case $j$ is given by the matrix (A.10) which yields the following
equations for $Z(j)$,
$$\align
p^{\alpha/2}\cdot T_0\big(\lambda T_0+2-2T_1\big)
&
=0\ ,\\
p^{\alpha/2}\cdot T_1\big(\lambda T_0+2-2T_1\big)
&
=0\ .
\endalign$$
Since we are in case 1, we may write $\lambda=2\lambda_0$ in the defining
equation $\varepsilon -1=4\lambda$. It follows that, after pulling a 2 out of the last factor
in both equations,
$Z(j)^{\roman{horiz}}$ is defined by the equation
$$\lambda_0T_0+1-T_1=0\ \ ,$$
and again ${\roman{pt}}_{\Delta}\not\in Z(j)^{\roman{horiz}}$.
\qed
\enddemo

\subheading{\Sec12. Archimedean contributions}

In this section, we compute the additional contribution
$$\kappa(m,v)={1\over 2}\int_{[\Gamma\setminus D]}\Xi(m,v)\,c_1(\hat\omega)\tag12.1$$
to the height pairing coming
from the fact that we are using nonstandard Green functions
defined in \cite{\annals} for the cycles $\Cal Z(m)$.
Recall that, by (5.8), for $z\in D$, we have
$$\Xi(m,v)(z) = \sum_{\matrix \scr x\in L(m)\endmatrix} \xi(v^{\frac12}x,z),\tag12.2$$
where
$$\xi(x,z) = - \Ei(-2\pi R(x,z)).$$
Note that the quantity
$$R(x,z) = - (\pr_z(x),\pr_z(x)),$$
and hence $\xi(x,z)$, is independent of the orientation of the plane $z$. Also
recall that $c_1(\hat\o)=\mu$.

\proclaim{Proposition 12.1} (i) If $m>0$, then
$$\kappa(m,v) = 2\,\delta(d,D)\, H_0(m;D)\cdot \frac12\,J(4\pi mv),$$
where
$$J(t) = \int_0^\infty e^{-t w } \,\big[\,(w+1)^{\frac12} -1\big]\,w^{-1}\,dw,$$
is as in Theorem~8.8,
and $H_0(m;D)$ is given by (8.19).
\hfill\break
(ii) If $m<0$, then
$$\kappa(m,v) = 2\,\delta(d;D)\,H_0(m;D)\,\frac1{4\pi}\,|m|^{-\frac12}\,v^{-\frac12}\,\int_1^{\infty}e^{-4\pi|m|vw}\,w^{-\frac32}\,dw,
$$
where $H_0(m;D)$ is given by (8.32).
\endproclaim

\demo{Proof}
We have
$$\align
\kappa(m,v) &= \frac14\int_{\Gamma\back D} \Xi(m,v)\cdot\mu\\
\nass
{}&= \frac14\int_{\Gamma\back D} \sum_{x\in L(m)} \xi(v^{\frac12}x,z)\,d\mu(z)\\
\nass
{}&= \frac14\sum_{\matrix \scr x\in L(m)\\ \scr \mod \Gamma\endmatrix}
\int_{\Gamma_x\back D} \xi(v^{\frac12}x,z)\,d\mu(z).
\endalign
$$
First suppose that $m= Q(x)>0$, so that $\Gamma_x$ is finite.
Then
$$
\int_{\Gamma_x\back D} \xi(x,z)\,d\mu(z) = 2\,|\Gamma_x|^{-1}\cdot \int_{D}\xi(x,z)\,d\mu(z)
=4\,|\Gamma_x|^{-1}\cdot \int_{D^+}\xi(x,z)\,d\mu(z).\tag12.3
$$
Here the factor of $2$ occurs since $\Gamma_x$ contains $\pm1$, but these elements act trivially
on $D$, while in the second step, we use the fact that $\xi(x,z)$ does not
depend on the orientation of $z$.
Since
$$\xi(gx,gz) = \xi(x,z),\tag12.4$$
for $g\in GL_2(\R)$,
we may assume that
$$x = m^{\frac12}\cdot x_0 = m^{\frac12}\cdot \pmatrix {}&1\\-1&{}\endpmatrix.\tag12.5$$
Then, writing $z = k_{\theta}(e^{t}i)\in \H\simeq D^+$, \cite{\annals}, p.601 , we have
$$R(x,z) = 2m\,\sinh^2(t).\tag12.6$$
Now, noting that $t$ runs from $0$ to $\infty$ and $\theta$ runs from $0$ to $\pi$,
$$\align
I:&=\int_{D^+} -\Ei(-2\pi R(v^{\frac12} x,z))\,d\mu(z)\\
\nass
{} &= \frac1{2\pi}\,\int_0^{\pi} \int_0^\infty -\Ei(-4\pi m v \sinh^2(t))\
2\sinh(t)\,dt\,d\theta\\
\nass
{}&=\frac12\,\int_0^\infty\left(\int_1^\infty e^{-4\pi m v \sinh^2(t) r} r^{-1}\,dr\right)\,2\sinh(t)\,dt\\
\nass
{}&= \frac12\,\int_0^\infty\left(\int_1^\infty e^{-4\pi m v w r} r^{-1}\,dr\right)\,(w+1)^{-\frac12}\,dw
\endalign
$$
But
$$\align
&\int_0^\infty e^{-4\pi m v w r} \,(w+1)^{-\frac12}\,dw\\
\nass
{} &= e^{-4\pi m v w r} \,2(w+1)^{\frac12}\bigg|_0^\infty
+4\pi m v r \int_0^\infty e^{-4\pi m v w r} \,2(w+1)^{\frac12}\,dw\\
\nass
{}&= -2 + 4\pi m v r \int_0^\infty e^{-4\pi m v w r} \,2(w+1)^{\frac12}\,dw\\
\nass
{}&= 8\pi m v r \int_0^\infty e^{-4\pi m v w r} \,\big[\,(w+1)^{\frac12} -1\big]\,dw
\endalign
$$
so that
$$\align
I&=\frac12\, \int_1^\infty 8\pi m v  \int_0^\infty e^{-4\pi m v w r} \,\big[\,(w+1)^{\frac12} -1\big]\,dw\,dr\\
\nass
{}&= \int_0^\infty e^{-4\pi m v w } \,\big[\,(w+1)^{\frac12} -1\big]\,w^{-1}\,dw\\
\nass
{}&=\,J(4\pi mv).
\endalign
$$
By Lemma~9.2, we have
$$\sum_{\matrix \scr x\in L(m)\\ \scr \mod \Gamma\endmatrix}
|\Gamma_x|^{-1} = \delta(d;D)\cdot H_0(m;D).\tag12.7$$
Collecting terms, we obtain (i).

Next suppose that $m<0$. Let $\Gamma^+ = \Gamma\cap \text{\rm GL}_2(\R)^+$, and let
$\delta_x = |\Gamma_x:\Gamma_x^+|$, where $\Gamma_x^+=\Gamma_x\cap \Gamma^+$.
Then
$$\int_{\Gamma_x\back D} \xi(x,z)\,d\mu(z) = \delta_x^{-1}\int_{\Gamma_x^+\back D} \xi(x,z)\,d\mu(z)
= 2 \delta_x^{-1}\int_{\Gamma_x^+\back D^+} \xi(x,z)\,d\mu(z).\tag12.8$$
By conjugating by a suitable
$g\in GL_2(\R)$, we can take
$$ g\cdot x = |m|^{\frac12}\cdot x_0 = |m|^{\frac12}\cdot \pmatrix 1&{}\\{}&-1\endpmatrix.\tag12.9$$
Let $\Gamma'$ be the corresponding conjugate of $\Gamma^+$ in $SL_2(\R)$, and
note that $\Gamma'_{g x}$ will then be generated by $\pm 1_2$ and a unique element
$$\pmatrix \e(x)&{}\\{}&\e(x)^{-1}\endpmatrix\tag12.10$$
for $\e(x)>1$ the fundamental unit {\it of norm 1} in the order $i_x^{-1}(O_B)$. If we write
$z = r e^{i\theta}\in \H\simeq D^+$, then $\Gamma'_{gx}$ acts
by multiplication by powers of $\e(x)^2$. Note that
$$R(gx,z) = \frac{2|m|}{\sin^2(\theta)}.\tag12.11$$
Then
$$
\align
\int_{\Gamma^+_x\back D^+} \xi(v^{\frac12}x,z)\,d\mu(z) & = \frac1{2\pi}\int_1^{\e(x)^2}\int_0^\pi
-\Ei\left(-\frac{4\pi |m|v}{\sin(\theta)^{2}} \right)\,
\,r^{-1}(\sin(\theta))^{-2}\,d\theta\,dr\\
\nass
\nass
{}&=\frac2{\pi}\log|\e(x)|\cdot \int_0^{\pi/2}
-\Ei\left(-\frac{4\pi |m|v}{\sin(\theta)^{2}} \right)\,
\,(\sin(\theta))^{-2}\,d\theta\\
\nass
\nass
{}&=\frac1{\pi}\log|\e(x)|\cdot \int_1^{\infty}\bigg(\int_1^\infty
e^{-4\pi |m|vtw} \,w^{-1}\,dw\,\bigg)
\,(t-1)^{-\frac12}\,dt\\
\nass
\nass
{}&=\frac1{\pi}\log|\e(x)|\cdot \int_1^{\infty}e^{-4\pi|m|vw}\bigg(\int_0^\infty
e^{-4\pi |m|vtw} \,t^{-\frac12}\,dt\bigg)
\,w^{-1}\,dw\\
\nass
\nass
{}&=\frac1{\pi}\log|\e(x)|\cdot \Gamma(\frac12)\,(4\pi|m|v)^{-\frac12}\, \int_1^{\infty}e^{-4\pi|m|vw}\,w^{-\frac32}\,dw\\
\nass
\nass
{}&=\frac1{2\pi}\log|\e(x)|\cdot (|m|v)^{-\frac12}\,\int_1^{\infty}e^{-4\pi|m|vw}\,w^{-\frac32}\,dw.
\endalign
$$

The analogue of Lemma~9.2 for $m<0$ is the following.
\proclaim{Lemma 12.3} If $m<0$, then
$$\bigg( \sum_{\matrix \scr x\in L(m)\\ \scr \mod \Gamma\endmatrix}
2\delta_x^{-1}\log|\e(x)|\bigg)
= 4\,\delta(d;D)\, H_0(m;D),$$
where $H_0(m;D)$ is as in (8.32).
\endproclaim

\demo{Proof} We proceed as in the proof of Lemma~9.2. Note that for $x$
of type $c$, $\Gamma_x\simeq O_{c^2d}^\times$ and $\Gamma_x^+\simeq O_{c^2d}^1$, the subgroup of norm $1$ elements,
so that
$\delta_x = \delta_c = |O_{c^2d}^\times:O_{c^2d}^1|$. Let $\e(c^2d)$  be a fundamental unit in
$O_{c^2d}$ with $\e(c^2d)>1$. Note that
$\e^+(c^2d) = \e(c^2d)^{\delta_c}$ is a generator of $O_{c^2d}^1/{\pm1}$. Then  we have
$$\align
&\bigg( \sum_{\matrix \scr x\in L(m)\\ \scr \mod \Gamma\endmatrix}
2\delta_x^{-1}\log|\e(x)|\bigg)\\
\nass
{} &=\sum_{\scr c|n \atop \scr (c,D) =1} 2 \delta_c^{-1}|\Optoc|\cdot \log|\e^+(c^2d)|\\
\nass
{}&=2\sum_{\scr c|n \atop \scr (c,D) =1} |\Optoc|\cdot \log|\e(c^2d)|\\
\nass
{}&=2\,\delta(d;D)\sum_{\scr c|n \atop \scr (c,D) =1} h(d)\cdot \frac{\log|\e(d)|}{\log|\e(c^2d)|}\cdot
\bigg(c\prod_{\ell|c}(1-\chi_d(\ell)\ell^{-1})\bigg) \log|\e(c^2d)|\\
\nass
{}&=4\,\delta(d;D)\, \frac{h(d)\,\log|\e(d)|}{w(d)}\sum_{\scr c|n \atop \scr (c,D) =1}
c\prod_{\ell|c}(1-\chi_d(\ell)\ell^{-1})\\
\nass
{}&=4\,\delta(d;D)\,H_0(m;D).
\endalign
$$
\qed
\enddemo
\enddemo

\demo{Proof of Theorem~7.2 (concluded)} We observe that the expression in Corollary~12.2 when $m>0$ coincides  with
the term (8.44) in the Fourier coefficient $\Cal E'_m(\tau,\frac12;D) = A'_m(\frac12,v)\,q^m$.
On the other hand, when $m<0$, the expression in Corollary~12.2 coincides with that in (iii) of Theorem~8.8.
\qed\enddemo

\vfill\eject

\subheading{\Sec13. Remarks about the constant term}

In this section we explain the motivation for our definition of $\hat\Cal Z(0,v)$.
The key point is the comparison of our expression for the constant term
of $\Cal E'(\tau,\frac12;D)$ with the result of Bost \cite{\bostumd} and K\"uhn \cite{\kuehn} concerning
$\lan \hat\o,\hat\o\ran$ in the case of the modular curve.

From Theorem~8.8, we have $$ \Cal
E'_0(\tau,\frac12;D)=c(D)\,\L_D(2)\,\bigg[\frac12\,\log(v) -
2\frac{\zeta'(-1)}{\zeta(-1)} -1 + 2C+ \sum_{p\mid
D}\frac{p\log(p)}{p-1}\,\bigg],\tag13.1 $$ where $C$ is as in
Definition~3.4.

Now suppose, for a moment, that $D=1$. Then, recalling (6.28), we have
$$c(D)\zeta_D(2) = -\frac1{12} = \zeta(-1) = -\vol(\M(\C)) = -\,\deg(\hat\o)$$
and so we obtain:
$$\Cal E'_0(\tau,\frac12;D)\big\vert_{D=1} = -2\big[\,\zeta'(-1) +\frac12\,\zeta(-1)\, \big] -\deg(\hat\o)\bigg[\,\frac12\,\log(v)
+2C \, \bigg].\tag13.2$$

Let $\hat\o_o$ denote the Hodge bundle with the metric defined by (3.11) or (10.15), so that,
as elements of $\CH^1(\M)$,
$$\hat\o = \hat\o_o + (0,2\,C).\tag13.3$$
Then, $\hat\o_o$ is the bundle of modular forms of weight $2$
with its Petersson metric and hence, by K\"uhn  and Bost,
$$\lan \hat\o_o,\hat\o_o\ran^\natural = 4 \,\big[\,\zeta'(-1) +\frac12\,\zeta(-1)\, \big].\tag13.4$$
Here $\lan\ ,\ \ran^\nat$ denotes the height pairing without the `stack' aspect!
Thus, for the renormalized metric, we have
$$\align
\lan \hat\o,\hat\o\ran^\natural &= \lan \hat\o_o,\hat\o_o\ran^\natural + 2\,C\,\deg^\nat(\o)\tag13.5\\
\nass
 {}&= 4 \,\big[\,\zeta'(-1) +\frac12\,\zeta(-1)\, \big] + \deg^\nat(\o)\,2C
\endalign
$$
Note that
$$\lan\ ,\ \ran = \frac12\,\lan\ ,\ \ran^\nat.\tag13.6$$

Then, recalling that
$$\hat\Cal Z(0,v) = -\,\bigg(\,\hat\o + (0,\log(v))\ \bigg),$$
we have
$$\align
&\lan \hat\Cal Z(0,v),\hat\o\ran\\
\nass
{} &= - \lan\hat\o,\hat\o\ran -\frac12\deg(\o)\,\log(v)\tag13.7\\
\nass
{}&= - 2\big[\,\zeta'(-1) +\frac12\,\zeta(-1)\, \big] - \deg(\o)\,\bigg[\frac12\log(v)
+2C\,\bigg]
\endalign
$$
This agrees perfectly with our constant term $\Cal E'_0(\tau,\frac12;D)\vert_{D=1}$.

Finally, for general $D=D(B)>1$, this discussion suggests that
$$\align
\lan\hat\o_o,\hat\o_o\ran &= -\lan \hat\Cal Z(0,v),\hat\o\ran -\frac12\deg(\o)\,\log(v)-\deg(\o)\,2C\\
\nass
{}&\overset{??}\to{=} -\Cal E'_0(\tau,\frac12;D)  -\frac12\deg(\o)\,\log(v)-\deg(\o)\,2C\tag13.8\\
\nass
{}&= -c(D)\,\L_D(2)\,\bigg[ - 2\frac{\zeta'(-1)}{\zeta(-1)} -1+ \sum_{p\mid D}\frac{p\log(p)}{p-1}\,\bigg]\\
\nass
{}&= \zeta_D(-1)\,\bigg[ 2\frac{\zeta'(-1)}{\zeta(-1)} +1 - \sum_{p\mid D}\frac{p\log(p)}{p-1}\,\bigg],
\endalign
$$
since
$c(D)\,\L_D(2) = \zeta_D(-1) = \zeta(-1) \,\prod_{p\mid D} (p-1)$.

\subheading{Part IV. Computations: analytic}

\subheading{\Sec14. Local Whittaker functions: the non-archimedean  case}

The main purpose of this section is to prove Proposition 8.1. We fix a prime $p$
and frequently drop the subscript $p$ to lighten the notation.
Recall that $B=B_p$
is a quaternion algebra which is a matrix algebra or a division algebra
depending on whether $p\nmid D$ or $p\mid D$. Here $D$ is a fixed square-free
positive integer. Let $O_B$ be a maximal order of $B$ and let
$$
V=\{\ x \in B\mid\, \tr^0 x =0 \}
$$
with the quadratic form $Q(x)= \kappa\, x^2$, where
$\kappa=\pm1$. Actually, only the case $\kappa=-1$ is needed in section 8, but we
treat the slightly more general case for future reference. Let $L=V \cap O_B$,
and let
$S\in \Sym_2(\Q_p)$ be the matrix associated to $L$ in the following sense. With
respect to a basis of $L$ over $\Bbb Z_p$, identify $L$ with
$\Bbb Z_p^3$. Then, for any $x \in L= \Bbb Z_p^3$,
$$
Q(x) =\frac{1}2(x, x) =  {}^tx S x.
\tag{14.1}
$$
Let $S_r= S \perp \frac{1}2\pmatrix 0 & I_r \\ I_r & 0 \endpmatrix$, and
let $L_r=\Bbb Z_p^{2r+3}$ be the associated quadratic lattice, viewed as
the direct sum of $L$ and $r$ hyperbolic planes. Let $dx=\prod dx_i$ be the
standard Haar measure on $L_r$, where
$$
\int_{\Bbb Z_p} dx_i =1.
$$
Let
$$
W(m, S_r) = \int_{\Bbb Q_p} \int_{L_r} \psi(b{}^tx S_r x )\, \psi(- mb)\, dx\, db
\tag{14.2}
$$
be the integral defined in \cite{\yangden}, (1.2). It is the same as the local quadratic
density polynomial $\alpha_p(X, m, S)$ with $X=p^{-r}$ defined in \cite{\yangden}, page
312. Here $\psi=\psi_p$ is the local component of our standard additive character of $\A/\Q$.

\proclaim{Lemma 14.1} With the notation as above, for any $r \ge 0$ and
any $m \in \Bbb Q_p$,
$$
W_{m, p} (1, r+\frac{1}2, \Phi_p)
 =\beta(V)\, |\det 2 S\,|_p^{\frac{1}2}\,  W(m, S_r).
$$
Here
 $$
\beta(V)=\bigg(\, \epsilon(V)\, \gamma(\frac{1}2\psi_p)^3 \,
         \gamma(\det V, \frac{1}2\psi)\,\bigg)^{-1}
$$
is the local splitting index defined in \cite{\splitting}, Theorem 3.1.
Here $\e(V)$ is the Hasse invariant of $V$ and $\gamma(\psi)$ and $\gamma(a,\psi)$ are the
local Weil indices as in \cite{\rao}.
\endproclaim
\demo{Proof} We remark that this proposition is true in general for any
3-dimensional quadratic space $V$ over $\Bbb Q_p$. Let $\phi=\cha(L)$ and
$\phi_r= \cha(L_r)$. Then \cite{\annals}, Appendix, asserts that
$$
\omega(g) \phi_r(0) = \Phi_p(g, r+\frac{1}2)
$$
where $\omega=\omega_{\psi}$ is the Weil representation of $G_p'$
(the metaplectic cover of $SL_2(\Bbb Q_p)$) on the space $S(V)$ of
Schwartz functions on $V$. Thus
$$
\align
&W_{m, p} (1, r+\frac{1}2, \Phi_p)
\\
&=\int_{\Bbb Q_p}\omega(w n(b)) \phi_r(0) \, \psi(-bm)\,db\tag14.3
\\
&=\beta(V) \int_{\Bbb Q_p} \int_{V} \phi_r(x)\, d_rx \, \psi(-bm)\,db,
\endalign
$$
where $d_rx$ is the self-dual Haar measure with respect to
the bi-character $(x, y) \mapsto \psi_p( (x, y)_r)$ with $(x, y)_r$ the
bilinear from on $V_r =L_r \otimes \Bbb Q_p$ associated to $S_r$. It is
easy to check
$$
d_r x = |\det 2 S_r|_p^{\frac{1}2} dx =|\det 2 S\,|_p^{\frac{1}2} dx
\tag{14.4}
$$
under the identification $L_r= \Bbb Z_p^{2r+3}$ as above. This proves the
proposition.
\qed\enddemo

The following lemma is standard.

\proclaim{Lemma 14.2} Let $V$ and $L$ be as above. \hfill\break
(i)  When $p \nmid D$, one has  $L \cong \Bbb Z_p^3$ with the quadratic
form
and symmetric matrix as follows:
$$
Q(x)=\kappa\,(x_1^2 +x_2 x_3), \qquad
S=\kappa \pmatrix  1 & &\\ & &\frac{1}2 \\ {} &\frac{1}2 &\endpmatrix.
$$
(ii)  When $2 \ne p \mid D$ , one has  $L \cong \Bbb Z_p^3$ with the quadratic
form
and symmetric matrix as follows:
$$
Q(x)=\kappa\,(\beta x_1^2 +p x_2^2 - \beta p x_3^2),
\qquad
S=\kappa\,\diag(\beta, p, -\beta p).
$$
Here $\beta\in \Bbb Z_p^\times$ with $(\beta, p)_p=-1$.\hfill\break
(iii) When $p=2 \mid D(B)$, one has $L \cong \Bbb Z_p^3$ with the quadratic form
and symmetric matrix as follows:
$$
q=\kappa\,(-3x_1^2+2x_2^2 + 2x_2 x_3 + 2x_3^2),
\qquad
S =\kappa\, \pmatrix -3 & & \\ & &2 &1\\ & &1 &2 \endpmatrix.
$$
(iv)  Finally,
$$
|\det 2 S\,|_p=|2|_p\cdot\cases
1 &\text{if $p \nmid D$,}\\
\nass
p^{-2} &\text{ if $p \mid D.$}
 \endcases
$$
\endproclaim

Notice that $\beta(V)$ in Lemma 14.1  depends only on $V$ and is
well-defined even when $p=\infty$.

\proclaim{Lemma 14.3}  Let $V$  be as above.
Let $\zeta_8 =e^{\frac{2 \pi i}8}$.  Then
$$
\beta(V)=\epsilon(B) \cdot
\cases
1 &\text{ if $p \nmid 2 \infty$,}\\
\nass
\zeta_8^{-\kappa} &\text{ if $ p =\infty$,}\\
\nass
\zeta_8^{\kappa} &\text{ if $ p =2$.}
 \endcases
$$
Here $\epsilon(B)=\pm 1$ depending on whether $B=B_p$ is split or not.
\endproclaim
\demo{Proof} It is easy to see from Lemma 14.2 that
$\det V \in -\kappa\ \Q_p^{\times, 2}$ in all cases, and so, by \cite{\rao},
 $$
\gamma(\det V, \frac{1}2\psi)
 =\gamma(-\kappa, \frac{1}2\psi)
 =\gamma(\frac{1}2\psi)^{-2a(\kappa)}.\tag14.5
$$
Here $a(\kappa) =(1+\kappa)/2$.
By \cite{\rao}, A.10, A.11, one has
$$
\gamma(\frac{1}2\psi_p)
 = \cases
1 &\text{ if  $p \nmid 2 \infty$,}\\
\nass
\zeta_8  &\text{ if $ p = \infty$.}
\endcases\tag14.6
$$
When $p=2$, following the principle in \cite{\rao}, page 370, one has
$$
\gamma(\frac{1}2\psi)
=\frac{1}2 \sum_{x \in \Bbb Z/4} \psi(\frac{1}8 x^2)
=\zeta_8^{-1}.\tag14.7
$$
As for the Hasse invariant $\epsilon(V)$, one has in the split
case
$$
\epsilon(V)
 = (\kappa, -1)_p
 =\gamma(\frac{1}2\psi_p)^{4(1-a(\kappa))},\tag14.8
$$
where $(\ ,\ )_p$ is the quadratic HIlbert symbol for $\Q_p$.
In the ramified case,  one has
$$
\epsilon(V)
 =(\eta \kappa, -\kappa^2 \eta p^2)_p ( \kappa p, - \kappa \eta p)_p
 =(\kappa, -1)_p (p, \eta)_p
 =-\gamma(\frac{1}2\psi_p)^{4(1-a(\kappa))}.\tag14.9
$$
Now the lemma follows from the formula
$$
\beta(V)=\{ \epsilon(V) \gamma(\frac{1}2\psi_p)^3
         \gamma(\det V, \frac{1}2\psi)\}^{-1}.\qed\tag14.10
$$
\enddemo

{\bf{Proof of Proposition 8.1}} By Lemmas 14.1 and 14.3, Proposition 8.1 is
equivalent to the following  proposition for $\kappa =-1$.

\proclaim{Proposition 14.4} For a nonzero integer $m$,  write $4 m=n^2 d$ such
that
$ \kappa d$ is a fundamental discriminant of a quadratic field. \hfill\break
(i) If $p \nmid D$,
$$
W_p(m, S_r) =\frac{L_p(r+1, \chi_{-\kappa d})\, b_p(n, r+1; D)}{\zeta_p(2r+2)}.
$$
(ii) If $p\mid D$,
$$
W_p(m, S_r) = L_p(r+1, \chi_{-\kappa d})\,b_p(n, r+1; D) .
$$
(iii) If $m=0$,
$$
W_p(0, S_r)
=\cases
\frac{\zeta_p(2r+1)}{\zeta_p(2r+2)} &\text{if $ p \nmid D$,}\\
\nass
\frac{\zeta_p(2r+1)}{\zeta_p(2r)} &\text{if $p \mid D$,}
\endcases
$$
\endproclaim
\demo{Proof}  Part (i)  is a  better reformulation of \cite{\yangden}, Propositions 8.3.
Part (ii)
follows from (i) and \cite{\yangden}, Proposition 8.2.  Part (iii) follows from (i) and (ii) when we let
$a=\ord_p m$ tends to infinity.  We verify the case $p=2 \nmid D$ and
leave  the
other (easier) cases  to the reader.
We recall again that $W_p(m, S_r)$ is just the local density polynomial
$\alpha_p(X, m, S)$ in \cite{\yangden} with $X=p^{-r}$. Write $m=\alpha p^a$ with
$a = \ord_p m=2 k + \ord_p\frac{d}4 $ and $\alpha \in \Z_p^\times$ where
$k=k_p(n) =\ord_pn$ as before. In the notation of \cite{\yangden}, Proposition 8.3,  one
has $(\frac{\alpha \kappa}p)=\chi_{-\kappa d}(p)$ if $p \nmid d$.
In our case $p=2 \nmid D$,
$a=\ord_2 m =2k-2 +\ord_2 d$ with $k=\ord_2 n$ as before.

{\bf{Subcase 1}.} First we assume $8 \mid d$. Then  $a= 2k+1$
is odd, and \cite{\yangden}, Proposition 8.3(1), implies
$$
\align
W_2(m, S_r)
 &= (1- 2^{-2} X^2) \sum_{l=0}^{k} (2^{-1} X^2)^l
\\
 &= \frac{ (1- 2^{-2} X^2)(1- (2^{-1}X^2)^{k+1})}{1- 2^{-1}X^2}\tag14.11
\\
 &=\frac{L_2(r+1, \chi_{-\kappa d})\,b_2(n, r+1; D)}{\zeta_2(2s+2)}
\endalign
$$
as desired.

{\bf{Subcase 2}.} Now assume that $\ord_2 d= 2$. Then $a = 2k$,
$\alpha \kappa =\frac{\kappa d}{4} (n 2^{-k})^2 \equiv  -1 \mod 4$, and thus
$(\frac{-1}{\alpha \kappa})=-1$ and $\delta_8(\alpha -\kappa) =0$
So \cite{\yangden}, Proposition 8.3(3),
($\frac{a-1}2$ in the summation there should be $\frac{a}2$) implies
$$
\align
W_2(m, S_r)
 &=1+2^{-1} \sum_{l=1}^{k} (2^{-1}X^2)^l-2^{-k-2} X^{2k+2}
\\
 &= \frac{ (1- 2^{-2} X^2)(1- (2^{-1}X^2)^{k+1})}{1- 2^{-1}X^2}\tag14.12
\\
 &=\frac{L_2(r+1, \chi_{-\kappa d})}{\zeta_2(2r+2)}\,b_2(n, r+1; D).
\endalign
$$

{\bf{Subcase 3}.} Finally if $ 2 \nmid d$, i.e., $\kappa d  \equiv 1 \mod 4$.
Then
$ a = 2k-2$ and $ \alpha=d(n 2^{-k+1})^2 \equiv d \mod 8 $. In this case,
 $(\frac{-1}{\alpha \kappa})=1$ and
$$
\delta_8(\alpha -\kappa) =\delta_8(d -\kappa ) =\chi_{-\kappa d}(2).\tag14.13
$$
Set $v_2= \chi_{-\kappa d}(2)$. Then \cite{\yangden}, Proposition 8.3(3), gives
 $$
\align
W_2(m, S_r)
 &=1+2^{-1} \sum_{l=1}^{k-1} (2^{-1}X^2)^l
   + 2^{-k-1} X^{2k}
   +v_2 2^{-k-1}X^{2k+1}\tag14.14
\\
 &=\frac{1-2^{-2}X^2 - 2^{-k-1}X^{2k+1}(-v_2 +2^{-1} X  +v_2 2^{-1} X^2)}
        {1-2^{-1}X^2}.
\endalign
$$
On the other hand,
$$
\align
&\frac{L_2(r+1, \chi_{-\kappa d}) b_2(n, r+1; D)}{\zeta_2(2r+2)}
\\
 &= \frac{(1+v_2 2^{-1}X)(1-v_2 X+ v_2 2^{-k-2} X^{2k+1} - 2^{-k-1}X^{2k+2})}
         {1-2^{-1}X^2}\tag14.15
\\
 &=\frac{1- 2^{-2}X^2  - 2^{-k-1}X^{2k+1} (1+v_2 2^{-1}X)(-v_2 + X)}
        {1-2^{-1}X^2}
\\
 &=\frac{1-2^{-2}X^2 - 2^{-k-1}X^{2k+1}(-v_2 +2^{-1} X  +v_2 2^{-1} X^2)}
        {1-2^{-1}X^2}.
\endalign
$$
Therefore
$$
W_2(m, S_r) = \frac{L_2(r+1, \chi_{-\kappa d})}{\zeta_2(2r+2)}\,b_2(n, r+1; D).\qed\tag14.16
$$
\enddemo

\subheading{\Sec15. Local Whittaker functions: the archimedean  case}

In this section, we compute the local Whittaker function
$$
W_{m, \infty}(\tau, s, \Phi_{\infty}^{\ell})
 =v^{-\frac{1}2\ell}  \int_{\Bbb R} \Phi_{\infty}^{\ell}(w n(b) g_{\tau}', s)\,
  \psi_{\infty}( -mb) \, db
\tag{15.1}
$$
and prove Lemmas~8.9, and 8.11. Here $ \ell \in \frac{1}2\Bbb Z$ is such
that $\ell\equiv \frac32 \mod 2\Z$. In this paper, we only need $\ell=\frac{3}2$.

  Let
 $$
\Psi(a, b; z)
 =\frac{1}{\Gamma(a)}\int_0^{\infty}
   e^{-z r} (r+1)^{b-a-1} r^{a-1} dr
\tag{15.2}
$$
be the standard confluent hypergeometric function of the second kind, \cite{\lebedev}, where
$a>0$, $z >0$ and $b$ is any real number. It satisfies the
functional equation, \cite{\lebedev}, p.~265
$$
\Psi(a, b; z)
 =z^{1-b} \Psi(1+a-b, 2-b; z).
\tag{15.3}
$$
For convenience, we also define
$$
\Psi(0, b;z)=\lim_{a\rightarrow 0+}  \Psi(a, b; z)=1.
\tag{15.4}
$$
So $\Psi(a, b; z)$ is well-defined for $z>0$, $a \ge \min\{0, b-1\}$.
Finally, for any number $n$, we define
$$
\Psi_n(s, z)
 =\Psi(\frac{1}2(1+ n+s), s+1; z).
\tag{15.5}
$$
Then $(15.3)$ implies
$$
\Psi_n(s, z)=z^{-s} \Psi_n(-s, z).
\tag{15.6}
$$

\proclaim{Proposition 15.1} Let $q=e(m \tau)$, $(-i)^\ell = e(-\ell/4)$, and
$$
\alpha = \frac{s+1+\ell}2, \qquad  \beta=\frac{s+1-\ell}2.
$$
(i) For $m > 0$,
$$
W_{m, \infty}(\tau, s, \Phi_{\infty}^\ell)
 =2 \pi\,(-i)^\ell\, v^\beta\, (2 \pi m)^s\, \frac{\Psi_{-\ell}(s, 4 \pi m v)
}{\Gamma(\alpha)}\cdot q^m.
$$
(ii) For $m  <0$,
$$
W_{m, \infty}(\tau, s, \Phi_{\infty}^\ell)
 =2 \pi\, (-i)^\ell\, v^\beta\, (2 \pi| m|)^s \, \frac{\Psi_\ell(s, 4 \pi |m| v)
}{\Gamma(\beta)}\, e^{-4 \pi |m| v}\cdot q^m.
$$
(iii) For $m=0$,
$$
W_{0, \infty}(\tau, s, \Phi_{\infty}^\ell)
=2 \pi\, (-i)^\ell\, v^{\frac{1}2(1-\ell-s)}\, \frac{2^{-
s}\Gamma(s)}{\Gamma(\alpha)\Gamma(\beta)}.
$$
(iv) The special value at $s = \ell-1$ is
$$
W_{m, \infty}(\tau, \ell-1, \Phi_{\infty}^\ell)
=\cases
 0 &\text{if $ m \le 0$,}\\
\nass
\frac{(-2 \pi i)^\ell}{ \Gamma(\ell)}\,m^{\ell-1}\,q^m &\text{ if $ m >0$.}
\endcases
$$
\endproclaim
\demo{Proof} A standard calculation, \cite{\shimura}, (see also \cite{\annals} pages 585-586, for
this special case) gives
$$
W_{m, \infty}(\tau, s, \Phi_{\infty}^\ell)
=(-i)^\ell v^\beta e(m \bar \tau) \frac{(2 \pi)^{s+1}}{\Gamma(\alpha)\Gamma(\beta)}
  \int_{ r >0 \atop r>m} e^{- 4 \pi v r} (r-m)^{\beta-1} r^{\alpha -1}
\, dr
\tag{15.7}
$$
When $m=0$, this gives (iii) immediately.

When $ m >0$, the integral equals
$$
\align
\int \Sb r > m \endSb e^{- 4 \pi v r} (r-m)^{\beta-1} r^{\alpha -1}   \, dr
&=m^s e^{- 4 \pi m v}
 \int_{0}^{\infty} e^{- 4 \pi m v r}
    r^{\beta-1} (r+1)^{\alpha -1} dr
\\
 &=m^s e^{- 4 \pi m v} \Gamma(\beta) \Psi_{-\ell}(s, 4 \pi m v).
\endalign
$$
This proves (i). The special value at $s=\ell-1$ is
$$
\align
W_{m, \infty}(\tau, \ell-1, \Phi_{\infty}^\ell)
&=2 \pi (-i)^\ell (2 \pi m)^{\ell-1} q^m \frac{\Psi_{-\ell}(\ell-1, 4 \pi m v)}{\Gamma(\ell)}
\\
 &=\frac{(- 2 \pi i m)^\ell}{m \Gamma (\ell)} q^m,
\endalign
$$
as claimed in (iv).

When  $m <0$, the integral is
$$
\align
\int \Sb r > 0 \endSb e^{- 4 \pi v r} (r-m)^{\beta-1} r^{\alpha -1}   \, dr
 &= |m|^s \int_{0}^{\infty} e^{- 4 \pi m v r}
    r^{\alpha -1} (r+1)^{\beta -1} dr
\\
 &=|m|^s  \Gamma(\alpha) \Psi_{\ell}(s, 4 \pi m v).
\endalign
$$
This proves (ii). The special value at $s=\ell-1$ is $0$ since
$\frac{1}{\Gamma(\beta)}=0$ at $s = \ell-1$ and $\Psi_\ell(\ell-1, 4 \pi |m| v)$ is
finite.
\qed\enddemo

{\bf{Proof of Lemma~8.9}.}  Since $m >0$, (i) of Proposition 15.1 implies
$$
\frac{W_{m, \infty}'(\tau, \ell-1, \Phi_\infty^\ell)}
     {W_{m, \infty}(\tau, \ell-1, \Phi_\infty^\ell)}
=\frac{1}2\log v + \log (2 \pi m) -\frac{1}2 \frac{\Gamma'(\ell)}{\Gamma(\ell)}
 +\frac{\Psi_{-\ell}'(\ell-1, 4 \pi m v)}{\Psi_{-\ell}(\ell-1, 4 \pi m v)}.
$$
Notice that, for any $z >0$,
$$
\Psi_{-\ell}(\ell-1, z)=\Psi(0, \ell; z) =1,
$$
by $(15.4)$. Observe that
$$
\Psi_{-\ell}(s, z)
=z^{-\beta} +
  \frac{1}{\Gamma(\beta)}\int_{0}^\infty e^{-z r}
   ((r+1)^{s-\beta}-1) r^{\beta -1} dr.
$$
The integral here is well-defined at $s=\ell-1$ and is equal to
$$
J(\ell-1, z):=\int_{0}^\infty e^{-z r}
   \frac{(r+1)^{\ell-1}-1}{ r}\, dr.
\tag{15.8}
$$
 Notice also that the function $\frac{1}{\Gamma(\beta)}$
vanishes  at $s=\ell-1$ and has the first derivative $\frac{1}2$ at $s=\ell-1$. Thus,
$$
\Psi_{-\ell}'(\ell-1, z)
  =-\frac{1}2 \log z + \frac{1}2\, J(\ell-1, z),
\tag{15.9}
$$
and we have
$$
\frac{W_{m, \infty}'(\tau, \ell-1, \Phi_\infty^\ell)}
     {W_{m, \infty}(\tau, \ell-1, \Phi_\infty^\ell)}
=\frac{1}2 \left[ \log( \pi m) -\frac{\Gamma'(\ell)}{\Gamma(\ell)}
           +J(\ell-1, 4 \pi m v)\right].
\tag{15.10}
$$
When $\ell=\frac{3}2$, $J(\frac12,4\pi m v) = J(4\pi m v)$ is the quantity defined in Theorem~8.8, so this gives Lemma~8.9.
\qed

{\bf{Proof of Lemma~8.11}.} (The case  $m< 0$, with derivative).
We now assume $m <0$. Since the function $\frac{1}{\Gamma(\beta)}$
vanishes  at $s=\ell-1$ and has the first derivative $\frac{1}2$ there,
one has, by (ii) of Proposition 15.1,
$$
W_{m, \infty}'(\tau, \ell-1, \Phi_{\infty}^\ell)
=2 \pi (-i)^\ell\, (2 \pi |m|)^{\ell-1}\, \frac{1}2\,
   \Psi_\ell(\ell-1, 4 \pi |m| v)\cdot e^{-4 \pi |m|v}\cdot q^m.
$$
By (15.6), one has
$$
\align
\Psi_\ell(\ell-1, 4 \pi |m| v)
 &=(4 \pi |m| v)^{1-\ell}\,\Psi_\ell(1-\ell, 4 \pi |m| v)
\\
 &=(4 \pi |m| v)^{1-\ell}\, \Psi(1, 2-\ell; 4 \pi |m| v)
\\
 &=(4 \pi |m| v)^{1-\ell}\,\int_{0}^\infty e^{-4 \pi |m| v r} (1+r)^{-\ell} \, dr
\\
 &=(4 \pi |m| v)^{1-\ell}\, e^{4 \pi |m| v}\,\int_{1}^\infty e^{-4 \pi |m| v r} r^{-\ell}
\, dr.
\endalign
$$
Therefore,
$$
W_{m, \infty}'(\tau, \ell-1, \Phi_{\infty}^\ell)
 =2 \pi\, (-i)^\ell\, 2^{-\ell}\, v^{1-\ell}\, q^m\, \int_{1}^\infty e^{-4 \pi |m| v r} r^{-\ell} \,
dr.
\tag{15.11}
$$
When $ \ell=\frac{3}2$, this gives Lemma~8.11. \qed

\subheading{\Sec16. The functional equation}

Let $D$ be a square-free positive integer, not necessarily the discriminant
of an indefinite quaternion algebra, and let
$$\Bbb E(\tau, s, \Phi^{\frac{3}2, D}) =c(D) \,(s+\frac12)\,\Lambda_D(2s+1)\,E(\tau, s, \Phi^{\frac{3}2, D})
\tag{16.1}
$$
be the renormalized Eisenstein series of (6.23). In this section, we prove
that it is invariant when $s$ goes to $-s$, i.e., that
$$
\Bbb E(\tau, s, \Phi^{\frac{3}2, D})
=\Bbb E(\tau, -s, \Phi^{\frac{3}2, D}).
\tag{16.2}
$$

First we need

\proclaim{Proposition 16.1}
 Set
$$
\Lambda(s, \chi_m; D)
 =\left(\,\frac{4 |m| D^2}\pi\,\right)^{\frac{1}2s}\, \Gamma(\frac{s+a}2)\, L(s, \chi_d)\,
  \prod_{p} b_p(n, s, D)
$$
with $a=(1+\sgn(m))/2$. Then $\Lambda(s, \chi_m, D)$ has a meromorphic
continuation to the whole complex $s$-plane with possible poles at $s=0$
and $1$, which occurs precisely when $D=-d=1$. Furthermore,  it satisfies the
following function equation
$$
\Lambda(s, \chi_m; D) =\Lambda(1-s, \chi_m; D),
$$
and
$$
\hbox{ord}_{s=0}\,\Lambda(s, \chi_m; D)
=\hbox{ord}_{s=1}\,\Lambda(s, \chi_m; D)
=\hbox{ord}_{s=1}\,L(s, \chi_d)+\#\{p|D: \chi_d(p)=1\}.
$$
\endproclaim
\demo{Proof} The functional equation follows from that of $L(s, \chi_d)$ and
$(8.10)$. The vanishing order at $s=0$ follows from
 $(8.14)$. We remark that $b_p(n, s;D)$ is a polynomial of
$p^{-s}$ even though it was written as a rational function and thus
is regular at $s =1/2$.\qed
\enddemo

\proclaim{Theorem 16.2} Let
$$
\Bbb E(\tau, s, \Phi^{\frac{3}2, D})
 =\sum_m \Cal A_m(v, s)\, q^m
$$
be the Fourier expansion of $\Bbb E(\tau, s, \Phi^{\frac{3}2, D})$.
\hfill\break
(i) For $m >0$, one has
$$
 \Cal A_m(v, s)
 =\frac{\Lambda(\frac{1}2+s, \chi_m; D)\,
   (4 \pi m v)^{\frac{s}2-\frac{1}4}\,\Psi_{-\frac{3}2}(s, 4 \pi m v)}
   {  \sqrt{\pi}\, \prod_{p|D} (1+p)}.
$$
(ii) For $m < 0$, one has
$$
\Cal A_m(v, s)
 =\frac{ (s^2 -\frac{1}4) \Lambda(\frac{1}2+s, \chi_m; D)\,
   (4 \pi |m| v)^{\frac{s}2-\frac{1}4}\,\Psi_{\frac{3}2}(s, 4 \pi | m| v)}
   {4 \sqrt{\pi} \prod_{p|D} (1+p)}\cdot e^{-4 \pi |m| v}.
$$
(iii) The constant term is
$$
\Cal A_0(v, s)=-\frac{D}{2 \pi \prod_{p|D}(p+1)} (G_D(s) +G_D(-s)),
$$
where
$$
G_D(s) =v^{-\frac{1}4 +\frac{s}2}\, \Lambda(1+2s)\, (s+\frac{1}2)\,
       \prod_{p|D}( p^{-\frac{1}2 -s} - p^{\frac{1}2 +s}).
$$
\endproclaim
\demo{Proof} When $m >0$, one has, by Propositions 8.1 and 15.1 and formula
(8.18),
$$
\align
\Cal A_m(v, s)
 &=c(D)  (\frac{D}\pi)^{s-\frac{1}2} \Gamma(s+\frac{3}2) \zeta_D(2s+1)
 \frac{C_{\infty}}{\sqrt 2} v^{\frac{s}2-\frac{1}4}
         (2 \pi m)^s
   \frac{\Psi_{-\frac{3}2}(s, 4 \pi m v)}{\Gamma(\frac{s}2+\frac{5}4)}\\
 &\qquad\qquad
  \times  C_f(D)\frac{L(s+\frac{1}2, \chi_d) \prod
b_p(n,s+\frac{1}2;D)}{\zeta_D(2s+1)}
\\
\nass
 &=\frac{\Lambda(s+\frac{1}2, \chi_m; D)}{  \prod_{p|D} (p+1)}
   \frac{\Gamma(s+\frac{3}2) 2^{-s-\frac{1}2}}
        {\Gamma(\frac{s}2+\frac{3}4)\Gamma(\frac{s}2+\frac{5}4)}
   (4 \pi m v)^{\frac{s}2-\frac{1}4} \Psi_{-\frac{3}2}(s, 4 \pi m v)
\endalign
$$
Now the doubling formula of the gamma function gives
$$
 \frac{\Gamma(s+\frac{3}2) }
        {\Gamma(\frac{s}2+\frac{3}4)\Gamma(\frac{s}2+\frac{5}4)}
 =2^{s+\frac{3}2 -1} \pi^{-\frac{1}2}
 =2^{s+\frac{1}2 } \pi^{-\frac{1}2}.
$$
This proves (i). The case $ m<0$ is the same and is left to the reader.
When $m=0$, one has by Corollary 8.2 and (8.12)
$$
\align
\Cal A_0(v, s)
 &=c(D) (\frac{D}\pi)^{s-\frac{1}2} \Gamma(s+\frac{3}2) \zeta_D(2s+1)
\\
 &\qquad\qquad \times
   \left[ v^{-\frac{1}4 +\frac{s}2}
  + \frac{C_\infty}{\sqrt 2}
    \frac{2^{-s}  v^{-\frac{1}4 -\frac{s}2} \Gamma(s)}
         {\Gamma(\frac{s}2-\frac{1}4)\Gamma(\frac{s}2+\frac{5}4)}
    \frac{C_f(D) \zeta(2s)}{\zeta_D(2s+1)}\prod_{p|D} (1-p^{1-2s})
   \right]
\\
\nass
\nass
 &=-\frac{  v^{-\frac{1}4 +\frac{s}2} \Lambda(1+2s) (\frac{1}2+s)}
        {2 \pi \prod_{p|D}(p+1)}
      (-1)^{\ord(D)} D^{s+\frac{3}2}\prod_{p|D}(1-p^{-1-2s})
 \\
 &\qquad\qquad -\frac{  v^{-\frac{1}4 -\frac{s}2} \Lambda(2s) (\frac{1}2-s)}
        {\sqrt{\pi} \prod_{p|D}(p+1)}
  \frac{\Gamma(s-\frac{1}2) 2^{\frac{1}2-s}}
       {\Gamma(\frac{s}2-\frac{1}4)\Gamma(\frac{s}2+\frac{1}4)}
       D^{\frac{1}2+s} \prod_{p|D} (1-p^{1-2s})
\\
\nass
\nass
 &=-\frac{D}{2 \pi \prod_{p|D}(p+1)}(G_D(s) + G_D(-s)).
\endalign
$$
Here we have used the doubling formula for the gamma functions again.
\qed
\enddemo

{\bf{Proof of the functional equation (16.2)}}. Now the functional equation
(16.2) follows immediately from Theorem~16.2, Proposition~16.1, and (15.6).
\qed

%
%

\vskip .5in

\redefine\vol{\oldvol}

\bye